\documentclass[a4paper]{article}

\usepackage[a4paper, total={6.5in, 8.5in}]{geometry}
\usepackage{amsmath, amssymb}   % the inevitable ams packages for beautiful math
\usepackage{mathtools}
\usepackage{graphicx,color}                   % graphicx package for graphics includes
\usepackage{float}  
\usepackage{xspace}         
\usepackage[scientific-notation=true]{siunitx}         

\usepackage{tikz}
\tikzstyle{startstop} = [rectangle, rounded corners,minimum width=3cm, minimum height=1cm, text centered, text width=5cm, draw=black]
\tikzstyle{process} = [rectangle, minimum width=3cm, minimum height=1cm, text width=6 cm,text centered, draw=black]
\tikzstyle{decision} = [diamond, minimum width=3cm, minimum height=1cm, text width=3cm, text centered, draw=black]
\tikzstyle{arrow} = [thick,->,>=stealth]
\usetikzlibrary{shapes.geometric, arrows}
\usetikzlibrary{positioning} 
\usetikzlibrary{calc}
\usepackage{pgfplots}
\pgfplotsset{width=10cm,compat=1.9}

\usepackage[justification=centering]{caption}

\usepackage{hyperref}
\usepackage[all]{hypcap}

\newcommand{\ep}{\varepsilon}    % ...and eps much better than the predefined.
\newcommand{\bb}{{\boldsymbol b}}

\newcommand{\me}{E\in \mathcal{E}_h}
\newcommand{\bn}{\boldsymbol n} 
\newcommand{\bt}{\boldsymbol t} 
\newcommand{\bx}{\boldsymbol x} 
\newcommand{\cF}{{\mathcal F}}

\newcommand{\afce}{{{AFC-energy}}\xspace}
\newcommand{\afcse}{{{AFC-SUPG-energy}}\xspace}

\definecolor{light_gray}{gray}{0.75}
\colorlet{light_blue}{blue!20}

\newcommand{\blist}{\begin{list}{}{\itemsep0.0ex\parsep0.1ex\topsep0.2ex\leftmargin1.6em\labelwidth1.3em}}

\newtheorem{lemma}{Lemma}
\newtheorem{theorem}[lemma]{Theorem}
\newtheorem{remark}[lemma]{Remark}
\newtheorem{definition}[lemma]{Definition}
\newtheorem{example}[lemma]{Example}
\definecolor{light_gray}{gray}{0.75}
\definecolor{lighter_gray}{gray}{0.5}
\colorlet{light_blue}{blue!20}
\definecolor{dark_green}{rgb}{0.0, 0.6, 0.0}
\definecolor{royal_blue}{rgb}{0.0, 0.22, 0.66}
\definecolor{salmon}{rgb}{1.0, 0.55, 0.41}
\definecolor{gold}{rgb}{0.8, 0.63, 0.21}
\definecolor{navy_blue}{rgb}{0.0, 0.0, 0.5}
\definecolor{crimson}{rgb}{0.79, 0.0, 0.09}
\definecolor{amethyst}{rgb}{0.6, 0.4, 0.8}
\definecolor{alizarin}{rgb}{0.82, 0.1, 0.26}
\definecolor{amaranth}{rgb}{0.9, 0.17, 0.31}
\definecolor{azure}{rgb}{0.0, 0.5, 1.0}
\definecolor{canaryyellow}{rgb}{0.82, 0.41, 0.12}
\definecolor{carrotorange}{rgb}{0.8, 0.33, 0.0}
\definecolor{cadmiumgreen}{rgb}{0.0, 0.42, 0.24}
\definecolor{copper}{rgb}{0.72, 0.45, 0.2}
\definecolor{aqua}{rgb}{0.5, 1.0, 0.83}
\definecolor{awesome}{rgb}{1.0, 0.13, 0.32}
\definecolor{candyapplered}{rgb}{1.0, 0.03, 0.0}
\definecolor{caribbeangreen}{rgb}{0.0, 0.8, 0.6}

\title{Adaptive Grids in the Context of Algebraic Stabilizations for Convection-Diffusion-Reaction Equations}
\author{Abhinav Jha
\footnote{RWTH Aachen University, Applied and Computational Mathematics, Schinkelstra\ss e 2, 52062, Aachen, Germany, Email: \texttt{jha@acom.rwth-aachen.de}},
Volker John
\footnote{Weierstrass Institute for Applied Analysis and Stochastics (WIAS), Mohrenstr.
39, 10117 Berlin, Germany and Freie Universit\"at Berlin, Department of Mathematics and
Computer Science, Arnimallee 6, 14195 Berlin, Germany,  Email: \texttt{john@wias-berlin.de}},
Petr Knobloch
\footnote{Department of Numerical Mathematics, Faculty of Mathematics
and Physics, Charles University, Sokolovsk\'a 83, Praha 8, 18675,
Czech Republic, Email: \texttt{knobloch@karlin.mff.cuni.cz}}}
\date{}

\begin{document}
\maketitle
\begin{abstract}
Three algebraically stabilized finite element schemes for discretizing convection-diffusion-reaction equations are 
studied on adaptively refined grids. These schemes are the algebraic flux correction (AFC) scheme 
with Kuzmin limiter, the AFC scheme with BJK limiter, and the recently proposed 
Monotone Upwind-type Algebraically Stabilized (MUAS) method. Both, conforming closure of the refined grids 
and grids with hanging vertices are considered. A non-standard algorithmic step becomes necessary before these 
schemes can be applied on grids with hanging vertices. 
The assessment of the schemes is performed with respect 
to the satisfaction of the global discrete maximum principle (DMP), the accuracy, e.g., smearing of layers, and the efficiency in solving the corresponding nonlinear problems. 
\end{abstract}


\textbf{AMS :}
65N12, 65N30

\section{Introduction}
The physical behavior of scalar quantities, like temperature (energy) or concentrations, in fluids is modeled 
by scalar convection-diffusion-reaction equations. Let $\Omega \subset \mathbb R^d$, $d\in\{2,3\}$, be
a bounded domain with a Lipschitz-continuous boundary $\partial\Omega$. In this paper, we consider the steady-state equations, 
which are given, already in non-dimensional form, as follows:
\begin{eqnarray}\label{eq:cdr_eqn}
-\varepsilon \Delta u+\bb\cdot \nabla u+cu& = & f\quad \mathrm{in}\ \Omega\nonumber\\
u&=&u_b\quad \mathrm{on}\ \Gamma_D,\\
\varepsilon \nabla u \cdot \bn & = & g\quad \mathrm{on}\ \ \Gamma_N.\nonumber
\end{eqnarray}
Here, $\varepsilon>0$ is the diffusion coefficient, $\bb\in
W^{1,\infty}(\Omega)^d$ is the convective field, $c\in
L^{\infty}(\Omega)$ is the reaction field, $f\in L^2(\Omega)$ is the
source or the sink term, $u_b\in H^{1/2}(\Gamma_D)$ and $g\in
L^2(\Gamma_N)$ specify the boundary conditions, $\bn$ is the unit outward normal to $\partial\Omega$,
$\Gamma_D\cup \Gamma_N=\partial \Omega$, $\Gamma_D\cap \Gamma_N=\emptyset$,
$\Gamma_D$ is the Dirichlet boundary and $\Gamma_N$ is the Neumann boundary. 
%It is assumed that $-(\nabla\cdot\bb(\bx))/2 + c(\bx) \ge \sigma_0 > 0$ in $\Omega$.
Under appropriate assumptions on the data, it is well known that problem \eqref{eq:cdr_eqn} possesses a 
unique weak solution. 

In practice, the convective transport usually dominates the diffusive transport. One speaks of the 
convection-dominated regime, given if $\varepsilon \ll L \|\bb\|_{L^\infty(\Omega)}$, where $L$ is a 
characteristic length scale of the problem. Then, a characteristic feature of solutions of  \eqref{eq:cdr_eqn} 
are layers, which are thin regions with a steep gradient. 
In general, computational grids cannot resolve layers. It is well known that one has to apply so-called 
stabilized discretizations in this situation, e.g., see \cite{RST08}. There are many proposals of 
such discretizations for convection-diffusion-reaction equations in the literature. 

For appropriate data, the solution of  \eqref{eq:cdr_eqn} takes only certain physical values, e.g., concentrations
are non-negative. The mathematical formulation of this physical feature is called maximum principle,  
see \cite{GT83}. For numerical simulations in practice, it is often of utmost importance that also the 
discrete solution possesses only physically consistent values, i.e., it satisfies a discrete maximum
principle (DMP). However, there are only very few among the stabilized discretizations with this property. 
The currently most promising class of methods seems to be the class of algebraically stabilized schemes. 

Algebraically stabilized discretizations have been becoming quite popular for a couple of years. Their 
construction relies (mainly) on the algebraic system of equations from the Galerkin finite element 
discretization of \eqref{eq:cdr_eqn} with conforming piecewise linear finite elements. Then, an 
algebraic stabilization term is introduced and certain coefficients (limiters) are computed that depend 
on the concrete discrete solution. Hence, these methods are nonlinear. The first comprehensive numerical 
analysis for the so-called algebraic flux correction (AFC) scheme with Kuzmin limiter, from \cite{Ku07},
was presented in \cite{BJK16}. The AFC scheme with BJK limiter was proposed and analyzed in \cite{BJK17}.
Recently, a new algebraically stabilized method was proposed in \cite{JK21}, which is called 
Monotone Upwind-type Algebraically Stabilized (MUAS) method. For all these methods, DMPs could be proved, 
sometimes under appropriate assumptions. For a detailed presentation of the methods and a discussion 
of the DMPs, we refer to Section~\ref{sec:as_schemes}.

Because of the presence of layers, where discrete solutions usually possess large errors, it is very 
attractive to use adaptively refined grids for the numerical solution of convection-diffusion-reaction 
equations. The control of adaptive grid refinement relies on an a posteriori error estimator or indicator. 
The first error estimator for the AFC schemes with Kuzmin and with BJK limiter has been proposed 
recently in \cite{Jha20}. On the basis of the error estimator or indicator, certain mesh cells are 
marked for refinement. 
Some of the common strategies to refine a grid can be found in \cite{Riv84, BSW83, KR89}. 
The first step of refining a grid, i.e., the refinement of the marked cells, leads to the formation of hanging vertices. 
In the framework of discontinuous finite elements, the handling of grids with hanging vertices is rather easy to understand, see \cite{AR10}. 
For continuous finite elements, the framework becomes more involved. 
A commonly used way around this issue is to use conforming closure or red-green refinements, 
see \cite{BSW83}, but this approach leads to the deterioration of angles. 
Also, while using hexahedral mesh cells in 3d, the green completion leads to the formation of pyramids or prisms, which are not easy to handle by many finite element codes. Hence, using grids with hanging vertices 
is attractive from the geometric point of view, because one can perform a simple grid refinement.

This paper explores the behavior of the three above-mentioned algebraically stabilized methods in simulations 
on adaptive grids in two dimensions. An initial comparison of the AFC schemes was performed in \cite{Jha20}, with the 
emphasis on studying the performance of two a posteriori error estimators in terms of their effectivity
indices and their control of the adaptive refinement process. Concerning the MUAS method, 
some of its properties are illustrated numerically with simulations on uniform 
grids in \cite{JK21}. In the current paper, first studies of this method on adaptively refined grids will be presented.
The goal of the numerical studies consists in comparing the methods with respect to accuracy, to the satisfaction
of the global DMP, and to efficiency in solving the nonlinear problems. A particular attention will be 
paid to the study of algebraic stabilizations on grids with hanging vertices. To the best of our 
knowledge, there is no such study in the literature so far.  We could find the use of an algebraic stabilization 
on grids with hanging nodes only for the linear transport equation in \cite{BK13}.
From the algorithmic point of view, 
it will be shown that 
compared with the standard approach of modifying a linear system of equations for discretizations on grids 
with hanging nodes, an additional step becomes necessary for algebraically stabilized schemes, namely the transform to conforming ansatz functions. Such a step is not reported in \cite{BK13}.

\iffalse
Throughout this paper, we use standard notations for Sobolev spaces and their
norms, see \cite{Ada75}.  The inner product in $L^2(\Omega)$ is denoted by
$\left( \cdot, \cdot\right)$. The norm (semi-norm) on $W^{m,p}(\Omega)$ is
denoted by $\|\cdot\|_{m,p,\Omega}$ ($|\cdot|_{m,p,\Omega}$), with the
convention $\|\cdot\|_{m,\Omega}=\|\cdot\|_{m,2,\Omega}$ and $\|\cdot \|_m$
($|\cdot |_m$) when the considered domain is clear.
\fi

The paper is organized as follows. Section~\ref{sec:prelim} introduces concepts of triangulations, in 
particular with hanging vertices, and corresponding finite element spaces. The algebraically stabilized 
methods are described in Section~\ref{sec:as_schemes}. Some information concerning the implementation of 
these discretizations on grids with hanging vertices are provided in Section~\ref{sec:hang_afc}. 
The numerical studies are presented in Section~\ref{sec:numres}. Section~\ref{sec:summary} summarizes
the findings of this paper.

\section{Triangulations and Finite Element Spaces}\label{sec:prelim}
This section introduces notations and recalls concepts with respect to 
triangulations and finite element spaces. A special emphasis is paid to 
triangulations with hanging vertices. Some work concerning this topic can be 
found in \cite{Car11}, where results are provided for the lowest order Lagrange 
elements in the framework of multigrid methods. Recently, in \cite{Jha21} the theory has been extended for higher order Lagrange elements. Most of the definitions in this
section follow standard texts, e.\,g., see \cite{Cia78, BS08}.

Let $\Omega \subset \mathbb{R}^d$, $d\in\{2,3\}$, be a polygonal
resp.~polyhedral domain that is decomposed into simplices (i.e., triangles
resp.~tetrahedra). This decomposition is referred to as triangulation
and is denoted by $\mathcal{T}_h$. As usual, it is 
assumed that the interiors of any two different elements of $\mathcal{T}_h$ are 
disjoint and that $\overline{\Omega}=\cup_{K\in \mathcal{T}_h}K$. A 
triangulation $\mathcal{T}_h$ of $\Omega$ is called conforming if, for any
$K_1,K_2\in \mathcal{T}_h$ with $K_1\neq K_2$, the intersection $K_1\cap K_2$ 
is either empty or a vertex or an edge or, in 3d, a face of both $K_1$ and 
$K_2$.  It is assumed that any edge or face lying on $\partial\Omega$ is
a subset of either $\overline{\Gamma_D}$ or $\overline{\Gamma_N}$.

For a given triangulation $\mathcal{T}_h$, we denote by $\mathcal{N}_h$ the set
of all vertices, by $\mathcal{E}_h$ the set of all edges, and by
$\mathcal{F}_h$ the set of all facets (i.e., all edges resp.~faces). Thus, 
in 2d, it holds that $\mathcal E_h = \cF_h$. The set of facets can be 
decomposed into $\mathcal{F}_h=\mathcal{F}_{h,\Omega}\cup\mathcal{F}_{h,D}\cup 
\mathcal{F}_{h,N}$, where $\mathcal{F}_{h,\Omega}, \mathcal{F}_{h,D}$, and 
$\mathcal{F}_{h,N}$ are the interior, Dirichlet, and Neumann facets, 
respectively. We denote the diameter of a mesh cell $K$ by $h_K$ and 
the diameter of an edge $E$ and a facet $F$ by $h_E$ and $h_F$, respectively.

\begin{definition}[Refinement,  \cite{Car11},  Def.~3.3]
Let $\mathcal{T}_{1}$ and $\mathcal{T}_{2}$ be triangulations of $\Omega$.
Then, $\mathcal{T}_{2}$ is called a refinement of $\mathcal{T}_{1}$ if for
all $K\in \mathcal{T}_{1}$ the set
$
\lbrace K'\in \mathcal{T}_{2}:K'\cap\mbox{\rm int}\,K\neq \emptyset\rbrace
$
is a triangulation of $K$, where $\mbox{\rm int}\,K$ is the interior of $K$.
\end{definition}

\begin{definition}[Grid hierarchy, \cite{Car11}, Def.~3.4] 
A family $\lbrace \mathcal{T}_i\rbrace_{i=0}^j$ is called a grid hierarchy on 
$\Omega$ if $\mathcal{T}_0$ is a conforming triangulation of $\Omega$ and if 
each $\mathcal{T}_i,i=1,\dots,j$, is a refinement of $\mathcal{T}_{i-1}$. 
%If the grid $\mathcal{T}_i$ is conforming, we call it conforming grid 
%refinement otherwise non-conforming grid refinement. 
% -- THESE NOTIONS ARE NOT USED IN THE PAPER
\end{definition}

\begin{definition}[Hanging vertex, \cite{Car11}, Def.~3.6] 
Let $\mathcal{T}_h$ be a triangulation of $\Omega$. Then, a vertex 
$p\in\mathcal{N}_h$ is called a hanging vertex if there is an element 
$K\in\mathcal{T}_h$ with $p\in\partial K$ but $p$ is not a vertex of $K$.
The set of all hanging vertices is denoted by $\mathcal{H}_h$.
\end{definition}

In this work, we will consider first order Lagrange finite element spaces
  $$
  S(\mathcal{T}_h):=\left\lbrace v\in \mathcal{C}(\overline{\Omega}):v|_K\in \mathbb{P}_1(K)\ \ \forall\ K\in\mathcal{T}_h\right\rbrace 
  $$
consisting of continuous functions on $\overline{\Omega}$ such that the 
restrictions to all cells $K\in \mathcal{T}_h$ are polynomials of degree at 
most $1$. It is well known that $S(\mathcal{T}_h)\subset H^1(\Omega)$. Degrees
of freedom which determine functions from $S(\mathcal{T}_h)$ are values at
vertices. Therefore, vertices are also called nodes. Due to the continuity
requirement, values at hanging nodes depend on the values at non-hanging nodes
as it is stated in the following lemma.

\begin{lemma}\emph{(\cite[Lemma~3.2]{Car11})} 
\label{lem:hang_comb}Let $\left\lbrace\mathcal{T}_0, \cdots, \mathcal{T}_j\right\rbrace$ be a grid 
hierarchy on $\Omega$. Let us denote $\mathcal{T}_h=\mathcal{T}_j$, i.e., the 
final refinement level. Then, for all $q\in\mathcal{H}_h$ there are 
coefficients $a_{qp}$ with $p\in\mathcal{N}_h\setminus\mathcal{H}_h$ such that 
all $v\in S(\mathcal{T}_h)$ can be represented as
\begin{equation}\label{eq:hanging_node_representation}
  v(q)=\sum_{p\in\mathcal{N}_h\setminus\mathcal{H}_h}a_{qp}v(p).
\end{equation}
\end{lemma}

For conforming triangulations, a basis of $S(\mathcal{T}_h)$ is given by the 
well-known nodal basis functions. To construct basis functions of
$S(\mathcal{T}_h)$ for a non-conforming triangulation, we first introduce 
non-conforming nodal basis functions that are generally not in
$S(\mathcal{T}_h)$.

\begin{definition}[Non-conforming nodal basis functions] 
\label{def:nc_basis_fct}Let 
$\mathcal{T}_h$ be a triangulation of $\Omega$. Then, the non-conforming nodal 
basis function $\varphi_p^{\mathrm{nc}}\in L^2(\Omega)$ associated with 
$p\in \mathcal{N}_h$ is defined as follows: For all $K\in \mathcal{T}_h$ there 
is a representative $\varphi_p^{\mathrm{nc}}|_K=\mu_{p,K}\in\mathbb{P}_1(K)$ 
with $\mu_{p,K}(q)=\delta_{pq}$ for all vertices $q$ of $K$.
\end{definition}

For a conforming mesh $\mathcal{T}_h$ this definition reduces to 
$\varphi_p^{\mathrm{nc}}\in S(\mathcal{T}_h)$ and
$
\varphi_p^{\mathrm{nc}}(q)=\delta_{pq}$ for all 
$p,q\in \mathcal{N}_h$, 
i.e., the set $\lbrace\varphi_p^{\mathrm{nc}}\rbrace_{p\in\mathcal{N}_h}$ is 
the conforming nodal basis of $S(\mathcal{T}_h)$. For a non-conforming 
triangulation, $S(\mathcal{T}_h)$ is in general only a subspace of the 
non-con\-forming finite element space
$$
S^{\mathrm{nc}}(\mathcal{T}_h):=\mathrm{span}\left\lbrace \varphi_p^{\mathrm{nc}}:p\in \mathcal{N}_h\right\rbrace.
$$
However, it is possible to construct a basis of $S(\mathcal{T}_h)$ from the 
non-conforming nodal basis of $S^{\mathrm{nc}}(\mathcal{T}_h)$.
%that resembles the usual nodal basis functions when $\mathcal{T}_h$ is
%conforming. -- THIS REPEATS WHAT WAS MENTIONED ALREADY ABOVE

\begin{theorem}\label{theorem:basis} \emph{(\cite[Theorem~3.1]{Car11})}
Let $\left\lbrace\mathcal{T}_0, \cdots, \mathcal{T}_j\right\rbrace$ be a grid 
hierarchy on $\Omega$. Let us denote $\mathcal{T}_h=\mathcal{T}_j$, i.e., the 
final refinement level. Then, a basis of $S(\mathcal{T}_h)$ is given by
  \begin{equation*}
  \left\lbrace \varphi_p=\varphi_p^{\mathrm{nc}}+\sum_{q\in\mathcal{H}_h}a_{qp}\varphi_q^{\mathrm{nc}}:p\in\mathcal{N}_h\setminus\mathcal{H}_h\right\rbrace,
  \end{equation*}
where the coefficients $a_{qp}$ are the same as in Lemma~\ref{lem:hang_comb}.
\end{theorem}

\section{Algebraically Stabilized Schemes}\label{sec:as_schemes}

As already mentioned, algebraic stabilizations are currently the most promising finite element discretizations for computing 
numerical solutions of steady-state convection-diffusion-reaction equations that satisfy DMPs.
This section presents the methods that will be studied. 

%Consider $\mathbb P_1$ finite elements. -- THIS WAS STATED IN THE PRECEDING
%SECTION
The first step of algebraically stabilized schemes consists in applying the 
standard Galerkin finite 
element method to the weak form of \eqref{eq:cdr_eqn}. Then, the discrete solution can be represented as a vector $U\in \mathbb{R}^N$, with the last $N-M$ components corresponding to the Dirichlet boundary conditions. The algebraic representation of the
method is given by
\begin{equation*}
AU=b,
\end{equation*}
where $A=(a_{ij})_{i,j=1}^N$ is the corresponding stiffness matrix and $b\in \mathbb{R}^N$ is the assembled right-hand side. In an algebraically stabilized method, an 
additional nonlinear stabilization term is added such that it takes the form 
\begin{equation}\label{eq:matrix_rep_asm}
\left(A +B(U)\right)U=b,
\end{equation}
with $B(U)=( b_{ij}(U))_{i,j=1}^N$.
For preserving conservation of the discrete solution, the stabilization has to 
be symmetric: $b_{ij}(U) = b_{ji}(U)$, $i,j=1,\ldots,M$. 

\subsection{AFC Scheme with Kuzmin Limiter}

AFC schemes consider in the first step the Galerkin finite element discretization 
in the case 
that Neumann boundary conditions are applied, i.e., it is $M=N$. The stabilization term in 
\eqref{eq:matrix_rep_asm} is of the form 
\begin{equation}\label{eq:stab_afc}
b_{ij}(U)=(1-\alpha_{ij}(U)) d_{ij}\quad \forall\ i\neq j,\qquad 
b_{ii}(U)=-\sum_{j\neq i}b_{ij}(U),
\end{equation}
where $D=(d_{ij})_{i,j=1}^N$ is an artificial diffusion matrix with entries 
\begin{equation}\label{eq:mat_art_diff}
d_{ij}=d_{ji}=-\max\lbrace a_{ij},0,a_{ji}\rbrace\quad \forall\ i\neq j,\qquad d_{ii}=-\sum_{j\neq i}d_{ij},
\end{equation}
and $(\alpha_{ij}(U))_{i,j=1}^N$ is the limiter matrix with $0\le \alpha_{ij}(U)\le 1$.
After having computed the limiters, Dirichlet boundary conditions are imposed in the usual way. 

In subregions where no layers appear, the standard Galerkin discretization can be applied. 
In this case, the corresponding limiters should be close to $1$. In a vicinity of layers, a 
stabilization is necessary, which is achieved by using values of the limiter that are much 
smaller than $1$. 

The Kuzmin limiter, proposed in \cite{Ku07}, is a monolithic upwind-type limiter and it is applicable to 
$\mathbb{P}_1$ and $\mathbb{Q}_1$ elements. For $\mathbb P_1$ elements, the existence  of a solution
is proved in \cite{BJK16}.  For a real number $a$, denote $a^+=\max\{a,0\}$ and 
$a^-=\min\{a,0\}$. Then, the limiters are computed as follows:
\begin{enumerate}
    \item Compute 
    \begin{equation*}
       P_i^+  =  \sum_{j=1,a_{ji}\leq a_{ij}}^N \left(
d_{ij}(u_j-u_i)\right)^+,\qquad 
       P_i^-  =  \sum_{j=1,a_{ji}\leq a_{ij}}^N \left(
d_{ij}(u_j-u_i)\right)^-.
    \end{equation*}
    \item Compute 
       \begin{equation*}
       Q_i^+ = -\sum_{j=1}^N \left( d_{ij}(u_j-u_i)\right)^-, \qquad 
       Q_i^- = -\sum_{j=1}^N \left( d_{ij}(u_j-u_i)\right)^+.
    \end{equation*}
    \item Compute
    \begin{equation*}
       R_i^+ = \min\left\{ 1,\frac{Q_i^+}{P_i^+} \right\}, \quad R_i^- =
\min\left\{ 1,\frac{Q_i^-}{P_i^-} \right\},\qquad i=1,\dots, M.
    \end{equation*}
     If \(P_i^+\) or \(P_i^-\) is zero, one sets \(R_i^+=1\) or \(R_i^-=1\), respectively.  The values of \(R_i^+\) and \(R_i^-\) are set to $1$ for Dirichlet nodes as well. 
    \item If $a_{ji} \le a_{ij}$, then set
    \[ \alpha_{ij} = \begin{cases} 
        R_i^+ & \mbox{ if } d_{ij}(u_j-u_i)>0,\\
        1 & \mbox{ if } d_{ij}(u_j-u_i)=0,\\
        R_i^- & \mbox{ if } d_{ij}(u_j-u_i)<0,
      \end{cases} \quad \alpha_{ji} := \alpha_{ij},
   \]
\end{enumerate}
for $i,j=1,\dots,N$. Note that the symmetry of the stabilization term follows 
from the symmetries of $D$ and the limiters. For the Kuzmin limiter, the local 
DMP is satisfied if the off-diagonal entries of $A$ possess a certain property, 
see \cite{Kno17} for details. It is also shown in this paper that this property 
and also the local DMP may be violated for certain types of 
triangulations, e.g., in two dimensions if the triangulation is not of 
Delaunay type

\subsection{AFC Scheme with BJK Limiter}

This method, proposed in \cite{BJK17}, starts in the same way as the previous method and the 
stabilization term has the form \eqref{eq:stab_afc}. 
It was derived for $\mathbb{P}_1$ elements. For this method, the existence of a solution of
the nonlinear problem and the satisfaction of a local and global DMP on arbitrary
conforming simplicial grids
can be proved. Moreover, it was shown in \cite{BJK17} that it is linearity
preserving, i.e., the stabilization term vanishes for any vector that represents a linear function.

The computation of the limiter starts with a pre-processing step, compare \cite[Eq.~(2.4)]{BJK17}. Then, the computation proceeds as follows:
\begin{enumerate}
    \item Compute 
    \begin{equation*}
       P_i^+  = \sum_{j=1}^N \left( d_{ij}(u_j-u_i)\right)^+,\qquad
       P_i^-  = \sum_{j=1}^N \left( d_{ij}(u_j-u_i)\right)^-.
    \end{equation*}
    \item\label{it:BJK2} Compute 
       \begin{equation*}
       Q_i^+ = q_i\left(u_i-u_i^{\max}\right), \qquad Q_i^- = q_i\left(u_i-u_i^{\min}\right),
    \end{equation*}
    with
    $$
    u_i^{\max}  =  \max_{j\in N_i\cup\{i\}}u_j,\qquad
    u_i^{\min}   = \min_{j\in N_i\cup\{i\}}u_j,\qquad
    q_i  = \sum_{j\in N_i}\gamma_i d_{ij},
    $$
    where 
$
   N_i=\{j\in\{1,\dots,N\}\setminus\{i\}\,:\,\,
             a_{ij}\neq0\,\,\,\mbox{or}\,\,\,a_{ji}>0\}
$
    and $\gamma_i$ is a positive constant 
%computed for interior nodes as given in \cite[Rem.~6.2]{BJK17}, 
% -- THE COMPUTATIONS DIFFERS FROM THE PAPER
which guanrantees the linearity preservation, see 
Section~\ref{sec:hang_afc} for details. 

    \item Compute
    \begin{equation*}
       R_i^+ = \min\left\{ 1,\frac{Q_i^+}{P_i^+} \right\}, \quad R_i^- =
\min\left\{ 1,\frac{Q_i^-}{P_i^-} \right\},\qquad i=1,\dots, M.
    \end{equation*}
If \(P_i^+\) or \(P_i^-\) is zero, one sets \(R_i^+=1\) or \(R_i^-=1\), respectively. The values for \(R_i^+\) and \(R_i^-\) are set to $1$ also for Dirichlet nodes. 
    \item Compute
    \[ \overline{\alpha}_{ij} = \begin{cases} 
        R_i^+ & \mbox{ if } d_{ij}(u_j-u_i)>0,\\
        1  & \mbox{ if } d_{ij}(u_j-u_i)=0, \\
        R_i^- & \mbox{ if } d_{ij}(u_j-u_i)<0,
      \end{cases}\quad i,j=1,\dots, N.
   \]
\end{enumerate}
Finally, one sets
\begin{equation*}
\alpha_{ij} = \min \left\lbrace\overline{\alpha}_{ij}, \overline{\alpha}_{ji}\right\rbrace,\quad i,j=1,\dots,N.
\end{equation*}
Again, the symmetry of the stabilization term follows from the symmetries of $D$ and of 
the limiters. 

\subsection{Monotone Upwind-type Algebraically Stabilized (MUAS) Meth\-od}

The MUAS method was recently proposed and analyzed in \cite{JK21}, where 
the solvability of the nonlinear discrete problem and the satisfaction of local and global DMPs 
on arbitrary conforming simplicial grids are proved. 

Also in this method, the matrix obtained for Neumann boundary conditions is considered in the 
first step. 
The stabilization term in \eqref{eq:matrix_rep_asm} is given by 
\begin{eqnarray*}
   b_{ij}(U)
   &= & -\max\{(1-\alpha_{ij}(U)) a_{ij},0, (1-\alpha_{ji}(U)) a_{ji}\},\quad
   i,j=1,\dots,N, \ i\neq j,\\
   b_{ii}(U) &=& -\sum_{j=1,j\neq i}^N\,b_{ij}(U),\quad i=1,\dots,N,
\end{eqnarray*}
which is clearly symmetric. The limiters $\alpha_{ij}(U)$ are computed as 
follows:
\begin{enumerate}
    \item Compute 
    \begin{equation*}
       P_i^+  =  \sum_{j=1, a_{ij}>0}^N a_{ij}(u_i-u_j)^+,  \quad
       P_i^-  =  \sum_{j=1, a_{ij}>0}^N  a_{ij}(u_i-u_j)^-.
    \end{equation*}
    \item Compute 
       \begin{equation*}
       Q_i^+ =\sum_{j=1}^N \max\left\{|a_{ij}|,a_{ji}\right\}(u_j-u_i)^+, \quad 
       Q_i^-=\sum_{j=1}^N \max\left\{|a_{ij}|,a_{ji}\right\}(u_j-u_i)^-.
    \end{equation*}
    \item Compute
    \begin{equation*}
       R_i^+ = \min\left\{ 1,\frac{Q_i^+}{P_i^+} \right\}, \quad R_i^- =
\min\left\{ 1,\frac{Q_i^-}{P_i^-} \right\},\qquad i=1,\dots, M.
    \end{equation*}
     If \(P_i^+\) or \(P_i^-\) is zero, one sets \(R_i^+=1\) or \(R_i^-=1\), respectively.  The values of \(R_i^+\) and \(R_i^-\) are set to $1$ for Dirichlet nodes as well. 
    \item Define
    \[ \alpha_{ij} = \begin{cases} 
        R_i^+ & \mbox{ if } u_i>u_j ,\\
        1 & \mbox{ if } u_i=u_j,\\
        R_i^- & \mbox{ if } u_i<u_j,
      \end{cases} \quad  i,j=1,\dots,N.
   \]
\end{enumerate}

\section{Hanging Nodes and Algebraically Stabilized Schemes}\label{sec:hang_afc}
For discret\-iza\-tions on grids with hanging nodes, first a linear system of equations for the non-conforming basis functions $\varphi_p^{\mathrm{nc}}$, 
see Definition~\ref{def:nc_basis_fct}, is assembled. The next step consists in transforming this 
system to a system corresponding to conforming test functions $\varphi_p$, introduced in Theorem~\ref{theorem:basis}.
Constraints are set for the values at the hanging nodes such that the finite element solution becomes continuous. 
An example will illustrate this approach.

\begin{example}[System corresponding to non-conforming ansatz and conforming test functions]
\label{ex:ImplementationP1}Consider a patch as defined in Figure~\ref{fig:row_sum_hanging}. The 
non-conforming space~$S^{\mathrm{nc}}(\mathcal{T}_h)$ is spanned from the following 
basis functions:
\begin{align*}
\varphi_{i_0}^{\mathrm{nc}}(x,y) & = \begin{cases}
                                      0 & \text{in } K_1,\\
                                      2-2y & \text{in } K_2,\\
                                      2x & \text{in } K_3,
                                     \end{cases}
\quad &
\varphi_{i_1}^{\mathrm{nc}}(x,y) & = \begin{cases}
                                      1-x & \text{in } K_1,\\
                                      0 & \text{in } K_2,\\
                                      1-x-y & \text{in } K_3,
                                     \end{cases}
\\
\varphi_{i_2}^{\mathrm{nc}}(x,y) & = \begin{cases}
                                      x-y & \text{in } K_1,\\
                                      0 & \text{in } K_2,\\
                                      0 & \text{in } K_3,
                                     \end{cases}
\quad &
\varphi_{i_3}^{\mathrm{nc}}(x,y) & = \begin{cases}
                                      y & \text{in } K_1,\\
                                      x+y-1 & \text{in } K_2,\\
                                      0 & \text{in } K_3,
                                     \end{cases}
\\
\varphi_{i_4}^{\mathrm{nc}}(x,y) & = \begin{cases}
                                      0 &\text{in } K_1,\\
                                      -x+y & \text{in } K_2,\\
                                      -x+y & \text{in } K_3.
                                     \end{cases}
\end{align*}
The conforming space is 
$S(\mathcal{T}_h) = \mathop{\mathrm{span}}\{\varphi_j \,|\, j\in\{i_1, \dots, i_4\} \}$,
where the continuous basis functions are given by 
$\varphi_{i_j}=\varphi_{i_j}^{\mathrm{nc}}$ for $j\in\{2,4\}$ and
\begin{equation*}
\varphi_{i_1} 
 = \varphi_{i_1}^{\mathrm{nc}} + \frac{1}{2}\varphi_{i_0}^{\mathrm{nc}}
 = \begin{cases}
    1-x & \text{in } K_1,\\
    1-y & \text{in } K_2,\\
    1-y & \text{in } K_3,
\end{cases}
\quad
\varphi_{i_3}
 = \varphi_{i_3}^{\mathrm{nc}} + \frac{1}{2}\varphi_{i_0}^{\mathrm{nc}}
 = \begin{cases}
    y & \text{in } K_1,\\
    x & \text{in } K_2,\\
    x & \text{in } K_3.
\end{cases}
\end{equation*}
This means, the coefficients $a_{qp}$ from Lemma~\ref{lem:hang_comb} (with
$q=0$), given by
$a_{qp} = \varphi_{i_p}^{\mathrm{nc}}\big|_{K_1}$ evaluated at $(0.5,0.5)$,
are zero for $p\in\{2,4\}$ and $1/2$ for $p\in\{1,3\}$.
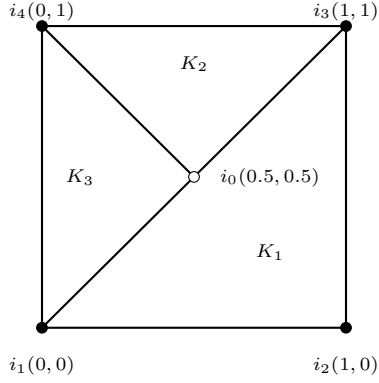
\begin{figure}[t!]
\centering
\begin{tikzpicture}
\fill[line width=0.8pt, color = white] (0,0) -- (4,0) -- (4,4) -- (0,4) -- cycle;
\draw [line width=0.8pt] (0,0)-- (4,0);
\draw [line width=0.8pt] (4,0)-- (4,4);
\draw [line width=0.8pt] (4,4)-- (0,4);
\draw [line width=0.8pt] (0,4)-- (0,0);
\draw [line width=0.8pt] (0,0)-- (4,4);
\draw [line width=0.8pt] (0,4)-- (2,2);
\begin{scriptsize}
\draw [fill=black] (0,0) circle (2pt);
\draw [fill=black] (4,0) circle (2pt);
\draw [fill=black] (4,4) circle (2pt);
\draw [fill=black] (0,4) circle (2pt);
\draw [fill=white] (2,2) circle (2pt);
\draw (0,-0.5) node {$i_1(0,0)$};
\draw (4,-0.5) node {$i_2(1,0)$};
\draw (4,4.2) node {$i_3(1,1)$};
\draw (0,4.2) node {$i_4(0,1)$};
\draw (3,2) node {$i_0(0.5,0.5)$};
\draw (3,1) node {$K_1$};
\draw (2,3.5) node {$K_2$};
\draw (0.5,2) node {$K_3$};
\end{scriptsize}
\end{tikzpicture}
\caption{Patch considered in Examples~\ref{ex:ImplementationP1} and~\ref{ex:ImplementationP1_0}.}\label{fig:row_sum_hanging}
\end{figure}

In standard finite element methods, 
the matrix and right-hand side are typically assembled cell-wise. This approach can be performed
also for the set of non-conforming basis
functions~$\varphi_{i_p}^{\mathrm{nc}}$, $p=0,\dots,4$, leading to 
\begin{equation*}
\begin{pmatrix}
a_{00} & a_{01} & a_{02} & a_{03} & a_{04} \\
a_{10} & a_{11} & a_{12} & a_{13} & a_{14} \\
a_{20} & a_{21} & a_{22} & a_{23} & a_{24} \\
a_{30} & a_{31} & a_{32} & a_{33} & a_{34} \\
a_{40} & a_{41} & a_{42} & a_{43} & a_{44} \\
\end{pmatrix},
\quad
\begin{pmatrix} b_0 \\ b_1 \\ b_2 \\ b_3 \\ b_4 \end{pmatrix}.
\end{equation*}
The $p$th equation of the corresponding linear system corresponds to the
non-con\-forming test function $\varphi_{i_p}^{\mathrm{nc}}$. In view of the
above relations between conforming and non-conforming basis functions, equations
corresponding to conforming test functions are obtained by adding $\frac12$ of
the 0th equation to the 1st and 3rd equations. To enforce continuity, the 0th
equation is then replaced by the relation
\eqref{eq:hanging_node_representation} with $q=0$. This leads to the following
matrix and right-hand side
%The step from the non-conforming to the conforming test functions yields
\begin{equation}\label{eq:matrix_conf_test}
\begin{pmatrix}
1 & -\frac{1}{2} & 0 & -\frac{1}{2} & 0  \\
a_{10}+\frac{a_{00}}{2} & a_{11}+\frac{a_{01}}{2} & a_{12}+\frac{a_{02}}{2} & a_{13}+\frac{a_{03}}{2} & a_{14}+\frac{a_{04}}{2} \\
a_{20} & a_{21} & a_{22} & a_{23} & a_{24} \\
a_{30}+\frac{a_{00}}{2} & a_{31}+\frac{a_{01}}{2} & a_{32}+\frac{a_{02}}{2} & a_{33}+\frac{a_{03}}{2} & a_{34}+\frac{a_{04}}{2} \\
a_{40} & a_{41} & a_{42} & a_{43} & a_{44} \\
\end{pmatrix},
\quad
 \begin{pmatrix}
   0 \\ b_1+\frac{b_0}{2} \\ b_2 \\ b_3+\frac{b_0}{2} \\ b_4 \\
 \end{pmatrix}.
\end{equation}
\hfill $\Box$ \end{example}

Usually, a system with matrix and right-hand side from \eqref{eq:matrix_conf_test} is used for computing the 
numerical solution on grids with hanging nodes. But for algebraically stabilized schemes there is a new
question: Which matrix should be used for computing the limiters? 
%A necessary assumption for performing the limiting procedure is 
The proofs of the DMP use the assumption
that the diagonal entries of the corresponding matrix are positive. However, this 
property cannot be guaranteed for the matrix from \eqref{eq:matrix_conf_test}. In fact, numerical studies, 
which are not reported here for the sake of brevity, that used the limiters computed with this matrix led in several cases
to unsatisfactory results, e.g., solutions obtained with the Kuzmin limiter showed spurious oscillation.
Consequently, an additional step has to be performed for algebraically stabilized schemes, namely 
%to transform the system 
a transformation of the system to a form corresponding
also to conforming ansatz functions. This means, the constraints for the hanging nodes are inserted in the 
other equations such that the corresponding matrix entries become zero.

Both steps, to the conforming test functions and to the conforming ansatz functions, extend the matrix stencil by few 
entries in rows that belong to test functions for non-hanging nodes which are located in a vicinity of hanging nodes. 

\begin{example}[System corresponding to conforming ansatz and test functions]
\label{ex:ImplementationP1_0}Consider the matrix and right-hand side from \eqref{eq:matrix_conf_test}. Inserting  
the equation for the finite element coefficient of
$\varphi_{i_0}^{\mathrm{nc}}$, which is the 0th equation, into 
the other equations, yields a matrix of the following form
\begin{equation}\label{eq:matrix_conf_test_ansatz}
\begin{pmatrix}
1 & -\frac{1}{2} & 0 & -\frac{1}{2} & 0  \\
0 & a_{11}+\frac{a_{01}}{2}+ \frac{a_{10}}2+\frac{a_{00}}{4} & a_{12}+\frac{a_{02}}{2} & a_{13}+\frac{a_{03}}{2} + \frac{a_{10}}2+\frac{a_{00}}{4} & a_{14}+\frac{a_{04}}{2} \\
0 & a_{21} + \frac{a_{20}}2 & a_{22} & a_{23} + \frac{a_{20}}2 & a_{24} \\
0 & a_{31}+\frac{a_{01}}{2} + \frac{a_{30}}2+\frac{a_{00}}{4} & a_{32}+\frac{a_{02}}{2} & a_{33}+\frac{a_{03}}{2}+ \frac{a_{30}}2+\frac{a_{00}}{4}  & a_{34}+\frac{a_{04}}{2} \\
0 & a_{41} + \frac{a_{40}}2 & a_{42} & a_{43} + \frac{a_{40}}2 & a_{44} \\
\end{pmatrix}.
\end{equation}
\hfill $\Box$ \end{example}

The computation of the limiters is performed for the submatrix from \eqref{eq:matrix_conf_test_ansatz} that corresponds to the rows and columns connected 
with non-hanging nodes. The Kuzmin and the MUAS limiter can be applied in a straightforward way. 
The set $N_i$ in Step~\ref{it:BJK2} of the BJK limiter is computed by exploring the entries of the $i$th row 
and taking all column indices of the corresponding sparsity pattern. Let $\Delta_i = \mathrm{conv}\left\{x_j\ : \ j\in N_i\right\}$
be the convex hull of the nodes belonging to $N_i$.  Then, the same definition as given in \cite{BJK17}
can be used:
\begin{equation*}
\gamma_i=\frac{\underset{x_j\in \partial \Delta_i}{\mathrm{max}}|x_i-x_j|}{\mathrm{dist}(x_i,\partial\Delta_i)},\quad i=1,\dots,M.
\end{equation*}
%can be used as given in \cite{BJK17}. 
%Note that if hanging nodes are considered, the set $\Delta_i$ is in 
%general different from the support of a conforming basis function, in 
%contrast to conforming triangulations.
% -- THIS CAN BE ALSO THE CASE WHEN NO HANGING NODES ARE CONSIDERED.

\section{Numerical Studies}\label{sec:numres}
This section presents numerical studies of algebraically stabilized schemes on adaptively 
refined grids. Both, grids with hanging nodes and grids with conforming closure will be
considered and the results will be compared. Given a grid with hanging nodes that should be 
closed in a conforming way, then the closure might
increase the 
largest angle or decrease the smallest angle of the triangles of the grid. The refinement with 
hanging nodes was performed such that there is not more than one hanging node per edge. 

Using adaptively refined grids requires some criterion for controlling the local refinement. 
Usually, a posteriori error estimators or indicators are utilized. For the considered methods
there is a residual-based a posteriori error estimator for the AFC schemes with Kuzmin and
BJK limiter on conforming grids, which was proposed and analyzed in \cite{Jha20}. 
In this paper, actually two different techniques for calculating an upper bound for the error in the energy norm of
solutions computed with AFC schemes on conforming grids are proposed.
One of them uses a residual-based approach, which is referred to as \afce technique, and the other one utilizes
the SUPG estimator from \cite{JN13}, which is referred to as \afcse technique.
It was observed in \cite{Jha20} that the \afce technique provides better results with respect to the refinement of the grids and hence we decided to use it as basis for our numerical studies. 

Denote by $\|\cdot\|_{0,\omega}$ the norm of $L^2(\omega)$ for some set $\omega$. 
In the \afce technique,  the error $u-u_h$ in the energy norm
is bounded, i.e.,
\begin{equation}\label{eq:post_upper_bound}
\|u-u_h\|_a^2\leq \eta^2=\eta_1^2+\eta_2^2+\eta_3^2,
\end{equation}
where $\|u\|_a^2 = \varepsilon \|\nabla u\|_{0,\Omega}^2+\sigma_0\|u\|_{0,\Omega}^2$, with $-(\nabla\cdot\bb(\bx))/2 + c(\bx) \ge \sigma_0 > 0$ being assumed in $\Omega$,  
and
\begin{eqnarray*}
\eta_1^2& := & \sum_{K\in \mathcal{T}_h}\mathrm{min}\left\{ \frac{4 C_I^2}{\sigma_0},\ \frac{4C_I^2h_K^2}{\varepsilon}\right\}\|R_K(u_h)\|_{0,K}^2, \\ 
\eta_2^2 & :=&   \sum_{F\in \cF_h} \min \left\{\frac{4 C_F^2 h_F}{\varepsilon}, \frac{4 C_F^2}{\sigma_0^{1/2}\varepsilon^{1/2}} \right\}\| R_F(u_h)\|_{0,F}^2,\\
\eta_3^2 & := & \sum_{\me}\min \Bigg\{\frac{4\kappa_1h_E^2}{
\varepsilon}, \frac{4\kappa_2}{\sigma_0}\Bigg\} (1-\alpha_E)^2|d_E|^2 h_E^{1-d} \|\nabla u_h \cdot\bt_E\|_{0,E}^2\nonumber,
\end{eqnarray*}
$u_h$ is the solution of the algebraically stabilized scheme,
$\alpha_E=\alpha_{ij}$ and $d_E=d_{ij}$ for an edge $E$ with endpoints $x_i$,
$x_j$, $\bt_E$ is the unit tangent vector along the edge $E$, 
$R_K(u_h)$ and $R_F(u_h)$ stand for the residuals on mesh cell $K$ and on the facet $F$  given by
\begin{eqnarray*}
R_K(u_h)&:=&(f+\varepsilon \Delta u_h-\boldsymbol{b}\cdot \nabla u_h-cu_h)|_K,\nonumber \\
R_F(u_h)&:=&
\left\lbrace
\begin{array}{lc}
-\varepsilon [|\nabla u_h \cdot \bn_F |]_F &  \mathrm{if}\ F\in \cF_{h,\Omega},\\
g-\varepsilon (\nabla u_h \cdot \bn_F)&
\mathrm{if}\ F\in \cF_{h,N},\\
0 & \mathrm{if}\ F\in \cF_{h,D},
\end{array}\right.\nonumber
\end{eqnarray*} 
$\bn_F$ is the unit normal on facet $F$, and $[|\cdot|]_F$ denotes the jump across  $F$.

The constants $C_I$ and $C_F$ appear from the interpolation and facet estimates and were set to unity in the simulations. The constants $\kappa_1$ and $\kappa_2$ are given by
$$
\kappa_1=CC_{\mathrm{edge,max}}\left(1+\left(1+C_I\right)^2\right),\quad \kappa_2=CC_{\mathrm{inv}}^2C_{\mathrm{edge,max}}\left( 1+\left(1+C_I\right)^2\right),
$$
where $C$ is a general constant independent of $h$, $C_{\mathrm{inv}}$ is an
inverse inequality constant, and $C_{\mathrm{edge,max}}$ is a computable
constant given by \cite[Remark~9]{Jha20}.  Likewise as the other constants, $C$
and $C_{\mathrm{inv}}$ were set to unity in our simulations. In
\cite{Jha20}, the error estimator $\eta$ was applied to the two
above-described AFC schemes on conforming grids. 

\iffalse
\begin{remark}
The residual based estimator presented above is based on the edge formulation of the stabilization term given in \cite{BJKR18}. Here, $d_E$ and $\alpha_E$ correspond to $d_{ij}$ and $\alpha_{ij}$ along the edge $E$ with end points $x_i$ and $x_j$. For the modified Kuzmin limiter, the edge formulation still holds as the matrix $(B-D)$ is symmetric and the derivation of the edge formulation can be applied. For the modified Kuzmin limiter, $\eta_3$ is given by
$$
\eta_3^2  :=  \sum_{\me}\min \Bigg\{\frac{4\kappa_1h_E^2}{
\varepsilon}, \frac{4\kappa_2}{\sigma_0}\Bigg\} (|b_E-d_E|)^2 h_E^{1-d} \|\nabla u_h \cdot\bt_E\|_{0,E}^2.
$$
\end{remark}
\fi

\begin{remark}\rm
The expression $\eta$ from \eqref{eq:post_upper_bound} can be computed also for the MUAS method and for 
all methods on grids with hanging nodes. Then, it is just an error indicator,  i.e., there is no analysis. 
In practice, often error indicators are used for controlling the adaptive grid refinement, 
like the popular gradient indicator. In preliminary studies, we could observe that for the AFC methods, 
the use of $\eta$ on grids with hanging nodes led to a quite similar adaptive refinement process as 
for grids with conforming closure, i.e., the refinement starts at the strongest singularities (exponential
layers) and regions with weaker singularities (parabolic layers) are refined somewhat later. For this reason, 
we applied $\eta$ also for the AFC methods on grids with hanging nodes. In contrast, we detected that 
applying $\eta$ for the MUAS method results in a simultaneous refinement in all regions with singularities
and considerably different adaptive grids compared with the AFC methods. This situation made it 
difficult to compare the computational results. Neglecting the term $\eta_3$ for the MUAS method, 
which results in a standard residual-based error indicator, led to a similar behavior of the 
adaptive grid refinement process as for the AFC schemes. For this reason, the adaptive grid 
refinement for the MUAS method was controlled  on all grids with $(\eta_1^2+\eta_2^2)^{1/2}$. 
\hfill$\Box$\end{remark}

A grid with conforming closure contains regularly refined cells and closure cells. Both types might be
marked for refinement by the error indicator. In the first step of the refinement process, parents of 
closure cells are marked for refinement if one of its children is marked for refinement. Note that parents
of closure cells are regularly refined cells on a coarser grid. Then, all closure cells are removed 
and all marked cells are refined regularly. Finally, the refined grid is closed. In the case of grids 
with hanging nodes, all marked cells are refined regularly. Then, a procedure is applied that 
refines all cells regularly that have an edge with more than one hanging node, until such cells are not contained 
any longer in the grid. 
The adaptive refinement process for the first two examples was stopped after the first adaptively
refined grid where 
the number of degrees of freedom ($\#\mathrm{dof}$) was $\gtrsim 2.5\times 10^{5}$. The given numbers 
$\#\mathrm{dof}$  contain always the hanging and the Dirichlet 
nodes.

\begin{remark}
Comparative studies for the solution of the nonlinear problem arising in the AFC schemes were performed in \cite{JJ18, JJ19}. It was found that the simplest fixed point iteration scheme was the most efficient one. A brief description of this scheme is as follows.
The matrix form of the algebraic stabilization given in \eqref{eq:matrix_rep_asm} is 
reformulated as 
$$
\left(A +D\right)U=b + \left(D-B(U)\right)U,
$$
with the artificial diffusion matrix $D$ from \eqref{eq:mat_art_diff}. The matrix on the left-hand 
side is by construction an M-matrix. Then, a fixed point iteration of the form 
\begin{equation}\label{eq:fp_equ}
\left(A +D\right)\tilde U^{\mu}=b + \left(D-B(U^{\mu})\right)U^{\mu}, \quad 
 U^{\mu+1} = \omega \tilde U^{\mu} + (1-\omega) U^{\mu},
\end{equation}
is applied, 
where $\mu$ denotes the $\mu^{\mathrm{th}}$ iterative step and $\omega\in \mathbb{R}^+$ is a damping parameter, which is chosen dynamically. Using a sparse direct solver for the linear 
systems of equations in \eqref{eq:fp_equ}
exploits that 
the matrix on the left-hand side does not change during the iteration and
hence its factorization needs to be computed only once. Also for iterative solvers, 
method \eqref{eq:fp_equ} is well suited, because they usually converge quickly since the matrix 
is an M-matrix, compare \cite{JJ19}.
A detailed description of this scheme, in particular of the dynamic damping procedure, 
can be found in \cite{JJ19}, where it is referred to as `fixed-point right-hand side'. 
The nonlinear loops  were stopped if $10,000$ iteration steps were reached or if $\mathtt{res}\leq \varepsilon_{\mathrm{thresh}}\sqrt{\#\mathrm{dof}}$, where $\mathtt{res}$ is the Euclidean norm of the residual vector and $\varepsilon_{\mathrm{thresh}}$ is the stopping threshold.
If not mentioned otherwise, then $\varepsilon_{\mathrm{thresh}} = 10^{-10}$.
\hfill $\Box$\end{remark}

All  schemes were used with $\mathbb{P}_1$ finite elements.
The matrices were assembled exactly and the arising systems of linear equations were solved using the sparse direct solver \textsc{UMFPACK},  \cite{Dav04}.  All simulations were performed with the in-house code \textsc{ParMooN}, \cite{WB16,GJM+16}.

The numerical results will be compared on the basis of the satisfaction of the global DMP, the 
accuracy of solutions, e.g., measured by sharpness of layers, and efficiency, measured by the 
number of iterations and rejections for the solver of the nonlinear problem. After having rejected
a step, the damping factor is decreased, but this step is computationally as expensive as an 
accepted step. 

\subsection{Solution Becoming Locally Diffusion-Dominated under Adaptive Grid Refinement}\label{ex:boundary_layer}
\begin{figure}[t!]
\centerline{\includegraphics[width=0.48\textwidth]{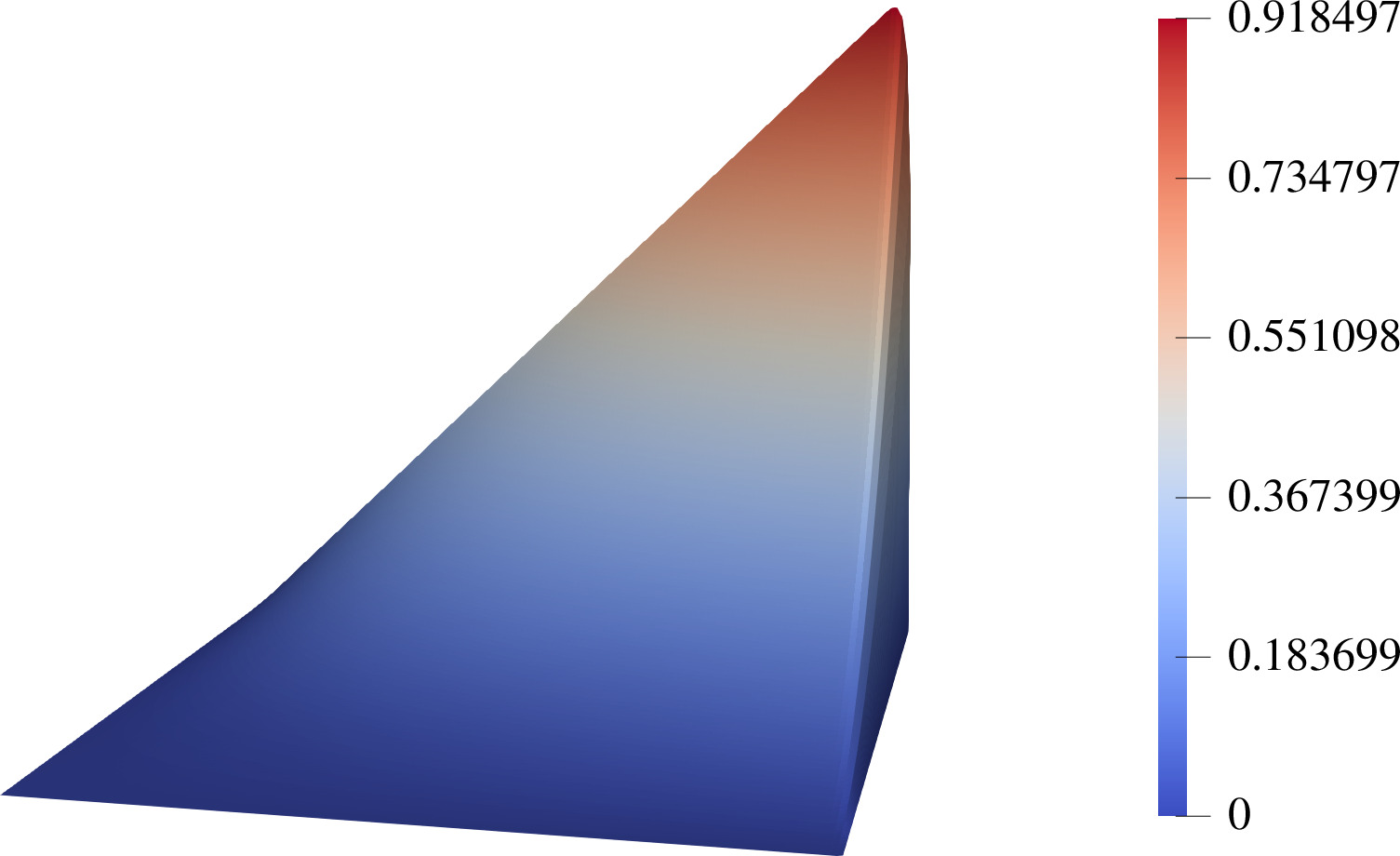}}
\caption{Example~\ref{ex:boundary_layer}: Solution computed with AFC scheme and Kuzmin limiter, level~7 with uniform refinement.}
\label{fig:boundary_layer}
\end{figure}

This example, presented in \cite{JMT96}, is given in $\Omega = (0,1)^2$ with 
$\bb = (2,3)^T$, $c=1$, and $\partial \Omega=\Gamma_D$.  The solution
\begin{eqnarray*}
u(x,y)& = & xy^2-y^2\exp \left(\frac{2(x-1)}{\varepsilon}\right)-x\exp \left(\frac{3(y-1)}{\varepsilon}\right)\\
&&+\exp \left(\frac{2(x-1)+3(y-1)}{\varepsilon}\right),
\end{eqnarray*}
defines the right-hand side $f$ and the Dirichlet boundary condition $u_b$. It possesses 
boundary layers at $x=1$ and $y=1$, see Figure~\ref{fig:boundary_layer}. We consider the case $\ep = 10^{-2}$, i.e., the 
discrete problem is convection-dominated on coarse grids (the layers are not resolved) and it becomes diffusion-dominated
on finer grids. 

The initial mesh (level 0) was defined by dividing the domain into two triangles by joining the points $(0,0)$ and $(1,1)$. 
The simulations were started with the level 2 grid obtained by uniform refinement (i.e., $\# \mathrm{dof}=25$) and initially uniform refinement was applied until level~5 (i.e., $\# \mathrm{dof}=1089$). After that, adaptive refinement was performed. 

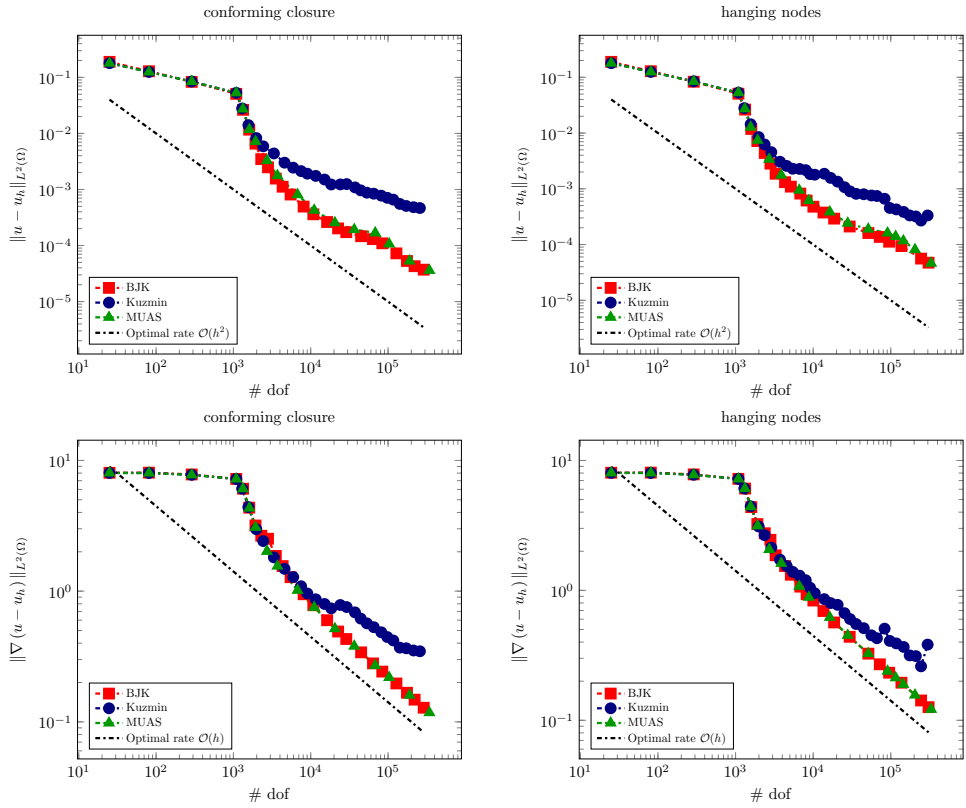
\begin{figure}[t!]
\centerline{
\begin{tikzpicture}[scale=0.6]
\begin{loglogaxis}[
    legend pos=south west, xlabel = $\#\ \mathrm{dof}$,, ylabel = $\|u-u_h\|_{L^2(\Omega)}$,
    legend cell align ={left},, title = {conforming closure},
    legend style={nodes={scale=0.75, transform shape}}]
\addplot[color=red,  mark=square*, line width = 0.5mm, dashdotted,,mark options = {scale= 1.5, solid}]
coordinates{( 25.0 , 0.18936692307798 )( 81.0 , 0.12858582339951 )( 289.0 , 0.083133041199614 )( 1089.0 , 0.050794543442462 )( 1333.0 , 0.026376741871119 )( 1602.0 , 0.011665181060873 )( 1931.0 , 0.0065781487615571 )( 2294.0 , 0.003479964451252 )( 2802.0 , 0.0024810648866703 )( 3536.0 , 0.001564758395006 )( 4333.0 , 0.0011193294730669 )( 5494.0 , 0.00081288590358186 )( 8098.0 , 0.00049524573663877 )( 10811.0 , 0.00036053711015347 )( 16299.0 , 0.00026258032679044 )( 22607.0 , 0.00020199143239416 )( 28839.0 , 0.0001730229237211 )( 44870.0 , 0.00014648503563946 )( 63932.0 , 0.00013011666121605 )( 84300.0 , 0.00010950445410602 )( 128333.0 , 7.2466982542363e-05 )( 175735.0 , 5.2960783158185e-05 )( 220676.0 , 4.2773045804573e-05 )( 289435.0 , 3.7067229110494e-05 )};
\addlegendentry{BJK} 
\addplot[color=navy_blue,  mark=oplus*, line width = 0.5mm, dashdotted,,mark options = {scale= 1.5, solid}]
coordinates{( 25.0 , 0.17913707011159 )( 81.0 , 0.12492901044605 )( 289.0 , 0.084504351507092 )( 1089.0 , 0.052989348097259 )( 1298.0 , 0.027871840264888 )( 1563.0 , 0.014077519211039 )( 1977.0 , 0.0082688121800735 )( 2433.0 , 0.0058830729999192 )( 3341.0 , 0.0043905657871804 )( 4587.0 , 0.0030074523845259 )( 5933.0 , 0.0024489439603241 )( 7505.0 , 0.0021274627963332 )( 9046.0 , 0.0019128424135109 )( 11611.0 , 0.0017460701272899 )( 15114.0 , 0.0015041223642896 )( 18415.0 , 0.0012266497504296 )( 24135.0 , 0.0012305426151133 )( 29166.0 , 0.001247636088881 )( 37373.0 , 0.0010892466880151 )( 44327.0 , 0.0009644080164438 )( 53249.0 , 0.00087650365180999 )( 65382.0 , 0.00084049314890724 )( 81816.0 , 0.00077603468814067 )( 97826.0 , 0.00071009713672393 )( 116962.0 , 0.00065886652117708 )( 142170.0 , 0.00055197481817266 )( 172137.0 , 0.00050793668431811 )( 213692.0 , 0.00048274956128489 )( 258076.0 , 0.00046650890681124 )};
\addlegendentry{Kuzmin} 
\addplot[color=dark_green,  mark=triangle*, line width = 0.5mm, dashdotted,,mark options = {scale= 1.5, solid}]
coordinates{( 25.0 , 0.17897867237439 )( 81.0 , 0.12461630959925 )( 289.0 , 0.084336193930147 )( 1089.0 , 0.052892559439686 )( 1327.0 , 0.027276020281876 )( 1610.0 , 0.011773609494832 )( 1935.0 , 0.0071941331645901 )( 2714.0 , 0.0032793162163063 )( 3711.0 , 0.0017662268697191 )( 6767.0 , 0.00080756771600921 )( 11089.0 , 0.00042643654521049 )( 20531.0 , 0.00025544969073557 )( 36565.0 , 0.00019267515734411 )( 68275.0 , 0.00016815413432938 )( 103436.0 , 0.00010665386573491 )( 190474.0 , 5.2992812768723e-05 )( 342988.0 , 3.5976717186205e-05 )};
\addlegendentry{MUAS} 
\addplot[color=black,  line width = 0.5mm, dashdotted,,mark options = {scale= 1.0, solid}]
coordinates{( 25.0 , 0.04 )( 81.0 , 0.012345679012345678 )( 289.0 , 0.0034602076124567475 )( 1089.0 , 0.0009182736455463728 )( 1333.0 , 0.0007501875468867217 )( 1602.0 , 0.0006242197253433209 )( 1931.0 , 0.0005178663904712584 )( 2294.0 , 0.00043591979075850045 )( 2802.0 , 0.00035688793718772306 )( 3536.0 , 0.0002828054298642534 )( 4333.0 , 0.00023078698361412417 )( 5494.0 , 0.00018201674554058973 )( 8098.0 , 0.00012348728081007656 )( 10811.0 , 9.249838127832763e-05 )( 16299.0 , 6.135345726731701e-05 )( 22607.0 , 4.4234086787278275e-05 )( 28839.0 , 3.467526613266757e-05 )( 44870.0 , 2.2286605749944284e-05 )( 63932.0 , 1.5641619220421698e-05 )( 84300.0 , 1.1862396204033214e-05 )( 128333.0 , 7.792228031761121e-06 )( 175735.0 , 5.6903860926963895e-06 )( 220676.0 , 4.531530388442785e-06 )( 289435.0 , 3.455007169139876e-06 )};
\addlegendentry{Optimal rate $\mathcal{O}(h^2)$} 
\end{loglogaxis}
\end{tikzpicture}\hspace*{1em}
\begin{tikzpicture}[scale=0.6]
\begin{loglogaxis}[
    legend pos=south west, xlabel = $\#\ \mathrm{dof}$, ylabel=$\|u-u_h\|_{L^2(\Omega)}$,
    legend cell align ={left},, title = hanging nodes,
    legend style={nodes={scale=0.75, transform shape}}]
\addplot[color=red,  mark=square*, line width = 0.5mm, dashdotted,,mark options = {scale= 1.5, solid}]
coordinates{( 25.0 , 0.18936692307798 )( 81.0 , 0.12858582339951 )( 289.0 , 0.083133041199614 )( 1089.0 , 0.050794543442462 )( 1329.0 , 0.026309193129987 )( 1594.0 , 0.011834695018292 )( 1910.0 , 0.0071723991847422 )( 2380.0 , 0.0043889744353469 )( 2816.0 , 0.0028239402269767 )( 3285.0 , 0.0018615179462317 )( 4344.0 , 0.0013135136316747 )( 5157.0 , 0.0010938437005513 )( 6827.0 , 0.0008207121513921 )( 8235.0 , 0.00061577999124851 )( 10068.0 , 0.00047745906244549 )( 13570.0 , 0.00037526864073644 )( 18642.0 , 0.00029171510061933 )( 29525.0 , 0.00021059415908463 )( 51452.0 , 0.00016097792468037 )( 72267.0 , 0.00013695202618703 )( 95250.0 , 0.00011242218771074 )( 137241.0 , 9.3842373747877e-05 )( 245784.0 , 5.5993719370482e-05 )( 306768.0 , 4.7268790554065e-05 )};
\addlegendentry{BJK} 
\addplot[color=navy_blue,  mark=oplus*, line width = 0.5mm, dashdotted,,mark options = {scale= 1.5, solid}]
coordinates{( 25.0 , 0.17913707011159 )( 81.0 , 0.12492901044605 )( 289.0 , 0.084504351507092 )( 1089.0 , 0.052989348097259 )( 1298.0 , 0.027977187946959 )( 1567.0 , 0.01436573907381 )( 1986.0 , 0.0084682185146477 )( 2369.0 , 0.006233178498921 )( 2887.0 , 0.0045434125399818 )( 3736.0 , 0.0030556254906906 )( 4560.0 , 0.0025713918829811 )( 5448.0 , 0.0022902060832118 )( 6684.0 , 0.0022703080527144 )( 8015.0 , 0.0021627376631793 )( 9113.0 , 0.0018336237361195 )( 10479.0 , 0.001789890716382 )( 13913.0 , 0.001861593001554 )( 16795.0 , 0.0015677999976555 )( 20647.0 , 0.0013433364985745 )( 25331.0 , 0.0010698753516306 )( 29856.0 , 0.00089956635492679 )( 35765.0 , 0.00081225586514596 )( 44862.0 , 0.00079857120345615 )( 56330.0 , 0.0007585137198415 )( 66248.0 , 0.00074293927754812 )( 83213.0 , 0.00066261602666805 )( 97131.0 , 0.00045290163716014 )( 117756.0 , 0.00042697684273625 )( 146554.0 , 0.00038463579979929 )( 175288.0 , 0.00033282803880996 )( 210243.0 , 0.00031910524635208 )( 245311.0 , 0.00026929681550672 )( 298373.0 , 0.00033269747814449 )};
\addlegendentry{Kuzmin} 
\addplot[color=dark_green,  mark=triangle*, line width = 0.5mm, dashdotted,,mark options = {scale= 1.5, solid}]
coordinates{( 25.0 , 0.17897867237439 )( 81.0 , 0.12461630959925 )( 289.0 , 0.084336193930147 )( 1089.0 , 0.052892559439686 )( 1319.0 , 0.027174615645002 )( 1566.0 , 0.012733536205319 )( 1923.0 , 0.0074196604451641 )( 2712.0 , 0.0033856159064852 )( 3881.0 , 0.0017609730260902 )( 6620.0 , 0.00094356393053168 )( 8806.0 , 0.00061423629397114 )( 16384.0 , 0.0003844260824191 )( 27902.0 , 0.00024155877190944 )( 51292.0 , 0.00018915040946489 )( 90849.0 , 0.00016086565501242 )( 115614.0 , 0.00013989901309822 )( 145041.0 , 0.00011658502968561 )( 205044.0 , 8.108802240879e-05 )( 330817.0 , 4.6204755206897e-05 )};
\addlegendentry{MUAS} 
\addplot[color=black,  line width = 0.5mm, dashdotted,,mark options = {scale= 1.0, solid}]
coordinates{( 25.0 , 0.04 )( 81.0 , 0.012345679012345678 )( 289.0 , 0.0034602076124567475 )( 1089.0 , 0.0009182736455463728 )( 1329.0 , 0.0007524454477050414 )( 1594.0 , 0.0006273525721455458 )( 1910.0 , 0.0005235602094240838 )( 2380.0 , 0.0004201680672268908 )( 2816.0 , 0.0003551136363636364 )( 3285.0 , 0.00030441400304414006 )( 4344.0 , 0.00023020257826887662 )( 5157.0 , 0.00019391118867558658 )( 6827.0 , 0.00014647722279185587 )( 8235.0 , 0.00012143290831815422 )( 10068.0 , 9.932459276916965e-05 )( 13570.0 , 7.369196757553427e-05 )( 18642.0 , 5.3642313056539e-05 )( 29525.0 , 3.386960203217612e-05 )( 51452.0 , 1.943559045323797e-05 )( 72267.0 , 1.3837574549932888e-05 )( 95250.0 , 1.0498687664041995e-05 )( 137241.0 , 7.286452299240023e-06 )( 245784.0 , 4.0686130911694824e-06 )( 306768.0 , 3.2597924164189223e-06 )};
\addlegendentry{Optimal rate $\mathcal{O}(h^2)$} 

\end{loglogaxis}
\end{tikzpicture}}
\centerline{
\begin{tikzpicture}[scale=0.6]
\begin{loglogaxis}[
    legend pos=south west, xlabel = $\#\ \mathrm{dof}$,ylabel = $\|\nabla\left(u-u_h\right)\|_{L^2(\Omega)}$,
    legend cell align ={left},, title = {conforming closure},
   legend style={nodes={scale=0.75, transform shape}}]
\addplot[color=red,  mark=square*, line width = 0.5mm, dashdotted,,mark options = {scale= 1.5, solid}]
coordinates{( 25.0 , 8.0424725406542 )( 81.0 , 8.0609946536039 )( 289.0 , 7.8105937015305 )( 1089.0 , 7.218507640822 )( 1333.0 , 6.0747681971772 )( 1602.0 , 4.3529433537147 )( 1931.0 , 3.17856381442 )( 2294.0 , 2.6518235262102 )( 2802.0 , 2.5185191302926 )( 3536.0 , 1.8595506255914 )( 4333.0 , 1.5552644661375 )( 5494.0 , 1.2787737863922 )( 8098.0 , 0.946357784552 )( 10811.0 , 0.78013226895274 )( 16299.0 , 0.59912260823659 )( 22607.0 , 0.49099046224222 )( 28839.0 , 0.43039869208635 )( 44870.0 , 0.34029445013933 )( 63932.0 , 0.28004272180811 )( 84300.0 , 0.24239559538406 )( 128333.0 , 0.19666870628084 )( 175735.0 , 0.16706628298625 )( 220676.0 , 0.14789570098187 )( 289435.0 , 0.12827061726442 )};
\addlegendentry{BJK} 
\addplot[color=navy_blue,  mark=oplus*, line width = 0.5mm, dashdotted,,mark options = {scale= 1.5, solid}]
coordinates{( 25.0 , 8.0205632126557 )( 81.0 , 8.0101826237912 )( 289.0 , 7.7559388638086 )( 1089.0 , 7.2026545553424 )( 1298.0 , 6.077623359542 )( 1563.0 , 4.3946349022518 )( 1977.0 , 2.9748789078349 )( 2433.0 , 2.415851033824 )( 3341.0 , 1.8124507586868 )( 4587.0 , 1.4809332813893 )( 5933.0 , 1.2820742605109 )( 7505.0 , 1.0906395084311 )( 9046.0 , 0.95946515131841 )( 11611.0 , 0.86639385472912 )( 15114.0 , 0.79766473306313 )( 18415.0 , 0.73831977708659 )( 24135.0 , 0.78355577124564 )( 29166.0 , 0.75570982654881 )( 37373.0 , 0.68820675106047 )( 44327.0 , 0.61690966794237 )( 53249.0 , 0.56697594323974 )( 65382.0 , 0.53015456827991 )( 81816.0 , 0.48452120597348 )( 97826.0 , 0.44560239185104 )( 116962.0 , 0.41930589709134 )( 142170.0 , 0.36952531402813 )( 172137.0 , 0.36548389310875 )( 213692.0 , 0.35161148900045 )( 258076.0 , 0.34678772627458 )};
\addlegendentry{Kuzmin} 
\addplot[color=dark_green,  mark=triangle*, line width = 0.5mm, dashdotted,,mark options = {scale= 1.5, solid}]
coordinates{( 25.0 , 8.0206175662346 )( 81.0 , 8.0108332048875 )( 289.0 , 7.7573416070453 )( 1089.0 , 7.2036916503474 )( 1327.0 , 6.0705994897884 )( 1610.0 , 4.3190636322318 )( 1935.0 , 3.0761881435576 )( 2714.0 , 2.0146257854485 )( 3711.0 , 1.5535890794866 )( 6767.0 , 1.0142080973916 )( 11089.0 , 0.75147264109627 )( 20531.0 , 0.516503250177 )( 36565.0 , 0.37938683967609 )( 68275.0 , 0.27011736593374 )( 103436.0 , 0.21858266531645 )( 190474.0 , 0.15994320068451 )( 342988.0 , 0.1176072280535 )};
\addlegendentry{MUAS} 
\addplot[color=black,  line width = 0.5mm, dashdotted,,mark options = {scale= 1.0, solid}]
coordinates{( 25.0 , 8.94427190999916 )( 81.0 , 4.969039949999533 )( 289.0 , 2.6306682088232822 )( 1089.0 , 1.3551927136362363 )( 1333.0 , 1.2248979932114525 )( 1602.0 , 1.1173358719233182 )( 1931.0 , 1.017709575931423 )( 2294.0 , 0.9337235037831065 )( 2802.0 , 0.8448525755274976 )( 3536.0 , 0.7520710469952335 )( 4333.0 , 0.6793923514643422 )( 5494.0 , 0.6033518799847892 )( 8098.0 , 0.4969653525349158 )( 10811.0 , 0.4301124998842224 )( 16299.0 , 0.35029546747657764 )( 22607.0 , 0.2974359991234359 )( 28839.0 , 0.26334489223323687 )( 44870.0 , 0.21112368768067824 )( 63932.0 , 0.17687068281895504 )( 84300.0 , 0.15402854413408715 )( 128333.0 , 0.12483771891348483 )( 175735.0 , 0.10668070202896482 )( 220676.0 , 0.09520010912223563 )( 289435.0 , 0.08312649600626597 )};
\addlegendentry{Optimal rate $\mathcal{O}(h)$} 
\end{loglogaxis}
\end{tikzpicture}\hspace*{1em}
\begin{tikzpicture}[scale=0.6]
\begin{loglogaxis}[
    legend pos=south west, xlabel = $\#\ \mathrm{dof}$,ylabel=$\|\nabla\left(u-u_h\right)\|_{L^2(\Omega)}$,
    legend cell align ={left},, title = {hanging nodes},
    legend style={nodes={scale=0.75, transform shape}}]
\addplot[color=red,  mark=square*, line width = 0.5mm, dashdotted,,mark options = {scale= 1.5, solid}]
coordinates{( 25.0 , 8.0424725406542 )( 81.0 , 8.0609946536039 )( 289.0 , 7.8105937015305 )( 1089.0 , 7.218507640822 )( 1329.0 , 6.0757678812282 )( 1594.0 , 4.3728853327837 )( 1910.0 , 3.2383568149604 )( 2380.0 , 2.7518541604003 )( 2816.0 , 2.4507643382804 )( 3285.0 , 1.8610175566629 )( 4344.0 , 1.5387748764759 )( 5157.0 , 1.3190953085567 )( 6827.0 , 1.0684216598307 )( 8235.0 , 0.93364739924626 )( 10068.0 , 0.83475084234687 )( 13570.0 , 0.69345758224009 )( 18642.0 , 0.56613551243336 )( 29525.0 , 0.43916845419131 )( 51452.0 , 0.32574841770827 )( 72267.0 , 0.26977714561922 )( 95250.0 , 0.23237502398895 )( 137241.0 , 0.19473874573061 )( 245784.0 , 0.14206270407487 )( 306768.0 , 0.12616813644844 )};
\addlegendentry{BJK} 
\addplot[color=navy_blue,  mark=oplus*, line width = 0.5mm, dashdotted,,mark options = {scale= 1.5, solid}]
coordinates{( 25.0 , 8.0205632126557 )( 81.0 , 8.0101826237912 )( 289.0 , 7.7559388638086 )( 1089.0 , 7.2026545553424 )( 1298.0 , 6.0872315302773 )( 1567.0 , 4.4341709480366 )( 1986.0 , 3.1016081242496 )( 2369.0 , 2.6595441765946 )( 2887.0 , 2.1316203895719 )( 3736.0 , 1.7247243301142 )( 4560.0 , 1.5384981764917 )( 5448.0 , 1.3882590510226 )( 6684.0 , 1.2979725574174 )( 8015.0 , 1.1969884930654 )( 9113.0 , 1.0482343248803 )( 10479.0 , 0.95004678433918 )( 13913.0 , 0.8559225488176 )( 16795.0 , 0.79831846242501 )( 20647.0 , 0.77380519580204 )( 25331.0 , 0.66782694172255 )( 29856.0 , 0.60169407113346 )( 35765.0 , 0.55066127755048 )( 44862.0 , 0.51159250412924 )( 56330.0 , 0.45095252279399 )( 66248.0 , 0.42817046765165 )( 83213.0 , 0.50737331906699 )( 97131.0 , 0.40823141395443 )( 117756.0 , 0.39064151719845 )( 146554.0 , 0.36681911042127 )( 175288.0 , 0.31594483783525 )( 210243.0 , 0.31187427847024 )( 245311.0 , 0.26020357734758 )( 298373.0 , 0.38260111184426 )};
\addlegendentry{Kuzmin} 
\addplot[color=dark_green,  mark=triangle*, line width = 0.5mm, dashdotted,,mark options = {scale= 1.5, solid}]
coordinates{( 25.0 , 8.0206175662346 )( 81.0 , 8.0108332048875 )( 289.0 , 7.7573416070453 )( 1089.0 , 7.2036916503474 )( 1319.0 , 6.0722382990134 )( 1566.0 , 4.4061251465044 )( 1923.0 , 3.1353848505752 )( 2712.0 , 2.0622558712326 )( 3881.0 , 1.6142252398098 )( 6620.0 , 1.0753861894191 )( 8806.0 , 0.88511704604191 )( 16384.0 , 0.62006119103729 )( 27902.0 , 0.45070959998104 )( 51292.0 , 0.32646390240696 )( 90849.0 , 0.23867517046579 )( 115614.0 , 0.212380889997 )( 145041.0 , 0.18898240974725 )( 205044.0 , 0.15640519903503 )( 330817.0 , 0.12136894098901 )};
\addlegendentry{MUAS} 
\addplot[color=black,  line width = 0.5mm, dashdotted,,mark options = {scale= 1.0, solid}]
coordinates{( 25.0 , 8.94427190999916 )( 81.0 , 4.969039949999533 )( 289.0 , 2.6306682088232822 )( 1089.0 , 1.3551927136362363 )( 1329.0 , 1.22673994612146 )( 1594.0 , 1.1201362168464564 )( 1910.0 , 1.0232890201933018 )( 2380.0 , 0.9166984970282113 )( 2816.0 , 0.8427498280790526 )( 3285.0 , 0.7802743146408705 )( 4344.0 , 0.6785316179351948 )( 5157.0 , 0.622753865785812 )( 6827.0 , 0.5412526633502247 )( 8235.0 , 0.4928141806363819 )( 10068.0 , 0.44570078027566806 )( 13570.0 , 0.38390615409376877 )( 18642.0 , 0.3275433194450438 )( 29525.0 , 0.26026756245132093 )( 51452.0 , 0.197157756394406 )( 72267.0 , 0.1663584957249427 )( 95250.0 , 0.14490471120044368 )( 137241.0 , 0.12071828609817173 )( 245784.0 , 0.090206575050486 )( 306768.0 , 0.0807439461064286 )};
\addlegendentry{Optimal rate $\mathcal{O}(h)$} 
\end{loglogaxis}
\end{tikzpicture}}
\caption{Example~\ref{ex:boundary_layer}: $L^2(\Omega)$ error (top) and $L^2(\Omega)$ error of the gradient (bottom); 
grids with conforming closure (left) and grids with hanging nodes (right).}\label{fig:l2_error_boundary_layer}
\end{figure}

Since the solution is known, errors of the discrete approximations computed with the algebraically
stabilized schemes can be computed. Figure~\ref{fig:l2_error_boundary_layer} presents the errors 
in the $L^2(\Omega)$ norm and in the $L^2(\Omega)$ norm of 
the gradient. It can be seen that  the solutions computed with the AFC scheme with BJK limiter
and with the MUAS scheme are likewise accurate. 
On both types of grids, the optimal convergence order of the error in the $L^2(\Omega)$ norm 
of the gradient can be seen. 
It has to be noted that the error estimator is for the error in the energy norm, which is dominated here 
by the $L^2(\Omega)$ error of the gradient, and not for the $L^2(\Omega)$ norm, such that the 
adaptive grids might not be always suitable for an optimal error convergence in the $L^2(\Omega)$ norm. 
The solutions 
obtained with the AFC scheme and Kuzmin limiter seem not to converge on grids with conforming closure and they 
converge slower on grids with hanging nodes. This behavior on conforming grids 
was already observed 
for a similar example in \cite{Jha20}. In fact, the analysis from \cite{BJK16} predicts that convergence
can be expected for this method in the diffusion-dominated case only 
if the grid satisfies certain conditions, e.g., if the grid is Delaunay.

\begin{figure}[t!]
\centerline{
\begin{tikzpicture}[scale=0.6]
\begin{semilogxaxis}[
    legend pos=north east, xlabel = $\#\ \mathrm{dof}$, ylabel = iterations+rejections,
    legend cell align ={left}, title = {conforming closure},
     legend style={nodes={scale=0.75, transform shape}}]
\addplot[color=red,  mark=square*, line width = 0.5mm, dashdotted,,mark options = {scale= 1.5, solid}]
coordinates{( 25.0 , 111.0 )( 81.0 , 102.0 )( 289.0 , 118.0 )( 1089.0 , 99.0 )( 1333.0 , 69.0 )( 1602.0 , 68.0 )( 1931.0 , 67.0 )( 2294.0 , 102.0 )( 2802.0 , 124.0 )( 3536.0 , 105.0 )( 4333.0 , 102.0 )( 5494.0 , 91.0 )( 8098.0 , 75.0 )( 10811.0 , 74.0 )( 16299.0 , 60.0 )( 22607.0 , 49.0 )( 28839.0 , 43.0 )( 44870.0 , 70.0 )( 63932.0 , 30.0 )( 84300.0 , 26.0 )( 128333.0 , 27.0 )( 175735.0 , 22.0 )( 220676.0 , 16.0 )( 289435.0 , 14.0 )};
\addlegendentry{BJK} 
\addplot[color=navy_blue,  mark=oplus*, line width = 0.5mm, dashdotted,,mark options = {scale= 1.5, solid}] 
coordinates{( 25.0 , 26.0 )( 81.0 , 31.0 )( 289.0 , 42.0 )( 1089.0 , 33.0 )( 1298.0 , 34.0 )( 1563.0 , 33.0 )( 1977.0 , 33.0 )( 2433.0 , 33.0 )( 3341.0 , 32.0 )( 4587.0 , 32.0 )( 5933.0 , 32.0 )( 7505.0 , 31.0 )( 9046.0 , 31.0 )( 11611.0 , 30.0 )( 15114.0 , 27.0 )( 18415.0 , 24.0 )( 24135.0 , 18.0 )( 29166.0 , 18.0 )( 37373.0 , 16.0 )( 44327.0 , 16.0 )( 53249.0 , 16.0 )( 65382.0 , 16.0 )( 81816.0 , 15.0 )( 97826.0 , 16.0 )( 116962.0 , 15.0 )( 142170.0 , 13.0 )( 172137.0 , 12.0 )( 213692.0 , 11.0 )( 258076.0 , 11.0 )};
\addlegendentry{Kuzmin} 
\addplot[color=dark_green,  mark=triangle*, line width = 0.5mm, dashdotted,,mark options = {scale= 1.5, solid}] 
coordinates{( 25.0 , 26.0 )( 81.0 , 36.0 )( 289.0 , 47.0 )( 1089.0 , 52.0 )( 1327.0 , 50.0 )( 1610.0 , 50.0 )( 1935.0 , 48.0 )( 2714.0 , 66.0 )( 3711.0 , 62.0 )( 6767.0 , 49.0 )( 11089.0 , 37.0 )( 20531.0 , 34.0 )( 36565.0 , 32.0 )( 68275.0 , 30.0 )( 103436.0 , 42.0 )( 190474.0 , 39.0 )( 342988.0 , 22.0 )};
\addlegendentry{MUAS} 
\end{semilogxaxis}
\end{tikzpicture}\hspace*{1em}
\begin{tikzpicture}[scale=0.6]
\begin{semilogxaxis}[
    legend pos=north east, xlabel = $\#\ \mathrm{dof}$, ylabel = iterations+rejections,
    legend cell align ={left}, title = {hanging nodes},
     legend style={nodes={scale=0.75, transform shape}, text opacity =10, fill opacity = 0.6}]
\addplot[color=red,  mark=square*, line width = 0.5mm, dashdotted,,mark options = {scale= 1.5, solid}]
coordinates{( 25.0 , 111.0 )( 81.0 , 102.0 )( 289.0 , 118.0 )( 1089.0 , 99.0 )( 1329.0 , 79.0 )( 1594.0 , 64.0 )( 1910.0 , 61.0 )( 2380.0 , 98.0 )( 2816.0 , 116.0 )( 3285.0 , 74.0 )( 4344.0 , 108.0 )( 5157.0 , 107.0 )( 6827.0 , 90.0 )( 8235.0 , 83.0 )( 10068.0 , 75.0 )( 13570.0 , 79.0 )( 18642.0 , 75.0 )( 29525.0 , 69.0 )( 51452.0 , 54.0 )( 72267.0 , 47.0 )( 95250.0 , 37.0 )( 137241.0 , 34.0 )( 245784.0 , 35.0 )( 306768.0 , 31.0 )};
\addlegendentry{BJK} 
\addplot[color=navy_blue,  mark=oplus*, line width = 0.5mm, dashdotted,,mark options = {scale= 1.5, solid}] 
coordinates{( 25.0 , 26.0 )( 81.0 , 31.0 )( 289.0 , 42.0 )( 1089.0 , 33.0 )( 1298.0 , 33.0 )( 1567.0 , 33.0 )( 1986.0 , 33.0 )( 2369.0 , 33.0 )( 2887.0 , 33.0 )( 3736.0 , 32.0 )( 4560.0 , 32.0 )( 5448.0 , 32.0 )( 6684.0 , 31.0 )( 8015.0 , 31.0 )( 9113.0 , 31.0 )( 10479.0 , 31.0 )( 13913.0 , 29.0 )( 16795.0 , 26.0 )( 20647.0 , 23.0 )( 25331.0 , 20.0 )( 29856.0 , 19.0 )( 35765.0 , 18.0 )( 44862.0 , 17.0 )( 56330.0 , 16.0 )( 66248.0 , 15.0 )( 83213.0 , 15.0 )( 97131.0 , 15.0 )( 117756.0 , 15.0 )( 146554.0 , 13.0 )( 175288.0 , 13.0 )( 210243.0 , 12.0 )( 245311.0 , 14.0 )( 298373.0 , 12.0 )};
\addlegendentry{Kuzmin} 
\addplot[color=dark_green,  mark=triangle*, line width = 0.5mm, dashdotted,,mark options = {scale= 1.5, solid}] 
coordinates{( 25.0 , 26.0 )( 81.0 , 36.0 )( 289.0 , 47.0 )( 1089.0 , 52.0 )( 1319.0 , 49.0 )( 1566.0 , 48.0 )( 1923.0 , 46.0 )( 2712.0 , 44.0 )( 3881.0 , 50.0 )( 6620.0 , 44.0 )( 8806.0 , 43.0 )( 16384.0 , 43.0 )( 27902.0 , 38.0 )( 51292.0 , 32.0 )( 90849.0 , 26.0 )( 115614.0 , 27.0 )( 145041.0 , 32.0 )( 205044.0 , 24.0 )( 330817.0 , 26.0 )};
\addlegendentry{MUAS} 
\end{semilogxaxis}
\end{tikzpicture}}
\caption{Example~\ref{ex:boundary_layer}: Number of iterations and rejections on grids with conforming closure (left) and on grids with hanging nodes (right).}\label{fig:iteration_2d_sol_boundary_layer}
\end{figure}
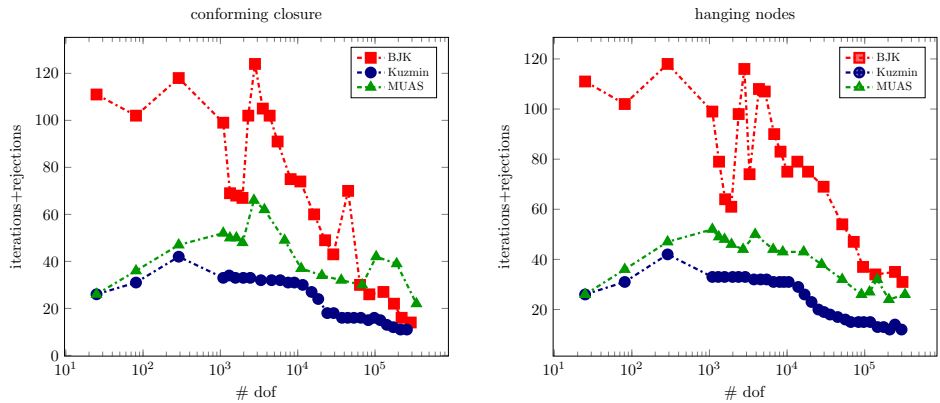

Figure~\ref{fig:iteration_2d_sol_boundary_layer} presents results concerning the 
efficiency of the methods. It can be observed that the AFC scheme with Kuzmin limiter needs usually
the smallest number of iterations and the AFC scheme with BJK limiter often the largest number. 
But altogether, no difficulties arose 
for solving the nonlinear problems.

\subsection{A Convection-Dominated Problem with Interior and Boundary Layers}\label{ex:hmm86}
\begin{figure}[t!]
\centerline{\includegraphics[width=0.48\textwidth]{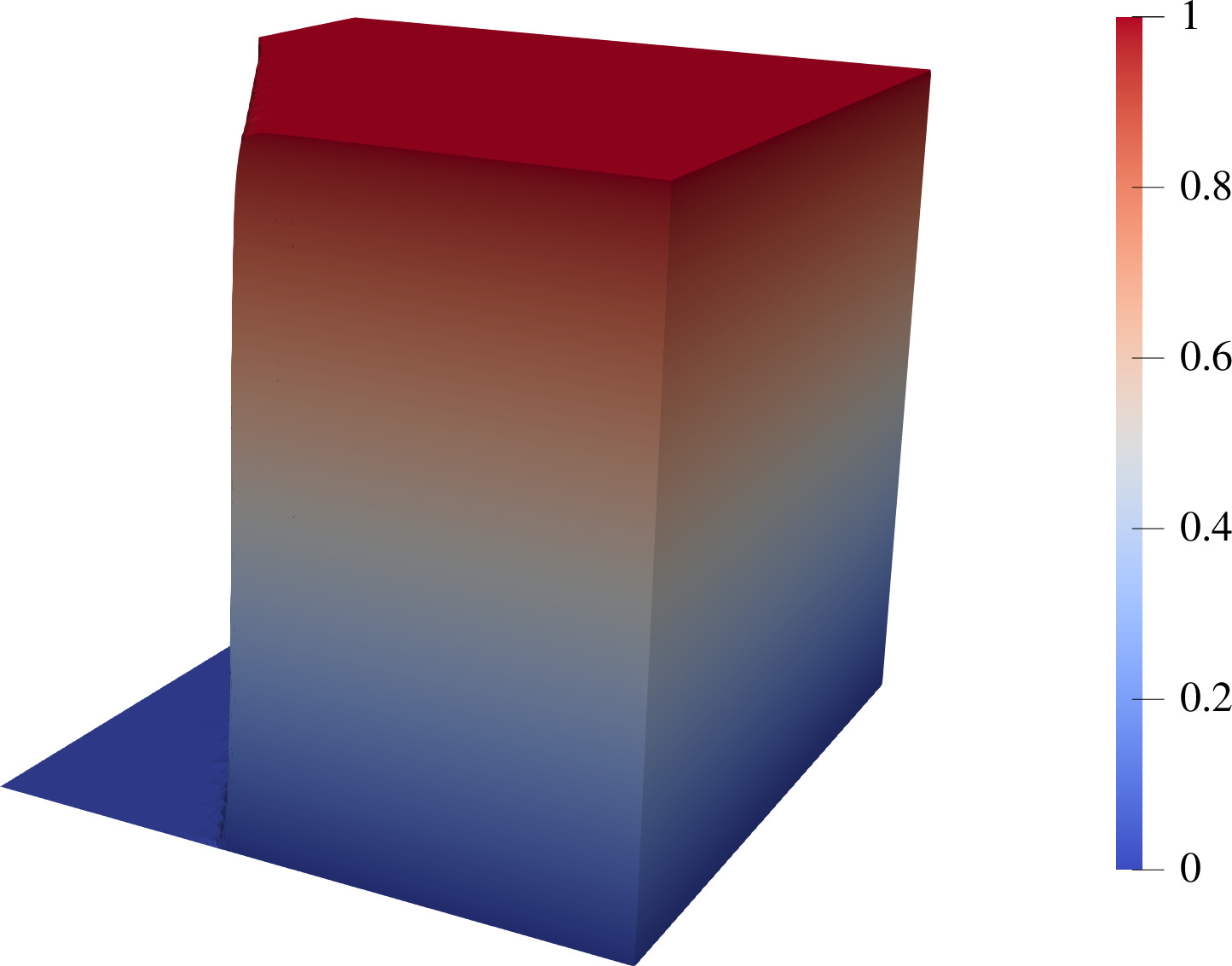}}
\caption{Example~\ref{ex:hmm86}: Solution to the interior and boundary layer example, computed with the
         BJK limiter, level~9.}
\label{fig:hmm86_sol}
\end{figure}

This standard example was proposed in \cite{HMM86}.
It is given in $\Omega = (0,1)^2$ with 
$\bb = (\cos(-\pi/3),$ $ \sin(-\pi/3))^T$, $c=f=0$, and the Dirichlet boundary condition
$$
u_b = \begin{cases} 1 & (y=1\wedge x>0) \mbox{ or } (x=0 \wedge y>0.7),\\
0 & \mbox{else}.
\end{cases}
$$
Here, the convection-dominated case $\ep = 10^{-6}$ is considered. The solution exhibits an interior layer in the
direction of the convection starting from the jump of the boundary condition at the left boundary and two exponential layers at the
right and the lower boundary, see Figure~\ref{fig:hmm86_sol}. An analytic solution to this problem is not available, but the 
solution satisfies the global maximum principle, i.e., $u\in [0,1]$. In the numerical studies, the satisfaction of the global DMP,
the accuracy by considering the width of the interior layer along a cut line, and the efficiency will be studied. 
In addition, the impact of relaxing the stopping criterion of the iteration on the quantities of interest will be investigated. 

The initial grid (level~0) was constructed by dividing the unit square with the diagonal from $(0,1)$ to $(1,0)$,
as advised in \cite{JK07_1}. The simulations were started on level~2 and uniform refinement was performed
until level~5.

The satisfaction of the global DMP is studied by evaluating the quantity
\begin{equation}\label{rem:var_def}
\mathrm{osc}_{\mathrm{max}}(u_h):= \max_{(x,y)\in\overline\Omega} u_h(x,y) -1 -
\min_{(x,y)\in\overline\Omega} u_h(x,y).
\end{equation}
It turned out that these values were for all schemes and all grids 
at most of the order of round-off errors. Hence, the corresponding numerical solutions satisfy the global DMP.

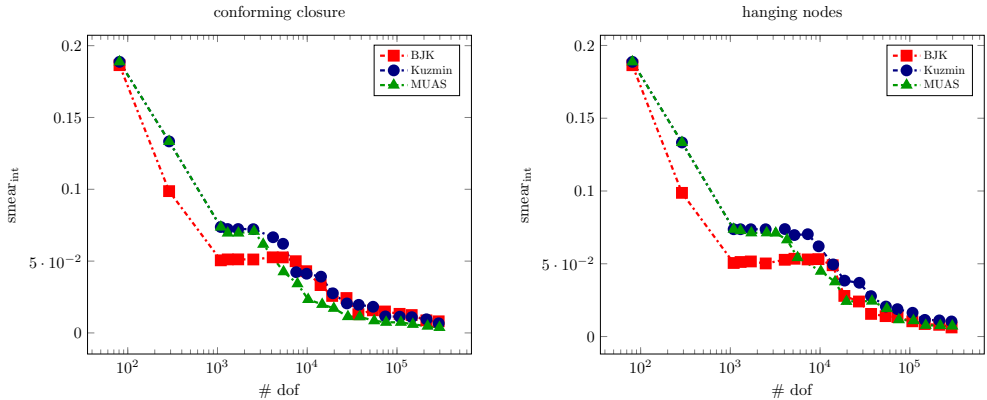
\begin{figure}[t!]
\centerline{
\begin{tikzpicture}[scale=0.6]
\begin{semilogxaxis}[
    legend pos=north east, xlabel = $\#\ \mathrm{dof}$,ylabel = $\mathrm{smear}_{\mathrm{int}}$,
    legend cell align ={left},, title = conforming closure,
   legend style={nodes={scale=0.75, transform shape}}]
\addplot[color=red,  mark=square*, line width = 0.5mm,dashdotted,, mark options = {scale= 1.5, solid}] 
coordinates{( 81.0 , 0.186423 )( 289.0 , 0.0987196 )( 1089.0 , 0.0505368 )( 1302.0 , 0.0511881 )( 1701.0 , 0.0512178 )( 2507.0 , 0.0511415 )( 4128.0 , 0.0526285 )( 5337.0 , 0.0525815 )( 7435.0 , 0.0499309 )( 9733.0 , 0.0429741 )( 14090.0 , 0.0333 )( 18961.0 , 0.0257425 )( 27369.0 , 0.0243013 )( 37169.0 , 0.0130342 )( 53950.0 , 0.0157404 )( 73406.0 , 0.0149364 )( 106897.0 , 0.0132817 )( 145753.0 , 0.0124581 )( 212612.0 , 0.0087442 )( 291208.0 , 0.00808405 )};
\addlegendentry{BJK} 
\addplot[color=navy_blue,  mark=oplus*, line width = 0.5mm,dashdotted,, mark options = {scale= 1.5, solid}] 
coordinates{( 81.0 , 0.188744 )( 289.0 , 0.133318 )( 1089.0 , 0.0737231 )( 1297.0 , 0.0723471 )( 1698.0 , 0.0723471 )( 2506.0 , 0.0722505 )( 4161.0 , 0.066619 )( 5375.0 , 0.0620606 )( 7541.0 , 0.0422916 )( 9845.0 , 0.0412139 )( 14232.0 , 0.0391718 )( 19283.0 , 0.0275836 )( 27668.0 , 0.020536 )( 37471.0 , 0.0195653 )( 54111.0 , 0.0182201 )( 74202.0 , 0.0115128 )( 107320.0 , 0.0112945 )( 146372.0 , 0.0109641 )( 213139.0 , 0.0094977 )( 291618.0 , 0.00668284 )};
\addlegendentry{Kuzmin} 
\addplot[color=dark_green,  mark=triangle*, line width = 0.5mm,dashdotted,, mark options = {scale= 1.5, solid}] 
coordinates{( 81.0 , 0.188744 )( 289.0 , 0.133313 )( 1089.0 , 0.0737168 )( 1301.0 , 0.0695135 )( 1710.0 , 0.0695135 )( 2553.0 , 0.070741 )( 3231.0 , 0.0617165 )( 5448.0 , 0.0425846 )( 7790.0 , 0.034233 )( 10246.0 , 0.0232862 )( 14678.0 , 0.0198739 )( 19790.0 , 0.0171224 )( 28407.0 , 0.011352 )( 38509.0 , 0.0112376 )( 55764.0 , 0.00846041 )( 75659.0 , 0.00738945 )( 110555.0 , 0.00739026 )( 150244.0 , 0.00596554 )( 219212.0 , 0.00468934 )( 297237.0 , 0.00381394 )};
\addlegendentry{MUAS} 
\end{semilogxaxis}
\end{tikzpicture}\hspace*{1em}
\begin{tikzpicture}[scale=0.6]
\begin{semilogxaxis}[
    legend pos=north east, xlabel = $\#\ \mathrm{dof}$,ylabel = $\mathrm{smear}_{\mathrm{int}}$,
    legend cell align ={left},, title = hanging nodes,
   legend style={nodes={scale=0.75, transform shape}}]
\addplot[color=red,  mark=square*, line width = 0.5mm,dashdotted,, mark options = {scale= 1.5, solid}] 
coordinates{( 81.0 , 0.186423 )( 289.0 , 0.0987196 )( 1089.0 , 0.0505368 )( 1301.0 , 0.051222 )( 1692.0 , 0.051598 )( 2479.0 , 0.0502288 )( 4040.0 , 0.0527386 )( 5189.0 , 0.0534455 )( 7264.0 , 0.0529795 )( 9598.0 , 0.0531824 )( 13811.0 , 0.0489896 )( 18650.0 , 0.0279052 )( 27133.0 , 0.024085 )( 37016.0 , 0.0155494 )( 53781.0 , 0.0140154 )( 72922.0 , 0.0133549 )( 106720.0 , 0.0104791 )( 146592.0 , 0.0085885 )( 212665.0 , 0.0080533 )( 291947.0 , 0.00625034 )};
\addlegendentry{BJK} 
\addplot[color=navy_blue,  mark=oplus*, line width = 0.5mm,dashdotted,, mark options = {scale= 1.5, solid}] 
coordinates{( 81.0 , 0.188744 )( 289.0 , 0.133318 )( 1089.0 , 0.0737231 )( 1297.0 , 0.0737287 )( 1690.0 , 0.0737287 )( 2473.0 , 0.0737287 )( 4049.0 , 0.0738539 )( 5174.0 , 0.0697279 )( 7294.0 , 0.0703198 )( 9646.0 , 0.0619963 )( 13979.0 , 0.0495638 )( 18753.0 , 0.038418 )( 27406.0 , 0.0369356 )( 36929.0 , 0.0277897 )( 54143.0 , 0.0206215 )( 72799.0 , 0.0188266 )( 107270.0 , 0.0163337 )( 147056.0 , 0.0114522 )( 213183.0 , 0.0111028 )( 291522.0 , 0.0103782 )};
\addlegendentry{Kuzmin} 
\addplot[color=dark_green,  mark=triangle*, line width = 0.5mm,dashdotted,, mark options = {scale= 1.5, solid}] 
coordinates{( 81.0 , 0.188744 )( 289.0 , 0.133313 )( 1089.0 , 0.0737168 )( 1300.0 , 0.0727994 )( 1715.0 , 0.071239 )( 2542.0 , 0.071255 )( 3206.0 , 0.0711018 )( 4281.0 , 0.0664564 )( 5617.0 , 0.0542361 )( 10201.0 , 0.0448649 )( 14681.0 , 0.0376189 )( 19681.0 , 0.0241483 )( 37908.0 , 0.0243262 )( 55475.0 , 0.0192385 )( 75556.0 , 0.0115461 )( 109948.0 , 0.0109112 )( 150447.0 , 0.00771116 )( 218442.0 , 0.00755285 )( 297418.0 , 0.00722165 )};
\addlegendentry{MUAS} 
\end{semilogxaxis}
\end{tikzpicture}}
\caption{Example~\ref{ex:hmm86}: Thickness of interior layer, $\mathrm{smear}_{\mathrm{int}}$.}
\label{fig:thickness_hmm86_grid_1}
\end{figure}

To check the thickness of the interior layer, we follow the idea described in \cite[Eq.~(48)]{JK07_1} and define
\begin{equation}\label{eq:smear}
\mathrm{smear}_{\mathrm{int}}= x_2-x_1,
\end{equation}
where $x_1$ is the $x$-coordinate on the cut line $(x,0.25)$ with
$u_h(x_1,0.25)= 0.1$ and $x_2$ is the $x-$coordinate with $u_h(x_2,0.25)=
0.9$. The cut line was discretized with 100,000 equidistant intervals, 
where the discrete solutions were evaluated at the nodes. Then, the values for $x_1$ and $x_2$ were computed
by linear interpolation. The results, presented in Figure~\ref{fig:thickness_hmm86_grid_1}, show 
that  there are only minor differences between the solutions obtained with the different methods. 
On grids with hanging nodes, the AFC method with BJK limiter and the MUAS method computed usually 
a little bit sharper layers than the AFC method with Kuzmin limiter.

\begin{figure}[t!]
\centerline{
\begin{tikzpicture}[scale=0.6]
\begin{semilogxaxis}[scaled y ticks = false,
    legend pos=north west, xlabel = $\#\ \mathrm{dof}$, ylabel = iterations+rejections,
    legend cell align ={left}, title = {conforming closure},
     legend style={nodes={scale=0.75, transform shape}}]
\addplot[color=red,  mark=square*, line width = 0.5mm, dashdotted,,mark options = {scale= 1.5, solid}] 
coordinates{( 25.0 , 48.0 )( 81.0 , 296.0 )( 289.0 , 389.0 )( 1089.0 , 603.0 )( 1302.0 , 810.0 )( 1701.0 , 647.0 )( 2507.0 , 462.0 )( 4128.0 , 501.0 )( 5337.0 , 497.0 )( 7435.0 , 524.0 )( 9733.0 , 1404.0 )( 14090.0 , 1427.0 )( 18961.0 , 10078.0 )( 27369.0 , 706.0 )( 37169.0 , 1061.0 )( 53950.0 , 10067.0 )( 73406.0 , 10113.0 )( 106897.0 , 4650.0 )( 145753.0 , 10110.0 )( 212612.0 , 10118.0 )( 291208.0 , 10139.0 )};
\addlegendentry{BJK17} 
\addplot[color=navy_blue,  mark=oplus*, line width = 0.5mm, dashdotted,,mark options = {scale= 1.5, solid}] 
coordinates{( 25.0 , 21.0 )( 81.0 , 45.0 )( 289.0 , 42.0 )( 1089.0 , 57.0 )( 1297.0 , 64.0 )( 1698.0 , 80.0 )( 2506.0 , 91.0 )( 4161.0 , 60.0 )( 5375.0 , 100.0 )( 7541.0 , 80.0 )( 9845.0 , 112.0 )( 14232.0 , 106.0 )( 19283.0 , 131.0 )( 27668.0 , 107.0 )( 37471.0 , 105.0 )( 54111.0 , 104.0 )( 74202.0 , 121.0 )( 107320.0 , 112.0 )( 146372.0 , 186.0 )( 213139.0 , 164.0 )( 291618.0 , 651.0 )};
\addlegendentry{Kuzmin} 
\addplot[color=dark_green,  mark=triangle*, line width = 0.5mm, dashdotted,,mark options = {scale= 1.5, solid}]
coordinates{( 25.0 , 21.0 )( 81.0 , 45.0 )( 289.0 , 42.0 )( 1089.0 , 57.0 )( 1301.0 , 53.0 )( 1710.0 , 83.0 )( 2553.0 , 55.0 )( 3231.0 , 66.0 )( 5448.0 , 65.0 )( 7790.0 , 122.0 )( 10246.0 , 79.0 )( 14678.0 , 128.0 )( 19790.0 , 115.0 )( 28407.0 , 4584.0 )( 38509.0 , 185.0 )( 55764.0 , 206.0 )( 75659.0 , 227.0 )( 110555.0 , 228.0 )( 150244.0 , 274.0 )( 219212.0 , 294.0 )( 297237.0 , 290.0 )};
\addlegendentry{MUAS} 
\end{semilogxaxis}
\end{tikzpicture}\hspace*{1em}
\begin{tikzpicture}[scale=0.6]
\begin{semilogxaxis}[scaled y ticks = false,
    legend pos=north west, xlabel = $\#\ \mathrm{dof}$, ylabel = iterations+rejections,
    legend cell align ={left}, title = {hanging nodes},
     legend style={nodes={scale=0.75, transform shape}}]
\addplot[color=red,  mark=square*, line width = 0.5mm, dashdotted,,mark options = {scale= 1.5, solid}] 
coordinates{( 25.0 , 48.0 )( 81.0 , 296.0 )( 289.0 , 389.0 )( 1089.0 , 603.0 )( 1301.0 , 432.0 )( 1692.0 , 517.0 )( 2479.0 , 679.0 )( 4040.0 , 502.0 )( 5189.0 , 761.0 )( 7264.0 , 488.0 )( 9598.0 , 746.0 )( 13811.0 , 840.0 )( 18650.0 , 588.0 )( 27133.0 , 1481.0 )( 37016.0 , 1067.0 )( 53781.0 , 10134.0 )( 72922.0 , 10127.0 )( 106720.0 , 10105.0 )( 146592.0 , 4776.0 )( 212665.0 , 10134.0 )( 291947.0 , 10161.0 )};
\addlegendentry{BJK17} 
\addplot[color=navy_blue,  mark=oplus*, line width = 0.5mm, dashdotted,,mark options = {scale= 1.5, solid}] 
coordinates{( 25.0 , 21.0 )( 81.0 , 45.0 )( 289.0 , 42.0 )( 1089.0 , 57.0 )( 1297.0 , 39.0 )( 1690.0 , 46.0 )( 2473.0 , 50.0 )( 4049.0 , 55.0 )( 5174.0 , 49.0 )( 7294.0 , 56.0 )( 9646.0 , 55.0 )( 13979.0 , 53.0 )( 18753.0 , 67.0 )( 27406.0 , 99.0 )( 36929.0 , 88.0 )( 54143.0 , 107.0 )( 72799.0 , 94.0 )( 107270.0 , 117.0 )( 147056.0 , 121.0 )( 213183.0 , 187.0 )( 291522.0 , 148.0 )};
\addlegendentry{Kuzmin} 
\addplot[color=dark_green,  mark=triangle*, line width = 0.5mm, dashdotted,,mark options = {scale= 1.5, solid}]
coordinates{( 25.0 , 21.0 )( 81.0 , 45.0 )( 289.0 , 42.0 )( 1089.0 , 57.0 )( 1300.0 , 43.0 )( 1715.0 , 54.0 )( 2542.0 , 61.0 )( 3206.0 , 53.0 )( 4281.0 , 61.0 )( 5617.0 , 85.0 )( 10201.0 , 95.0 )( 14681.0 , 99.0 )( 19681.0 , 131.0 )( 37908.0 , 150.0 )( 55475.0 , 190.0 )( 75556.0 , 193.0 )( 109948.0 , 177.0 )( 150447.0 , 300.0 )( 218442.0 , 174.0 )( 297418.0 , 192.0 )};
\addlegendentry{MUAS} 
\end{semilogxaxis}
\end{tikzpicture}}
\centerline{
\begin{tikzpicture}[scale=0.6]
\begin{semilogxaxis}[
    legend pos=north west, xlabel = $\#\ \mathrm{dof}$, ylabel = iterations+rejections,
    legend cell align ={left}, title = {$\varepsilon_{\mathrm{thresh}}=10^{-8}$, conforming closure},
     legend style={nodes={scale=0.75, transform shape}}]
\addplot[color=red,  mark=square*, line width = 0.5mm, dashdotted,,mark options = {scale= 1.5, solid}] 
coordinates{( 25.0 , 41.0 )( 81.0 , 212.0 )( 289.0 , 278.0 )( 1089.0 , 341.0 )( 1302.0 , 644.0 )( 1701.0 , 457.0 )( 2507.0 , 271.0 )( 4128.0 , 282.0 )( 5337.0 , 320.0 )( 7435.0 , 328.0 )( 9731.0 , 1160.0 )( 14087.0 , 760.0 )( 18964.0 , 543.0 )( 27372.0 , 424.0 )( 37174.0 , 485.0 )( 53957.0 , 435.0 )( 73739.0 , 572.0 )( 106951.0 , 1689.0 )( 146421.0 , 1345.0 )( 212466.0 , 813.0 )( 291307.0 , 549.0 )};
\addlegendentry{BJK17} 
\addplot[color=navy_blue,  mark=oplus*, line width = 0.5mm, dashdotted,,mark options = {scale= 1.5, solid}] 
coordinates{( 25.0 , 15.0 )( 81.0 , 32.0 )( 289.0 , 29.0 )( 1089.0 , 38.0 )( 1297.0 , 31.0 )( 1698.0 , 39.0 )( 2506.0 , 39.0 )( 4161.0 , 32.0 )( 5375.0 , 50.0 )( 7541.0 , 37.0 )( 9845.0 , 47.0 )( 14232.0 , 59.0 )( 19283.0 , 67.0 )( 27668.0 , 51.0 )( 37473.0 , 49.0 )( 54111.0 , 55.0 )( 74202.0 , 55.0 )( 107320.0 , 53.0 )( 146360.0 , 61.0 )( 213142.0 , 57.0 )( 291630.0 , 63.0 )};
\addlegendentry{Kuzmin} 
\addplot[color=dark_green,  mark=triangle*, line width = 0.5mm, dashdotted,,mark options = {scale= 1.5, solid}]
coordinates{( 25.0 , 15.0 )( 81.0 , 32.0 )( 289.0 , 29.0 )( 1089.0 , 38.0 )( 1301.0 , 30.0 )( 1710.0 , 29.0 )( 2553.0 , 31.0 )( 3231.0 , 33.0 )( 5448.0 , 38.0 )( 7790.0 , 63.0 )( 10246.0 , 47.0 )( 14678.0 , 57.0 )( 19790.0 , 58.0 )( 28407.0 , 65.0 )( 38512.0 , 71.0 )( 55776.0 , 71.0 )( 75659.0 , 81.0 )( 110605.0 , 64.0 )( 150548.0 , 77.0 )( 219344.0 , 74.0 )( 298608.0 , 63.0 )};
\addlegendentry{MUAS} 
\end{semilogxaxis}
\end{tikzpicture}\hspace*{1em}
\begin{tikzpicture}[scale=0.6]
\begin{semilogxaxis}[scaled y ticks = false,
    legend pos=north west, xlabel = $\#\ \mathrm{dof}$, ylabel = iterations+rejections,
    legend cell align ={left}, title = {$\varepsilon_{\mathrm{thresh}}=10^{-8}$, hanging nodes},
     legend style={nodes={scale=0.75, transform shape}}]
\addplot[color=red,  mark=square*, line width = 0.5mm, dashdotted,,mark options = {scale= 1.5, solid}] 
coordinates{( 25.0 , 41.0 )( 81.0 , 212.0 )( 289.0 , 278.0 )( 1089.0 , 341.0 )( 1301.0 , 294.0 )( 1692.0 , 328.0 )( 2479.0 , 401.0 )( 4040.0 , 320.0 )( 5189.0 , 398.0 )( 7264.0 , 313.0 )( 9598.0 , 469.0 )( 13811.0 , 476.0 )( 18650.0 , 369.0 )( 27133.0 , 663.0 )( 37020.0 , 668.0 )( 53770.0 , 1034.0 )( 72936.0 , 548.0 )( 106735.0 , 1412.0 )( 146401.0 , 3101.0 )( 212494.0 , 821.0 )( 291858.0 , 1186.0 )};
\addlegendentry{BJK17} 
\addplot[color=navy_blue,  mark=oplus*, line width = 0.5mm, dashdotted,,mark options = {scale= 1.5, solid}] 
coordinates{( 25.0 , 15.0 )( 81.0 , 32.0 )( 289.0 , 29.0 )( 1089.0 , 38.0 )( 1297.0 , 28.0 )( 1690.0 , 29.0 )( 2473.0 , 33.0 )( 4049.0 , 29.0 )( 5174.0 , 27.0 )( 7294.0 , 29.0 )( 9646.0 , 33.0 )( 13979.0 , 32.0 )( 18753.0 , 32.0 )( 27406.0 , 42.0 )( 36929.0 , 51.0 )( 54145.0 , 44.0 )( 72805.0 , 44.0 )( 107270.0 , 48.0 )( 147056.0 , 51.0 )( 213220.0 , 58.0 )( 292315.0 , 53.0 )};
\addlegendentry{Kuzmin} 
\addplot[color=dark_green,  mark=triangle*, line width = 0.5mm, dashdotted,,mark options = {scale= 1.5, solid}]
coordinates{( 25.0 , 15.0 )( 81.0 , 32.0 )( 289.0 , 29.0 )( 1089.0 , 38.0 )( 1300.0 , 29.0 )( 1715.0 , 31.0 )( 2542.0 , 33.0 )( 3206.0 , 28.0 )( 4281.0 , 34.0 )( 5617.0 , 54.0 )( 10201.0 , 50.0 )( 14681.0 , 44.0 )( 19681.0 , 75.0 )( 37911.0 , 67.0 )( 55475.0 , 77.0 )( 75568.0 , 84.0 )( 109926.0 , 61.0 )( 150484.0 , 70.0 )( 218643.0 , 56.0 )( 297004.0 , 57.0 )};
\addlegendentry{MUAS} 
\end{semilogxaxis}
\end{tikzpicture}}
\centerline{
\begin{tikzpicture}[scale=0.6]
\begin{semilogxaxis}[
    legend pos=north west, xlabel = $\#\ \mathrm{dof}$, ylabel = iterations+rejections,
    legend cell align ={left}, title = {$\varepsilon_{\mathrm{thresh}}=10^{-6}$, conforming closure},
     legend style={nodes={scale=0.75, transform shape}}]
\addplot[color=red,  mark=square*, line width = 0.5mm, dashdotted,,mark options = {scale= 1.5, solid}] 
coordinates{( 25.0 , 34.0 )( 81.0 , 124.0 )( 289.0 , 176.0 )( 1089.0 , 144.0 )( 1302.0 , 411.0 )( 1701.0 , 191.0 )( 2507.0 , 128.0 )( 4128.0 , 129.0 )( 5337.0 , 169.0 )( 7435.0 , 212.0 )( 9825.0 , 194.0 )( 14077.0 , 70.0 )( 18849.0 , 92.0 )( 27431.0 , 89.0 )( 37243.0 , 75.0 )( 53915.0 , 57.0 )( 73399.0 , 58.0 )( 106850.0 , 50.0 )( 146007.0 , 28.0 )( 212464.0 , 18.0 )( 290753.0 , 10.0 )};
\addlegendentry{BJK17} 
\addplot[color=navy_blue,  mark=oplus*, line width = 0.5mm, dashdotted,,mark options = {scale= 1.5, solid}] 
coordinates{( 25.0 , 10.0 )( 81.0 , 20.0 )( 289.0 , 17.0 )( 1089.0 , 20.0 )( 1297.0 , 17.0 )( 1698.0 , 19.0 )( 2506.0 , 16.0 )( 4161.0 , 15.0 )( 5375.0 , 15.0 )( 7541.0 , 15.0 )( 9839.0 , 17.0 )( 14229.0 , 18.0 )( 18627.0 , 16.0 )( 27679.0 , 16.0 )( 37325.0 , 14.0 )( 54127.0 , 16.0 )( 74021.0 , 14.0 )( 107520.0 , 13.0 )( 147158.0 , 6.0 )( 213695.0 , 3.0 )( 292544.0 , 2.0 )};
\addlegendentry{Kuzmin} 
\addplot[color=dark_green,  mark=triangle*, line width = 0.5mm, dashdotted,,mark options = {scale= 1.5, solid}]
coordinates{( 25.0 , 10.0 )( 81.0 , 20.0 )( 289.0 , 17.0 )( 1089.0 , 20.0 )( 1301.0 , 19.0 )( 1710.0 , 18.0 )( 2553.0 , 16.0 )( 3229.0 , 16.0 )( 5446.0 , 16.0 )( 7790.0 , 15.0 )( 10246.0 , 18.0 )( 14673.0 , 16.0 )( 19050.0 , 16.0 )( 37068.0 , 15.0 )( 55702.0 , 12.0 )( 75471.0 , 10.0 )( 109740.0 , 3.0 )( 148946.0 , 3.0 )( 215616.0 , 3.0 )( 294308.0 , 2.0 )};
\addlegendentry{MUAS} 
\end{semilogxaxis}
\end{tikzpicture}\hspace*{1em}
\begin{tikzpicture}[scale=0.6]
\begin{semilogxaxis}[scaled y ticks = false,
    legend pos=north west, xlabel = $\#\ \mathrm{dof}$, ylabel = iterations+rejections,
    legend cell align ={left}, title = {$\varepsilon_{\mathrm{thresh}}=10^{-6}$, hanging nodes},
     legend style={nodes={scale=0.75, transform shape}}]
\addplot[color=red,  mark=square*, line width = 0.5mm, dashdotted,,mark options = {scale= 1.5, solid}] 
coordinates{( 25.0 , 34.0 )( 81.0 , 124.0 )( 289.0 , 176.0 )( 1089.0 , 144.0 )( 1301.0 , 153.0 )( 1692.0 , 134.0 )( 2479.0 , 159.0 )( 4040.0 , 145.0 )( 5189.0 , 120.0 )( 7264.0 , 135.0 )( 9602.0 , 231.0 )( 13824.0 , 104.0 )( 18653.0 , 126.0 )( 27099.0 , 126.0 )( 37166.0 , 115.0 )( 53771.0 , 151.0 )( 72903.0 , 85.0 )( 106709.0 , 87.0 )( 146098.0 , 30.0 )( 212260.0 , 38.0 )( 291053.0 , 13.0 )};
\addlegendentry{BJK17} 
\addplot[color=navy_blue,  mark=oplus*, line width = 0.5mm, dashdotted,,mark options = {scale= 1.5, solid}] 
coordinates{( 25.0 , 10.0 )( 81.0 , 20.0 )( 289.0 , 17.0 )( 1089.0 , 20.0 )( 1297.0 , 18.0 )( 1690.0 , 18.0 )( 2473.0 , 17.0 )( 4049.0 , 17.0 )( 5174.0 , 17.0 )( 7294.0 , 15.0 )( 9644.0 , 15.0 )( 13979.0 , 15.0 )( 18753.0 , 15.0 )( 27405.0 , 14.0 )( 36922.0 , 13.0 )( 54162.0 , 13.0 )( 73167.0 , 12.0 )( 107351.0 , 13.0 )( 147071.0 , 10.0 )( 213363.0 , 6.0 )( 292275.0 , 3.0 )};
\addlegendentry{Kuzmin} 
\addplot[color=dark_green,  mark=triangle*, line width = 0.5mm, dashdotted,,mark options = {scale= 1.5, solid}]
coordinates{( 25.0 , 10.0 )( 81.0 , 20.0 )( 289.0 , 17.0 )( 1089.0 , 20.0 )( 1300.0 , 17.0 )( 1713.0 , 18.0 )( 2543.0 , 18.0 )( 3204.0 , 15.0 )( 5169.0 , 14.0 )( 7732.0 , 17.0 )( 10217.0 , 15.0 )( 14647.0 , 16.0 )( 19201.0 , 13.0 )( 37341.0 , 12.0 )( 73259.0 , 10.0 )( 145048.0 , 4.0 )( 216521.0 , 3.0 )( 295712.0 , 3.0 )};
\addlegendentry{MUAS} 
\end{semilogxaxis}
\end{tikzpicture}}
\caption{Example~\ref{ex:hmm86} Number of iterations and rejections on grids with conforming closure (left) and on grids with hanging nodes (right),  $\varepsilon_{\mathrm{thresh}}=10^{-10}$ (top),
 $\varepsilon_{\mathrm{thresh}}=10^{-8}$ (middle),  $\varepsilon_{\mathrm{thresh}}=10^{-6}$ (bottom).}\label{fig:iteration_hmm86_grid_1}
\end{figure}
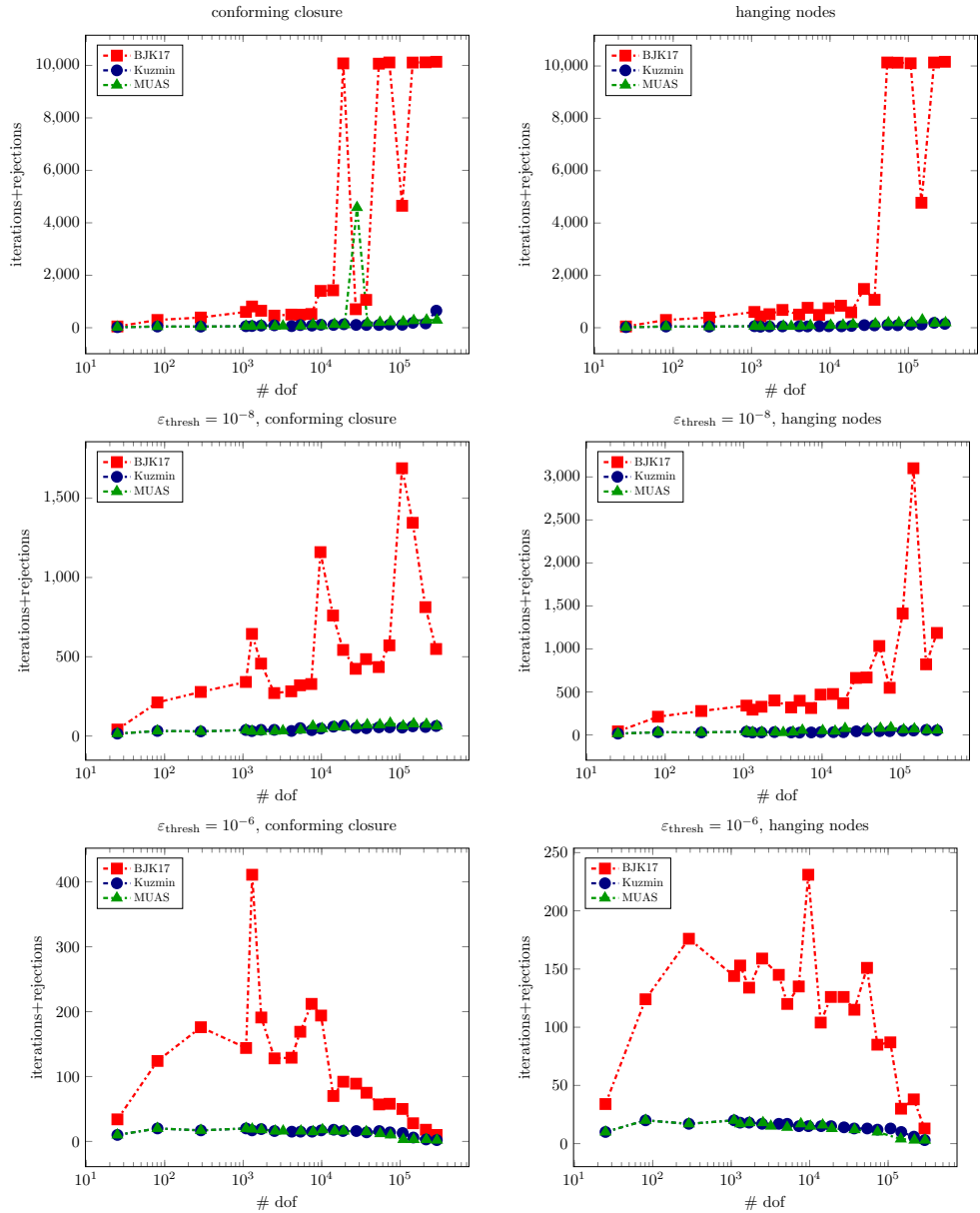

Figure~\ref{fig:iteration_hmm86_grid_1} presents the number of iterations and rejections. It can be 
observed that  the AFC method with BJK limiter sometimes stopped 
because the maximal number was reached, in particular on fine grids. The other two methods needed usually a similar
and much smaller number of iterations.
The rationale for choosing the hard stopping criterion with
 $\varepsilon_{\mathrm{thresh}}=10^{-10}$
is that analytic results, like the satisfaction of DMPs, can be proved only for the solution 
of the nonlinear discrete problem and thus an accurate solution seems to be advisable. 

The number of iterations and rejections for the weaker stopping criteria with $\varepsilon_{\mathrm{thresh}} = 10^{-6}$ 
and $\varepsilon_{\mathrm{thresh}} = 10^{-8}$ are depicted also in Figure~\ref{fig:iteration_hmm86_grid_1}.
It can be seen that in all situations 
the stopping criterion with respect to the residual could be satisfied now. The AFC scheme with Kuzmin 
limiter and the MUAS method require generally notably less iterations than the AFC scheme with BJK limiter. 

\begin{figure}[t!]
\centerline{
\begin{tikzpicture}[scale=0.6]
\begin{semilogxaxis}[
    legend pos=north west, xlabel = $\#\ \mathrm{dof}$,
    legend cell align ={left},, title = conforming closure,
    ymin = -0.01, ymax = 0.1,
    legend style={nodes={scale=0.75, transform shape}, }]
\addplot[color=black,  mark=diamond*, line width = 0.5mm, dashdotted,, mark options = {scale= 1.5, solid}] 
coordinates{( 25.0 , 1.0000000000065512e-05 )( 81.0 , 0.0 )( 289.0 , 0.0 )( 1089.0 , 2.719179992283216e-09 )( 1302.0 , 0.0 )( 1701.0 , 1.7323209533515183e-10 )( 2507.0 , 1.0000164297085945e-05 )( 4128.0 , 0.0 )( 5337.0 , 0.0 )( 7435.0 , 0.0 )( 9825.0 , 2.31980989973124e-09 )( 14077.0 , 0.00028002607370014587 )( 18849.0 , 1.0028210799983839e-05 )( 27431.0 , 1.0066435200073443e-05 )( 37243.0 , 0.00011012640000007679 )( 53915.0 , 0.0054748306400000946 )( 73399.0 , 0.00285009148149995 )( 106850.0 , 0.0003014474900000508 )( 146007.0 , 0.010482840000000104 )( 212464.0 , 0.014073789999999864 )( 290753.0 , 0.07561010000000001 )};
\addlegendentry{$\mathrm{osc}_{\max}(u_h), \varepsilon_{\mathrm{thresh}}=10^{-6}$} 
\addplot[color=gold,  mark=pentagon*, line width = 0.5mm, dashdotted,, mark options = {scale= 1.5, solid}] 
coordinates{( 25.0 , 0.0 )( 81.0 , 0.0 )( 289.0 , 0.0 )( 1089.0 , 6.17239592770602e-12 )( 1302.0 , 0.0 )( 1701.0 , 0.0 )( 2507.0 , 0.0 )( 4128.0 , 0.0 )( 5337.0 , 0.0 )( 7435.0 , 0.0 )( 9731.0 , 0.0 )( 14087.0 , 0.0 )( 18964.0 , 0.0 )( 27372.0 , 0.0 )( 37174.0 , 0.0 )( 53957.0 , 9.786615962070755e-12 )( 73739.0 , 1.5473700010204539e-07 )( 106951.0 , 0.0 )( 146421.0 , 0.0 )( 212466.0 , 0.0 )( 291307.0 , 1.0680345496894006e-13 )};
\addlegendentry{$\mathrm{osc}_{\max}(u_h), \varepsilon_{\mathrm{thresh}}=10^{-8}$} 
\addplot[color=magenta,  mark=star, line width = 0.5mm, dashdotted,, mark options = {scale= 1.5, solid}] 
coordinates{( 25.0 , 0.0 )( 81.0 , 0.0 )( 289.0 , 0.0 )( 1089.0 , 0.0 )( 1302.0 , 0.0 )( 1701.0 , 0.0 )( 2507.0 , 0.0 )( 4128.0 , 0.0 )( 5337.0 , 0.0 )( 7435.0 , 0.0 )( 9733.0 , 0.0 )( 14090.0 , 9.985567928083583e-12 )( 18961.0 , 0.0 )( 27369.0 , 0.0 )( 37169.0 , 1.2968737195251379e-11 )( 53950.0 , 0.0 )( 73406.0 , 0.0 )( 106897.0 , 3.707700813038173e-12 )( 145753.0 , 0.0 )( 212612.0 , 0.0 )( 291208.0 , 0.0 )};
\addlegendentry{$\mathrm{osc}_{\max}(u_h), \varepsilon_{\mathrm{thresh}}=10^{-10}$} 
\end{semilogxaxis}
\end{tikzpicture}\hspace*{1em}
\begin{tikzpicture}[scale=0.6]
\begin{semilogxaxis}[
    legend pos=north west, xlabel = $\#\ \mathrm{dof}$,
    legend cell align ={left},, title = hanging nodes,
    ymin = -0.01, ymax = 0.25,
    legend style={nodes={scale=0.75, transform shape}, }]
\addplot[color=black,  mark=diamond*, line width = 0.5mm, dashdotted,, mark options = {scale= 1.5, solid}] 
coordinates{( 25.0 , 1.0000000000065512e-05 )( 81.0 , 0.0 )( 289.0 , 0.0 )( 1089.0 , 2.719179992283216e-09 )( 1301.0 , 7.73070496506989e-11 )( 1692.0 , 0.0010100002980799605 )( 2479.0 , 4.266031972122164e-11 )( 4040.0 , 0.0008100005012399158 )( 5189.0 , 0.002780009398529959 )( 7264.0 , 0.00021000000000004349 )( 9602.0 , 6.022626841684087e-11 )( 13824.0 , 0.0001706440860000935 )( 18653.0 , 1.8363200005033775e-07 )( 27099.0 , 6.042628199987021e-05 )( 37166.0 , 1.023998400007109e-05 )( 53771.0 , 2.03189999936626e-07 )( 72903.0 , 0.00022054247600000743 )( 106709.0 , 0.0003749004999999972 )( 146098.0 , 0.10732492949999983 )( 212260.0 , 0.008960929999999978 )( 291053.0 , 0.12214309999999995 )};
\addlegendentry{$\mathrm{osc}_{\max}(u_h), \varepsilon_{\mathrm{thresh}}=10^{-6}$} 
\addplot[color=gold,  mark=pentagon*, line width = 0.5mm, dashdotted,, mark options = {scale= 1.5, solid}] 
coordinates{( 25.0 , 0.0 )( 81.0 , 0.0 )( 289.0 , 0.0 )( 1089.0 , 6.17239592770602e-12 )( 1301.0 , 0.0 )( 1692.0 , 0.0 )( 2479.0 , 0.0 )( 4040.0 , 0.0 )( 5189.0 , 0.0 )( 7264.0 , 1.999999999990898e-05 )( 9598.0 , 0.0 )( 13811.0 , 0.0 )( 18650.0 , 0.0 )( 27133.0 , 4.747739978938625e-10 )( 37020.0 , 0.0 )( 53770.0 , 0.0 )( 72936.0 , 1.4858700936315472e-09 )( 106735.0 , 0.0 )( 146401.0 , 0.0 )( 212494.0 , 6.167109933841175e-09 )( 291858.0 , 3.941291737419306e-12 )};
\addlegendentry{$\mathrm{osc}_{\max}(u_h), \varepsilon_{\mathrm{thresh}}=10^{-8}$} 
\addplot[color=magenta,  mark=star, line width = 0.5mm, dashdotted,, mark options = {scale= 1.5, solid}] 
coordinates{( 25.0 , 0.0 )( 81.0 , 0.0 )( 289.0 , 0.0 )( 1089.0 , 0.0 )( 1301.0 , 0.0 )( 1692.0 , 0.0 )( 2479.0 , 0.0 )( 4040.0 , 0.0 )( 5189.0 , 0.0 )( 7264.0 , 0.0 )( 9598.0 , 0.0 )( 13811.0 , 1.3593126624300567e-11 )( 18650.0 , 0.0 )( 27133.0 , 1.2987388942065081e-12 )( 37016.0 , 3.878009025015672e-12 )( 53781.0 , 0.0 )( 72922.0 , 0.0 )( 106720.0 , 0.0 )( 146592.0 , 6.661338147750939e-15 )( 212665.0 , 0.0 )( 291947.0 , 0.0 )};
\addlegendentry{$\mathrm{osc}_{\max}(u_h), \varepsilon_{\mathrm{thresh}}=10^{-10}$} 
\end{semilogxaxis}
\end{tikzpicture}}
\centerline{
\begin{tikzpicture}[scale=0.6]
\begin{semilogxaxis}[
    legend pos=north west, xlabel = $\#\ \mathrm{dof}$,
    legend cell align ={left},, title = conforming closure,
    ymin = -0.01, ymax = 0.2,
    legend style={nodes={scale=0.75, transform shape}, }]
\addplot[color=black,  mark=diamond*, line width = 0.5mm, dashdotted,, mark options = {scale= 1.5, solid}] 
coordinates{( 25.0 , 0.0 )( 81.0 , 0.0 )( 289.0 , 0.0 )( 1089.0 , 1.397239999256783e-07 )( 1297.0 , 3.21494000088407e-07 )( 1698.0 , 3.081819999195545e-08 )( 2506.0 , 9.192580008310358e-08 )( 4161.0 , 2.939370000065722e-06 )( 5375.0 , 2.1524700000163932e-05 )( 7541.0 , 4.709689999993216e-05 )( 9839.0 , 0.00013908630000014632 )( 14229.0 , 0.00011496420000001173 )( 18627.0 , 0.000273480699999995 )( 27679.0 , 0.0007569250000001304 )( 37325.0 , 0.0003147599999999695 )( 54127.0 , 0.0004175814000000333 )( 74021.0 , 0.001830492999999933 )( 107520.0 , 0.0020227299999999726 )( 147158.0 , 0.026892499999999986 )( 213695.0 , 0.04589759999999998 )( 292544.0 , 0.1733465999999999 )};
\addlegendentry{$\mathrm{osc}_{\max}(u_h), \varepsilon_{\mathrm{thresh}}=10^{-6}$} 
\addplot[color=gold,  mark=pentagon*, line width = 0.5mm, dashdotted,, mark options = {scale= 1.5, solid}] 
coordinates{( 25.0 , 0.0 )( 81.0 , 0.0 )( 289.0 , 0.0 )( 1089.0 , 0.0 )( 1297.0 , 0.0 )( 1698.0 , 0.0 )( 2506.0 , 0.0 )( 4161.0 , 0.0 )( 5375.0 , 0.0 )( 7541.0 , 0.0 )( 9845.0 , 0.0 )( 14232.0 , 0.0 )( 19283.0 , 0.0 )( 27668.0 , 0.0 )( 37473.0 , 0.0 )( 54111.0 , 0.0 )( 74202.0 , 0.0 )( 107320.0 , 4.5098391687758976e-10 )( 146360.0 , 5.405009773085112e-12 )( 213142.0 , 1.0620400114902395e-09 )( 291630.0 , 2.1639801062178776e-11 )};
\addlegendentry{$\mathrm{osc}_{\max}(u_h), \varepsilon_{\mathrm{thresh}}=10^{-8}$} 
\addplot[color=magenta,  mark=star, line width = 0.5mm, dashdotted,, mark options = {scale= 1.5, solid}] 
coordinates{( 25.0 , 0.0 )( 81.0 , 0.0 )( 289.0 , 0.0 )( 1089.0 , 0.0 )( 1297.0 , 0.0 )( 1698.0 , 0.0 )( 2506.0 , 0.0 )( 4161.0 , 0.0 )( 5375.0 , 0.0 )( 7541.0 , 0.0 )( 9845.0 , 0.0 )( 14232.0 , 0.0 )( 19283.0 , 0.0 )( 27668.0 , 0.0 )( 37471.0 , 0.0 )( 54111.0 , 0.0 )( 74202.0 , 0.0 )( 107320.0 , 0.0 )( 146372.0 , 0.0 )( 213139.0 , 0.0 )( 291618.0 , 0.0 )};
\addlegendentry{$\mathrm{osc}_{\max}(u_h), \varepsilon_{\mathrm{thresh}}=10^{-10}$} 
\end{semilogxaxis}
\end{tikzpicture}\hspace*{1em}
\begin{tikzpicture}[scale=0.6]
\begin{semilogxaxis}[
    legend pos=north west, xlabel = $\#\ \mathrm{dof}$,
    legend cell align ={left},, title = hanging nodes,
    ymin = -0.01, ymax = 0.2,
    legend style={nodes={scale=0.75, transform shape}, }]
\addplot[color=black,  mark=diamond*, line width = 0.5mm, dashdotted,, mark options = {scale= 1.5, solid}] 
coordinates{( 25.0 , 0.0 )( 81.0 , 0.0 )( 289.0 , 0.0 )( 1089.0 , 1.397239999256783e-07 )( 1297.0 , 1.5644299189432331e-09 )( 1690.0 , 1.7032399934890918e-09 )( 2473.0 , 1.3236399998284298e-07 )( 4049.0 , 4.898509999851086e-07 )( 5174.0 , 1.9631699998967633e-06 )( 7294.0 , 6.063520000099132e-06 )( 9644.0 , 4.4912999999979775e-05 )( 13979.0 , 3.447629999997481e-05 )( 18753.0 , 7.155739999986643e-05 )( 27405.0 , 0.0001484175999999504 )( 36922.0 , 0.0002592429999999091 )( 54162.0 , 0.0002526060000000996 )( 73167.0 , 0.0013317899999998328 )( 107351.0 , 0.0017151600000000489 )( 147071.0 , 0.0023278530000001574 )( 213363.0 , 0.03074359999999987 )( 292275.0 , 0.05070169999999985 )};
\addlegendentry{$\mathrm{osc}_{\max}(u_h), \varepsilon_{\mathrm{thresh}}=10^{-6}$} 
\addplot[color=gold,  mark=pentagon*, line width = 0.5mm, dashdotted,, mark options = {scale= 1.5, solid}] 
coordinates{( 25.0 , 0.0 )( 81.0 , 0.0 )( 289.0 , 0.0 )( 1089.0 , 0.0 )( 1297.0 , 0.0 )( 1690.0 , 0.0 )( 2473.0 , 0.0 )( 4049.0 , 0.0 )( 5174.0 , 0.0 )( 7294.0 , 0.0 )( 9646.0 , 0.0 )( 13979.0 , 0.0 )( 18753.0 , 0.0 )( 27406.0 , 0.0 )( 36929.0 , 0.0 )( 54145.0 , 0.0 )( 72805.0 , 0.0 )( 107270.0 , 4.219469218469385e-11 )( 147056.0 , 0.0 )( 213220.0 , 7.492451103985331e-12 )( 292315.0 , 1.858970000157001e-08 )};
\addlegendentry{$\mathrm{osc}_{\max}(u_h), \varepsilon_{\mathrm{thresh}}=10^{-8}$} 
\addplot[color=magenta,  mark=star, line width = 0.5mm, dashdotted,, mark options = {scale= 1.5, solid}] 
coordinates{( 25.0 , 0.0 )( 81.0 , 0.0 )( 289.0 , 0.0 )( 1089.0 , 0.0 )( 1297.0 , 0.0 )( 1690.0 , 0.0 )( 2473.0 , 0.0 )( 4049.0 , 0.0 )( 5174.0 , 0.0 )( 7294.0 , 0.0 )( 9646.0 , 0.0 )( 13979.0 , 0.0 )( 18753.0 , 0.0 )( 27406.0 , 0.0 )( 36929.0 , 0.0 )( 54143.0 , 0.0 )( 72799.0 , 0.0 )( 107270.0 , 0.0 )( 147056.0 , 0.0 )( 213183.0 , 0.0 )( 291522.0 , 0.0 )};
\addlegendentry{$\mathrm{osc}_{\max}(u_h), \varepsilon_{\mathrm{thresh}}=10^{-10}$} 
\end{semilogxaxis}
\end{tikzpicture}}
\centerline{
\begin{tikzpicture}[scale=0.6]
\begin{semilogxaxis}[
    legend pos=north west, xlabel = $\#\ \mathrm{dof}$,
    legend cell align ={left},, title = conforming closure,
    ymin = -0.01, ymax = 0.19,
    legend style={nodes={scale=0.75, transform shape}, }]
\addplot[color=black,  mark=diamond*, line width = 0.5mm, dashdotted,, mark options = {scale= 1.5, solid}] 
coordinates{( 25.0 , 0.0 )( 81.0 , 0.0 )( 289.0 , 0.0 )( 1089.0 , 1.4367400003045816e-07 )( 1301.0 , 4.2663006460941233e-10 )( 1710.0 , 3.246059998573969e-08 )( 2553.0 , 2.28254000000927e-05 )( 3229.0 , 6.400029999964474e-06 )( 5446.0 , 4.4543300000032815e-05 )( 7790.0 , 3.8297599999959075e-05 )( 10246.0 , 1.1068400000091572e-05 )( 14673.0 , 4.804280000003658e-05 )( 19050.0 , 0.0001934329999999651 )( 37068.0 , 0.0012694740000001037 )( 55702.0 , 0.0022295800000000643 )( 75471.0 , 0.00258402000000002 )( 109740.0 , 0.03882129999999995 )( 148946.0 , 0.04444849999999989 )( 215616.0 , 0.040556000000000036 )( 294308.0 , 0.18382149999999986 )};
\addlegendentry{$\mathrm{osc}_{\max}(u_h), \varepsilon_{\mathrm{thresh}}=10^{-6}$} 
\addplot[color=gold,  mark=pentagon*, line width = 0.5mm, dashdotted,, mark options = {scale= 1.5, solid}] 
coordinates{( 25.0 , 0.0 )( 81.0 , 0.0 )( 289.0 , 0.0 )( 1089.0 , 0.0 )( 1301.0 , 0.0 )( 1710.0 , 0.0 )( 2553.0 , 0.0 )( 3231.0 , 0.0 )( 5448.0 , 0.0 )( 7790.0 , 0.0 )( 10246.0 , 0.0 )( 14678.0 , 0.0 )( 19790.0 , 0.0 )( 28407.0 , 0.0 )( 38512.0 , 0.0 )( 55776.0 , 0.0 )( 75659.0 , 0.0 )( 110605.0 , 3.330689057889913e-10 )( 150548.0 , 8.737510714951213e-10 )( 219344.0 , 8.401350726217061e-10 )( 298608.0 , 2.9316199912088337e-08 )};
\addlegendentry{$\mathrm{osc}_{\max}(u_h), \varepsilon_{\mathrm{thresh}}=10^{-8}$} 
\addplot[color=magenta,  mark=star, line width = 0.5mm, dashdotted,, mark options = {scale= 1.5, solid}] 
coordinates{( 25.0 , 0.0 )( 81.0 , 0.0 )( 289.0 , 0.0 )( 1089.0 , 0.0 )( 1301.0 , 0.0 )( 1710.0 , 0.0 )( 2553.0 , 0.0 )( 3231.0 , 0.0 )( 5448.0 , 0.0 )( 7790.0 , 0.0 )( 10246.0 , 0.0 )( 14678.0 , 0.0 )( 19790.0 , 0.0 )( 28407.0 , 0.0 )( 38509.0 , 0.0 )( 55764.0 , 0.0 )( 75659.0 , 0.0 )( 110555.0 , 0.0 )( 150244.0 , 0.0 )( 219212.0 , 0.0 )( 297237.0 , 0.0 )};
\addlegendentry{$\mathrm{osc}_{\max}(u_h), \varepsilon_{\mathrm{thresh}}=10^{-10}$} 
\end{semilogxaxis}
\end{tikzpicture}\hspace*{1em}
\begin{tikzpicture}[scale=0.6]
\begin{semilogxaxis}[scaled y ticks = false,
    legend pos=north west, xlabel = $\#\ \mathrm{dof}$,
    legend cell align ={left},, title = hanging nodes,
    ymin = -0.01, ymax = 0.05,
    legend style={nodes={scale=0.75, transform shape}, }]
\addplot[color=black,  mark=diamond*, line width = 0.5mm, dashdotted,, mark options = {scale= 1.5, solid}] 
coordinates{( 25.0 , 0.0 )( 81.0 , 0.0 )( 289.0 , 0.0 )( 1089.0 , 1.4367400003045816e-07 )( 1300.0 , 1.1679400002684304e-07 )( 1713.0 , 2.9239299998984336e-06 )( 2543.0 , 1.412110000043043e-06 )( 3204.0 , 6.29712000099758e-07 )( 5169.0 , 4.1658700000724025e-06 )( 7732.0 , 7.977289999683279e-07 )( 10217.0 , 1.7179270000111657e-05 )( 14647.0 , 7.952519999987473e-05 )( 19201.0 , 0.0004457570000000466 )( 37341.0 , 0.0003436359999999805 )( 73259.0 , 0.0022881100000000654 )( 145048.0 , 0.026917700000000044 )( 216521.0 , 0.042079999999999895 )( 295712.0 , 0.04710250000000005 )};
\addlegendentry{$\mathrm{osc}_{\max}(u_h), \varepsilon_{\mathrm{thresh}}=10^{-6}$} 
\addplot[color=gold,  mark=pentagon*, line width = 0.5mm, dashdotted,, mark options = {scale= 1.5, solid}] 
coordinates{( 25.0 , 0.0 )( 81.0 , 0.0 )( 289.0 , 0.0 )( 1089.0 , 0.0 )( 1300.0 , 0.0 )( 1715.0 , 0.0 )( 2542.0 , 0.0 )( 3206.0 , 0.0 )( 4281.0 , 0.0 )( 5617.0 , 0.0 )( 10201.0 , 0.0 )( 14681.0 , 0.0 )( 19681.0 , 0.0 )( 37911.0 , 7.466471885209103e-12 )( 55475.0 , 0.0 )( 75568.0 , 0.0 )( 109926.0 , 3.663649383867096e-10 )( 150484.0 , 3.4861002973229915e-14 )( 218643.0 , 1.1787899545367964e-09 )( 297004.0 , 1.653460035555554e-09 )};
\addlegendentry{$\mathrm{osc}_{\max}(u_h), \varepsilon_{\mathrm{thresh}}=10^{-8}$} 
\addplot[color=magenta,  mark=star, line width = 0.5mm, dashdotted,, mark options = {scale= 1.5, solid}] 
coordinates{( 25.0 , 0.0 )( 81.0 , 0.0 )( 289.0 , 0.0 )( 1089.0 , 0.0 )( 1300.0 , 0.0 )( 1715.0 , 0.0 )( 2542.0 , 0.0 )( 3206.0 , 0.0 )( 4281.0 , 0.0 )( 5617.0 , 0.0 )( 10201.0 , 0.0 )( 14681.0 , 0.0 )( 19681.0 , 0.0 )( 37908.0 , 0.0 )( 55475.0 , 0.0 )( 75556.0 , 0.0 )( 109948.0 , 0.0 )( 150447.0 , 0.0 )( 218442.0 , 0.0 )( 297418.0 , 0.0 )};
\addlegendentry{$\mathrm{osc}_{\max}(u_h), \varepsilon_{\mathrm{thresh}}=10^{-10}$} 
\end{semilogxaxis}
\end{tikzpicture}}
\caption{Example~\ref{ex:hmm86}: Dependency of the spurious oscillations on the stopping 
criterion in the solver for 
%scheme for solving 
the nonlinear problem: AFC scheme with BJK limiter (top), with Kuzmin limiter (middle), and MUAS method (bottom).}
\label{fig:var_hmm86_differentl}
\end{figure}
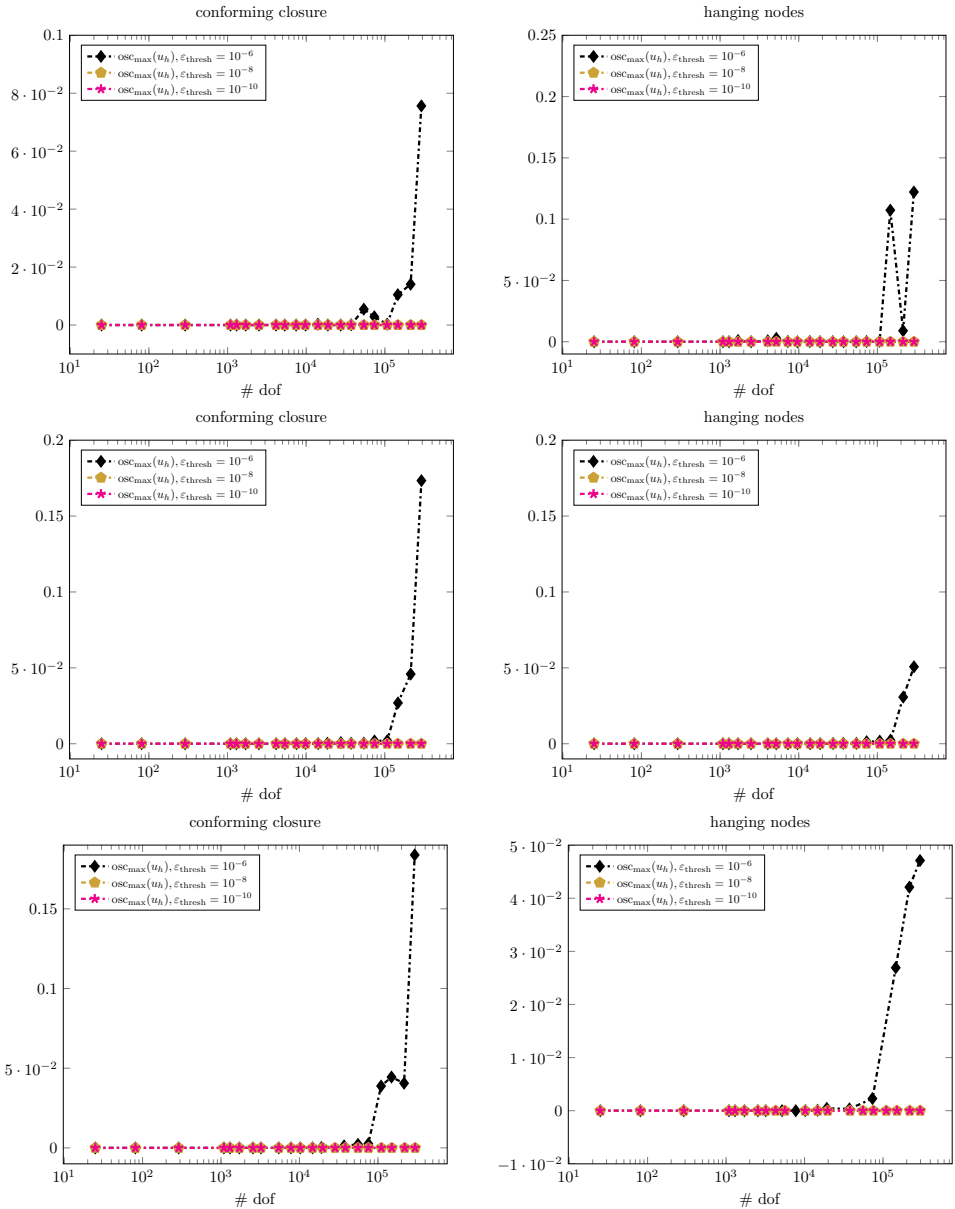

Figure~\ref{fig:var_hmm86_differentl} provides information on the impact of the weaker stopping criteria
on the satisfaction of the global DMP. Only for the weakest stopping criterion  $\varepsilon_{\mathrm{thresh}} = 10^{-6}$ 
and on fine grids there are notable spurious oscillations.

Concerning the width of the interior layer, we usually could not observe visible differences 
between the results from Figure~\ref{fig:thickness_hmm86_grid_1} and the results for  $\varepsilon_{\mathrm{thresh}} = 10^{-8}$.
Often, also the layer width of the solutions computed with $\varepsilon_{\mathrm{thresh}} = 10^{-6}$ 
is similar. Only on very fine grids, we could see more smearing with this stopping criterion.
For the sake of brevity, the results with respect to the layer width are not presented in detail. 

\begin{remark} Continuing the adaptive refinement in this example creates very small mesh cells. We
could observe that the sparse direct solver failed, giving \texttt{nan}, if cells with a diameter 
of around $10^{-6}$ occurred. In contrast, a standard iterative solver, GMRES with SSOR preconditioner, still 
worked well in this situation.
\hfill$\Box$\end{remark}

\subsection{Hemker Problem}\label{ex:hemker}
\begin{figure}[t!]
\centerline{\includegraphics[width=0.5\textwidth]{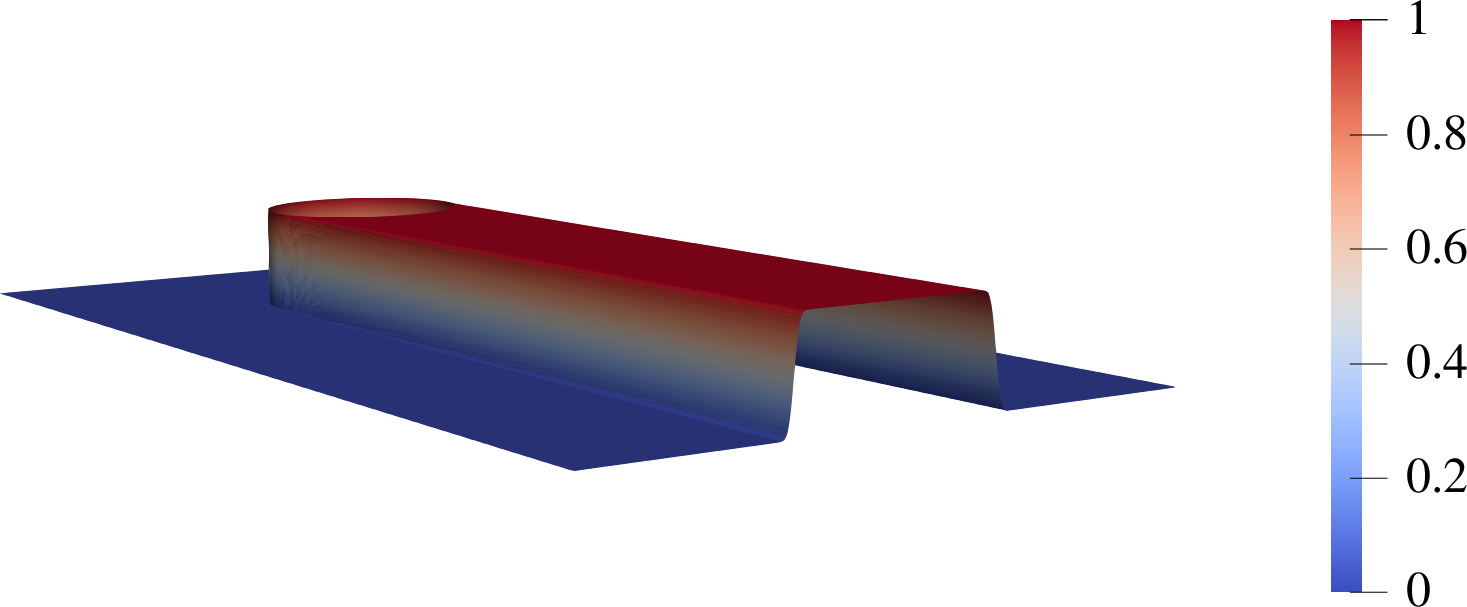}
\hspace*{1em} \includegraphics[width=0.4\textwidth]{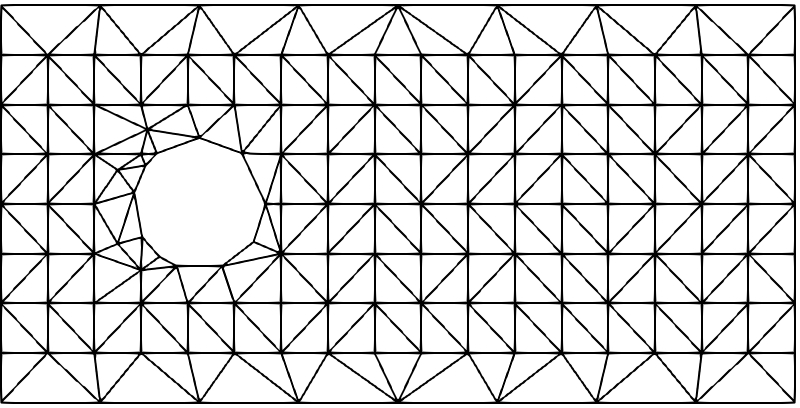}}
\caption{Example~\ref{ex:hemker}: Solution for $\varepsilon = 10^{-4}$ (left), computed with the BJK limiter, level~6; initial grid (right), level~0.}
\label{fig:sol_hemker2d}
\end{figure}

The Hemker problem is a standard benchmark problem defined in \cite{Hem96}. The domain  is given by
$
\Omega = \{(-3,9)\times(-3,3)\} \setminus \{(x,y) \ : \ x^2+y^2\le 1 \}
$,
the convection field by 
$\bb = (1,0)^T$, and the reaction field and right-hand side in Eq.~\eqref{eq:cdr_eqn} vanish:
$c=f=0$. Dirichlet boundary conditions are set
at $x=-3$, with $u_b=0$ and at the circular boundary with $u_b=1$. On all other 
boundaries, homogeneous Neumann conditions are prescribed.
This problem was studied comprehensively for $\varepsilon = 10^{-4}$ in
\cite{ACJ11} and  reference values are 
available for some quantities of interest. This diffusion parameter was used also in our studies, see 
Figure~\ref{fig:sol_hemker2d} for an illustration of the solution, which takes values in $[0,1]$.

Figure~\ref{fig:sol_hemker2d} presents the initial grid with 
$\#\mathrm{dof}=151$.  The adaptive refinement was started after having computed the solution on the initial grid. It 
was stopped when $\#\mathrm{dof}\gtrsim 5\times 10^5$. During refinement, the approximation of the circular
boundary was improved. Based on the experience from the previous example, the threshold for stopping the 
iterative solution of the nonlinear problem was set to be 
$\varepsilon_{\mathrm{thresh}} = 10^{-8}$.

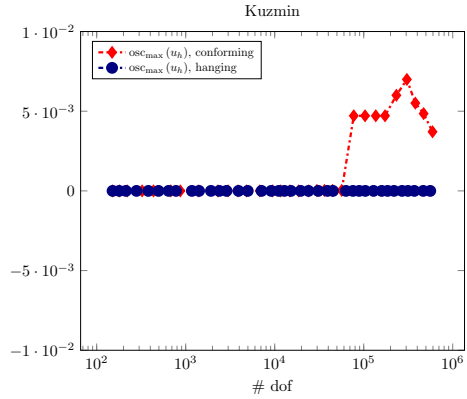
\begin{figure}[t!]
\centerline{
\begin{tikzpicture}[scale=0.6]
\begin{semilogxaxis}[scaled y ticks = false,
    legend pos=north west, xlabel = $\#\ \mathrm{dof}$,
    legend cell align ={left},, title = Kuzmin,
    ymin = -0.01, ymax = 0.01,
    legend style={nodes={scale=0.75, transform shape}, }]
\addplot[color=red,  mark=diamond*, line width = 0.5mm, dashdotted,, mark options = {scale= 1.5, solid}] 
coordinates{( 151.0 , 6.061817714453355e-14 )( 180.0 , 1.1857181902996672e-13 )( 220.0 , 2.220446049250313e-16 )( 324.0 , 1.9984014443252818e-15 )( 435.0 , 1.1102230246251565e-15 )( 681.0 , 0.0 )( 871.0 , 0.0 )( 1402.0 , 0.0 )( 2299.0 , 0.0 )( 2997.0 , 0.0 )( 3872.0 , 0.0 )( 4833.0 , 0.0 )( 6798.0 , 0.0 )( 9437.0 , 0.0 )( 11688.0 , 0.0 )( 14830.0 , 0.0 )( 18577.0 , 0.0 )( 23974.0 , 0.0 )( 29852.0 , 2.9913899999955973e-05 )( 36191.0 , 2.9913899999955973e-05 )( 44648.0 , 2.9913899999955973e-05 )( 56593.0 , 2.9913899999955973e-05 )( 77334.0 , 0.004712529999999937 )( 103535.0 , 0.004712529999999937 )( 136608.0 , 0.004712529999999937 )( 173635.0 , 0.004712529999999937 )( 232836.0 , 0.005999850000000029 )( 305545.0 , 0.006993870000000069 )( 380436.0 , 0.005506470000000041 )( 471317.0 , 0.004848310000000078 )( 593075.0 , 0.003710719999999945 )};
\addlegendentry{$\mathrm{osc}_{\mathrm{max}}\left(u_h\right)$, conforming} 
\addplot[color=navy_blue,  mark=oplus*, line width = 0.5mm, dashdotted,, mark options = {scale= 1.5, solid}] 
coordinates{( 151.0 , 6.061817714453355e-14 )( 179.0 , 2.361222328772783e-12 )( 213.0 , 8.426592756904938e-13 )( 280.0 , 0.0 )( 380.0 , 0.0 )( 495.0 , 0.0 )( 642.0 , 0.0 )( 778.0 , 0.0 )( 1174.0 , 0.0 )( 1412.0 , 0.0 )( 1923.0 , 0.0 )( 2388.0 , 0.0 )( 2873.0 , 0.0 )( 3922.0 , 0.0 )( 4977.0 , 0.0 )( 7124.0 , 0.0 )( 9128.0 , 0.0 )( 11099.0 , 0.0 )( 12860.0 , 0.0 )( 15228.0 , 0.0 )( 19571.0 , 0.0 )( 23982.0 , 0.0 )( 30866.0 , 0.0 )( 39171.0 , 0.0 )( 45141.0 , 0.0 )( 63237.0 , 0.0 )( 75748.0 , 0.0 )( 88913.0 , 0.0 )( 104176.0 , 0.0 )( 128866.0 , 0.0 )( 159279.0 , 0.0 )( 185642.0 , 0.0 )( 219126.0 , 0.0 )( 269617.0 , 0.0 )( 314219.0 , 0.0 )( 370647.0 , 0.0 )( 466214.0 , 0.0 )( 558852.0 , 0.0 )};
\addlegendentry{$\mathrm{osc}_{\mathrm{max}}\left(u_h\right)$, hanging} 
\end{semilogxaxis}
\end{tikzpicture}}
\caption{Example~\ref{ex:hemker}: Spurious oscillations for different grids with conforming closure and hanging nodes for the AFC scheme with Kuzmin limiter. 
There are no spurious oscillations for the solutions computed with the two other schemes.}
\label{fig:var_hemker_different}
\end{figure}

The satisfaction of the global DMP was measured again by $\mathrm{osc}_{\mathrm{max}}(u_h)$
defined in \eqref{rem:var_def}. As in the previous example, for the AFC scheme with BJK limiter and the MUAS method, only 
unphysical values of the order of the stopping criterion for solving the nonlinear problems could be 
observed. Hence, these methods satisfy the global DMP. In contrast, there are small but notable spurious
oscillations for the AFC scheme with Kuzmin limiter on fine conforming grids, compare Figure~\ref{fig:var_hemker_different}.
We think that the reason is the appearance of non-Delaunay closure cells in combination with the fact that the 
discrete problem becomes locally diffusion-dominated in strongly refined regions. 

\begin{figure}[t!]
\centerline{
\begin{tikzpicture}[scale=0.6]
\begin{semilogxaxis}[
    legend pos=north east, xlabel = $\#\ \mathrm{dof}$,ylabel = $\mathrm{smear}_{\mathrm{int}}$,
    legend cell align ={left},xmin=500, xmax= 650000,ymin = 0.02, ymax = 0.6, title = conforming closure,
   legend style={nodes={scale=0.75, transform shape}}]
\addplot[color=red,  mark=square*, line width = 0.5mm,dashdotted,, mark options = {scale= 1.5, solid}] 
coordinates{( 180.0 , 1.22588 )( 220.0 , 1.19307 )( 306.0 , inf )( 429.0 , 0.97238 )( 627.0 , 0.605724 )( 942.0 , 0.570385 )( 1205.0 , 0.60169 )( 1675.0 , 0.463656 )( 2143.0 , 0.322207 )( 3038.0 , 0.329577 )( 4129.0 , 0.31345 )( 6107.0 , 0.303 )( 7674.0 , 0.175465 )( 10512.0 , 0.172546 )( 12814.0 , 0.160637 )( 15202.0 , 0.157707 )( 20405.0 , 0.141938 )( 24476.0 , 0.102078 )( 29150.0 , 0.0954961 )( 37127.0 , 0.0931501 )( 44503.0 , 0.0797588 )( 54415.0 , 0.0766658 )( 71993.0 , 0.0755837 )( 92520.0 , 0.0750083 )( 115006.0 , 0.0744322 )( 140054.0 , 0.0751272 )( 170485.0 , 0.0750016 )( 209713.0 , 0.0746181 )( 270877.0 , 0.0742318 )( 336760.0 , 0.0742325 )( 447822.0 , 0.0731091 )( 564487.0 , 0.0730544 )};
\addlegendentry{BJK} 
\addplot[color=navy_blue,  mark=oplus*, line width = 0.5mm,dashdotted,, mark options = {scale= 1.5, solid}] 
coordinates{( 180.0 , inf )( 220.0 , inf )( 324.0 , inf )( 435.0 , 1.45736 )( 681.0 , 1.08377 )( 871.0 , 0.83538 )( 1402.0 , 0.611372 )( 2299.0 , 0.540408 )( 2997.0 , 0.449022 )( 3872.0 , 0.391676 )( 4833.0 , 0.362138 )( 6798.0 , 0.302258 )( 9437.0 , 0.217042 )( 11688.0 , 0.198294 )( 14830.0 , 0.19036 )( 18577.0 , 0.175051 )( 23974.0 , 0.125985 )( 29852.0 , 0.11421 )( 36191.0 , 0.108763 )( 44648.0 , 0.108747 )( 56593.0 , 0.0927532 )( 77334.0 , 0.0827715 )( 103535.0 , 0.0794993 )( 136608.0 , 0.0773821 )( 173635.0 , 0.0770665 )( 232836.0 , 0.0776442 )( 305545.0 , 0.0758397 )( 380436.0 , 0.0762465 )( 471317.0 , 0.0741041 )( 593075.0 , 0.0745353 )};
\addlegendentry{Kuzmin} 
\addplot[color=dark_green,  mark=triangle*, line width = 0.5mm,dashdotted,, mark options = {scale= 1.5, solid}] 
coordinates{( 185.0 , inf )( 230.0 , inf )( 287.0 , inf )( 377.0 , inf )( 523.0 , 1.20554 )( 743.0 , 0.955361 )( 995.0 , 0.898526 )( 1315.0 , 0.760484 )( 1789.0 , 0.549112 )( 2502.0 , 0.550505 )( 3142.0 , 0.512208 )( 4833.0 , 0.44909 )( 6638.0 , 0.338412 )( 8325.0 , 0.311956 )( 11034.0 , 0.283197 )( 16794.0 , 0.202766 )( 20449.0 , 0.196295 )( 25525.0 , 0.187635 )( 32206.0 , 0.179287 )( 42990.0 , 0.165572 )( 55403.0 , 0.145417 )( 68619.0 , 0.11021 )( 90754.0 , 0.0973033 )( 134973.0 , 0.0812435 )( 166271.0 , 0.0778485 )( 232164.0 , 0.0762917 )( 324063.0 , 0.0753247 )( 468358.0 , 0.073809 )( 694708.0 , 0.073745 )};
\addlegendentry{MUAS} 
\addplot[color=black,  line width = 0.25mm, dotted] 
coordinates{( 1.0 , 0.0723 )( 600000, 0.0723 )};
\addlegendentry{reference}
\end{semilogxaxis}
\end{tikzpicture}\hspace*{1em}
\begin{tikzpicture}[scale=0.6]
\begin{semilogxaxis}[
    legend pos=north east, xlabel = $\#\ \mathrm{dof}$,ylabel = $\mathrm{smear}_{\mathrm{int}}$,
    legend cell align ={left},, title = hanging nodes,xmin=500, xmax= 650000,ymin = 0.02, ymax = 0.9, 
   legend style={nodes={scale=0.75, transform shape}}]
\addplot[color=red,  mark=square*, line width = 0.5mm,dashdotted,, mark options = {scale= 1.5, solid}] 
coordinates{( 180.0 , 1.2562 )( 226.0 , 1.18516 )( 278.0 , 1.40899 )( 374.0 , inf )( 556.0 , 1.43846 )( 678.0 , 0.674452 )( 807.0 , 0.617483 )( 1117.0 , 0.613201 )( 1475.0 , 0.586475 )( 1994.0 , 0.578814 )( 2518.0 , 0.534038 )( 2982.0 , 0.542781 )( 3872.0 , 0.306865 )( 4834.0 , 0.302958 )( 5851.0 , 0.287246 )( 7583.0 , 0.190609 )( 8991.0 , 0.14517 )( 10930.0 , 0.139165 )( 14089.0 , 0.16963 )( 19211.0 , 0.155646 )( 23382.0 , 0.133148 )( 26980.0 , 0.142673 )( 34123.0 , 0.129287 )( 40367.0 , 0.101808 )( 46530.0 , 0.0974224 )( 54801.0 , 0.0893923 )( 65219.0 , 0.0848154 )( 75631.0 , 0.0806323 )( 91513.0 , 0.0748113 )( 109441.0 , 0.0749665 )( 126845.0 , 0.0751264 )( 149633.0 , 0.0753319 )( 189801.0 , 0.0745206 )( 218893.0 , 0.07359 )( 257942.0 , 0.0731232 )( 316908.0 , 0.0728889 )( 382326.0 , 0.0727795 )( 445494.0 , 0.0729016 )( 536367.0 , 0.072914 )};
\addlegendentry{BJK} 
\addplot[color=navy_blue,  mark=oplus*, line width = 0.5mm,dashdotted,, mark options = {scale= 1.5, solid}] 
coordinates{( 179.0 , 2.03797 )( 213.0 , inf )( 280.0 , inf )( 380.0 , inf )( 495.0 , inf )( 642.0 , 1.44698 )( 778.0 , 1.25933 )( 1174.0 , 1.00242 )( 1412.0 , 0.950838 )( 1923.0 , 0.695389 )( 2388.0 , 0.648 )( 2873.0 , 0.594404 )( 3922.0 , 0.4843 )( 4977.0 , 0.421517 )( 7124.0 , 0.373737 )( 9128.0 , 0.31432 )( 11099.0 , 0.2353 )( 12860.0 , 0.226301 )( 15228.0 , 0.215355 )( 19571.0 , 0.20257 )( 23982.0 , 0.173705 )( 30866.0 , 0.130657 )( 39171.0 , 0.123002 )( 45141.0 , 0.119412 )( 63237.0 , 0.112555 )( 75748.0 , 0.095214 )( 88913.0 , 0.0863842 )( 104176.0 , 0.0839128 )( 128866.0 , 0.0841345 )( 159279.0 , 0.0822109 )( 185642.0 , 0.0841384 )( 219126.0 , 0.080478 )( 269617.0 , 0.0810177 )( 314219.0 , 0.0752269 )( 370647.0 , 0.0741174 )( 466214.0 , 0.0739017 )( 558852.0 , 0.0739417 )};
\addlegendentry{Kuzmin} 
\addplot[color=dark_green,  mark=triangle*, line width = 0.5mm,dashdotted,, mark options = {scale= 1.5, solid}] 
coordinates{( 183.0 , 2.04201 )( 223.0 , inf )( 279.0 , inf )( 364.0 , inf )( 500.0 , inf )( 638.0 , inf )( 799.0 , 1.87656 )( 942.0 , 1.36934 )( 1211.0 , 1.35393 )( 1416.0 , 1.02848 )( 1658.0 , 0.944558 )( 2274.0 , 0.667732 )( 2969.0 , 0.652978 )( 3614.0 , 0.630744 )( 4903.0 , 0.575977 )( 6191.0 , 0.506537 )( 9087.0 , 0.356509 )( 11118.0 , 0.332958 )( 13401.0 , 0.331062 )( 17823.0 , 0.296584 )( 23202.0 , 0.255661 )( 29496.0 , 0.200765 )( 36529.0 , 0.176618 )( 49974.0 , 0.1346 )( 64742.0 , 0.123906 )( 79739.0 , 0.119729 )( 94504.0 , 0.108188 )( 110910.0 , 0.0987091 )( 128613.0 , 0.085802 )( 158021.0 , 0.0830357 )( 184392.0 , 0.0821496 )( 227287.0 , 0.0793954 )( 297832.0 , 0.078906 )( 363196.0 , 0.0769846 )( 431915.0 , 0.0771061 )( 519291.0 , 0.075085 )};
\addlegendentry{MUAS} 
\addplot[color=black,  line width = 0.25mm, dotted] 
coordinates{( 1.0 , 0.0723 )( 600000, 0.0723 )};
\addlegendentry{reference}
\end{semilogxaxis}
\end{tikzpicture}}
\caption{Example~\ref{ex:hemker}: Thickness of the internal layer at $x=4$ , $\mathrm{smear}_{\mathrm{int}}$.}
\label{fig:smear_hemker_different}
\end{figure}
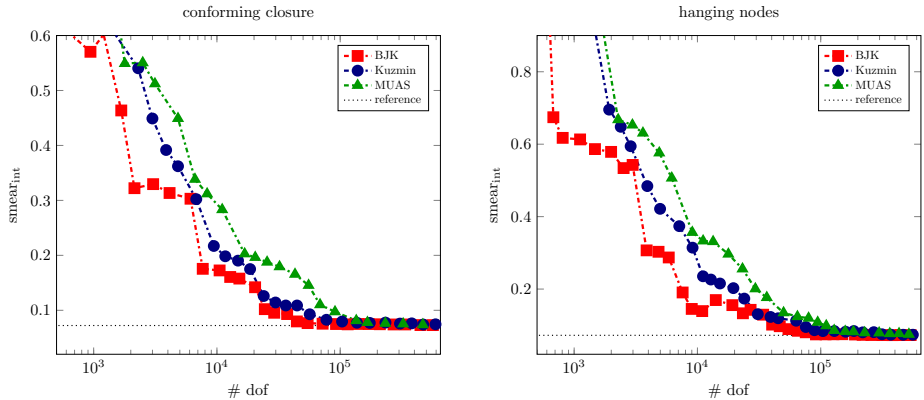

\begin{figure}[t!]
\centerline{\includegraphics[width=0.47\textwidth]{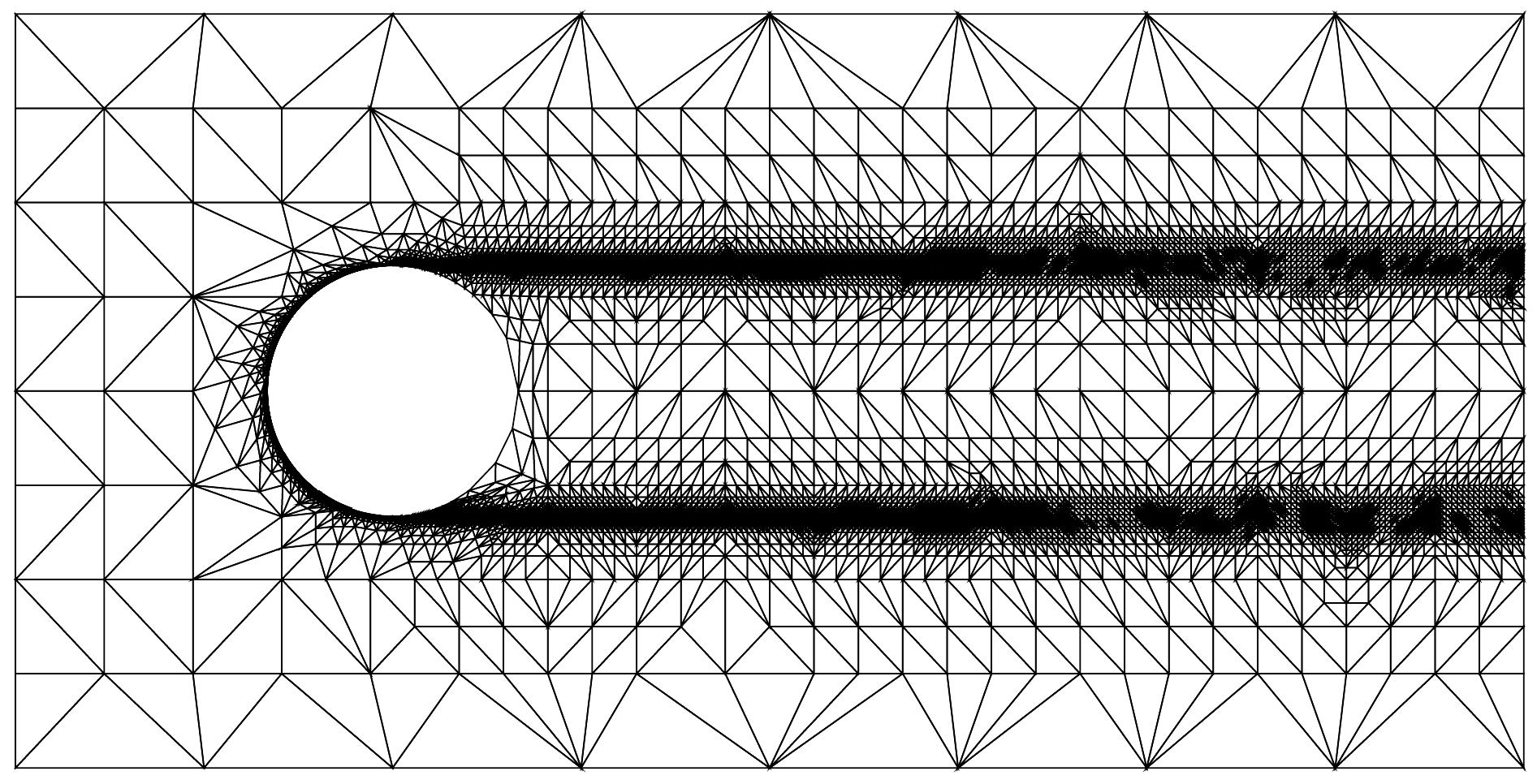}
\hspace*{1em} \includegraphics[width=0.47\textwidth]{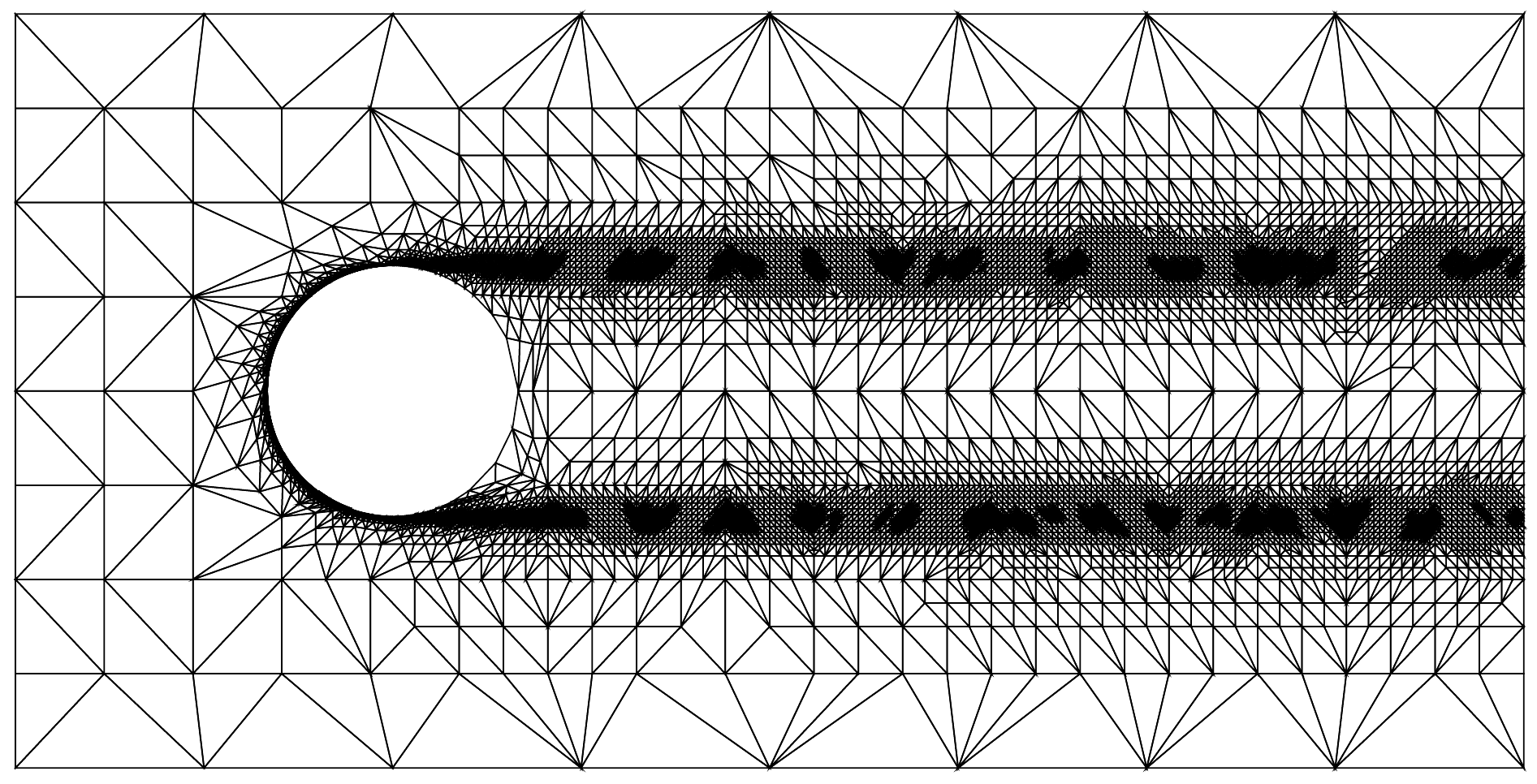}}
\caption{Example~\ref{ex:hemker}: Adaptively refined conforming grids with $\approx 25,000 \ \#\mathrm{dof}$,
left with AFC method and Kuzmin limiter, right with MUAS method.}
\label{fig:hemker_grids}
\end{figure}

For assessing the accuracy of the solutions in \cite{ACJ11}, the width of the internal layer at $y=1$ on 
the cut line at $x=4$ was considered. The definition of the layer width is similar like for the quantity 
$\mathrm{smear}_{\mathrm{int}}$ from \eqref{eq:smear}. In \cite{ACJ11}, the reference value $0.0723$ is 
provided. The results obtained with the considered schemes are presented in 
Figure~\ref{fig:smear_hemker_different}. In general, the sharpest layer was computed with AFC scheme with 
BJK limiter. On sufficiently fine grids, the results for all methods are very close to the reference value.
Up to around 100,000 $\#\mathrm{dof}$, the results for the MUAS method 
are notably less accurate than for the other two methods. The reason is that the 
adaptive grid refinement occurred for this method in a somewhat different way, see Figure~\ref{fig:hemker_grids}
for a representative example. For the AFC methods, the region of this cut line is already much stronger refined. This situation shows that there is the need of an improved mechanism for controlling the 
adaptive grid refinement for the MUAS method, i.e., the need of developing an a posteriori error estimator. 

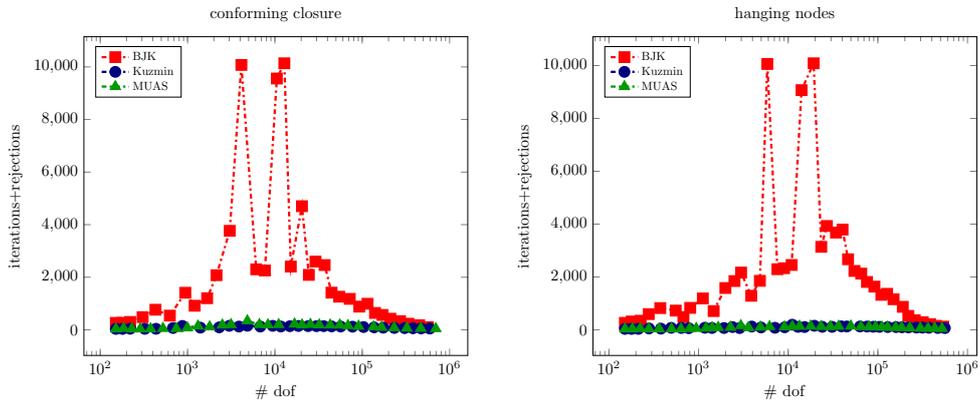
\begin{figure}[t!]
\centerline{
\begin{tikzpicture}[scale=0.6]
\begin{semilogxaxis}[scaled y ticks = false,
    legend pos=north west, xlabel = $\#\ \mathrm{dof}$, ylabel = iterations+rejections,
    legend cell align ={left}, title = {conforming closure},
     legend style={nodes={scale=0.75, transform shape}}]
\addplot[color=red,  mark=square*, line width = 0.5mm, dashdotted,,mark options = {scale= 1.5, solid}]
coordinates{( 151.0 , 264.0 )( 180.0 , 266.0 )( 220.0 , 299.0 )( 306.0 , 482.0 )( 429.0 , 764.0 )( 627.0 , 540.0 )( 942.0 , 1409.0 )( 1205.0 , 916.0 )( 1675.0 , 1202.0 )( 2143.0 , 2077.0 )( 3038.0 , 3765.0 )( 4129.0 , 10074.0 )( 6107.0 , 2300.0 )( 7674.0 , 2256.0 )( 10512.0 , 9558.0 )( 12814.0 , 10137.0 )( 15202.0 , 2406.0 )( 20405.0 , 4702.0 )( 24476.0 , 2089.0 )( 29150.0 , 2596.0 )( 37127.0 , 2463.0 )( 44503.0 , 1414.0 )( 54415.0 , 1259.0 )( 71993.0 , 1177.0 )( 92520.0 , 882.0 )( 115006.0 , 991.0 )( 140054.0 , 641.0 )( 170485.0 , 566.0 )( 209713.0 , 422.0 )( 270877.0 , 331.0 )( 336760.0 , 232.0 )( 447822.0 , 177.0 )( 564487.0 , 102.0 )};
\addlegendentry{BJK} 
\addplot[color=navy_blue,  mark=oplus*, line width = 0.5mm, dashdotted,,mark options = {scale= 1.5, solid}] 
coordinates{( 151.0 , 35.0 )( 180.0 , 35.0 )( 220.0 , 41.0 )( 324.0 , 38.0 )( 435.0 , 39.0 )( 681.0 , 69.0 )( 871.0 , 148.0 )( 1402.0 , 82.0 )( 2299.0 , 93.0 )( 2997.0 , 152.0 )( 3872.0 , 124.0 )( 4833.0 , 155.0 )( 6798.0 , 124.0 )( 9437.0 , 151.0 )( 11688.0 , 126.0 )( 14830.0 , 138.0 )( 18577.0 , 146.0 )( 23974.0 , 141.0 )( 29852.0 , 139.0 )( 36191.0 , 140.0 )( 44648.0 , 130.0 )( 56593.0 , 125.0 )( 77334.0 , 119.0 )( 103535.0 , 106.0 )( 136608.0 , 94.0 )( 173635.0 , 82.0 )( 232836.0 , 71.0 )( 305545.0 , 64.0 )( 380436.0 , 58.0 )( 471317.0 , 54.0 )( 593075.0 , 49.0 )};
\addlegendentry{Kuzmin} 
\addplot[color=dark_green,  mark=triangle*, line width = 0.5mm, dashdotted,,mark options = {scale= 1.5, solid}] 
coordinates{( 151.0 , 35.0 )( 185.0 , 54.0 )( 230.0 , 35.0 )( 287.0 , 36.0 )( 377.0 , 38.0 )( 523.0 , 51.0 )( 743.0 , 48.0 )( 995.0 , 75.0 )( 1315.0 , 142.0 )( 1789.0 , 113.0 )( 2502.0 , 181.0 )( 3142.0 , 197.0 )( 4833.0 , 340.0 )( 6638.0 , 199.0 )( 8325.0 , 190.0 )( 11034.0 , 197.0 )( 16794.0 , 219.0 )( 20449.0 , 187.0 )( 25525.0 , 204.0 )( 32206.0 , 187.0 )( 42990.0 , 183.0 )( 55403.0 , 182.0 )( 68619.0 , 160.0 )( 90754.0 , 137.0 )( 134973.0 , 115.0 )( 166271.0 , 95.0 )( 232164.0 , 78.0 )( 324063.0 , 66.0 )( 468358.0 , 57.0 )( 694708.0 , 45.0 )};
\addlegendentry{MUAS} 
\end{semilogxaxis}
\end{tikzpicture}\hspace*{1em}
\begin{tikzpicture}[scale=0.6]
\begin{semilogxaxis}[scaled y ticks = false,
    legend pos=north west, xlabel = $\#\ \mathrm{dof}$, ylabel = iterations+rejections,
    legend cell align ={left}, title = {hanging nodes},
     legend style={nodes={scale=0.75, transform shape}}]
\addplot[color=red,  mark=square*, line width = 0.5mm, dashdotted,,mark options = {scale= 1.5, solid}]
coordinates{( 151.0 , 264.0 )( 180.0 , 316.0 )( 226.0 , 345.0 )( 278.0 , 594.0 )( 374.0 , 823.0 )( 556.0 , 728.0 )( 678.0 , 484.0 )( 807.0 , 830.0 )( 1117.0 , 1188.0 )( 1475.0 , 710.0 )( 1994.0 , 1582.0 )( 2518.0 , 1846.0 )( 2982.0 , 2167.0 )( 3872.0 , 1288.0 )( 4834.0 , 1853.0 )( 5851.0 , 10060.0 )( 7583.0 , 2294.0 )( 8991.0 , 2331.0 )( 10930.0 , 2454.0 )( 14089.0 , 9070.0 )( 19211.0 , 10086.0 )( 23382.0 , 3140.0 )( 26980.0 , 3929.0 )( 34123.0 , 3679.0 )( 40367.0 , 3793.0 )( 46530.0 , 2670.0 )( 54801.0 , 2228.0 )( 65219.0 , 2127.0 )( 75631.0 , 1811.0 )( 91513.0 , 1634.0 )( 109441.0 , 1324.0 )( 126845.0 , 1363.0 )( 149633.0 , 1156.0 )( 189801.0 , 870.0 )( 218893.0 , 526.0 )( 257942.0 , 364.0 )( 316908.0 , 289.0 )( 382326.0 , 216.0 )( 445494.0 , 151.0 )( 536367.0 , 122.0 )};
\addlegendentry{BJK} 
\addplot[color=navy_blue,  mark=oplus*, line width = 0.5mm, dashdotted,,mark options = {scale= 1.5, solid}] 
coordinates{( 151.0 , 35.0 )( 179.0 , 35.0 )( 213.0 , 36.0 )( 280.0 , 52.0 )( 380.0 , 45.0 )( 495.0 , 76.0 )( 642.0 , 71.0 )( 778.0 , 62.0 )( 1174.0 , 78.0 )( 1412.0 , 69.0 )( 1923.0 , 69.0 )( 2388.0 , 105.0 )( 2873.0 , 66.0 )( 3922.0 , 124.0 )( 4977.0 , 90.0 )( 7124.0 , 78.0 )( 9128.0 , 96.0 )( 11099.0 , 194.0 )( 12860.0 , 112.0 )( 15228.0 , 106.0 )( 19571.0 , 154.0 )( 23982.0 , 114.0 )( 30866.0 , 112.0 )( 39171.0 , 131.0 )( 45141.0 , 127.0 )( 63237.0 , 127.0 )( 75748.0 , 124.0 )( 88913.0 , 115.0 )( 104176.0 , 112.0 )( 128866.0 , 105.0 )( 159279.0 , 104.0 )( 185642.0 , 96.0 )( 219126.0 , 92.0 )( 269617.0 , 85.0 )( 314219.0 , 83.0 )( 370647.0 , 78.0 )( 466214.0 , 70.0 )( 558852.0 , 67.0 )};
\addlegendentry{Kuzmin} 
\addplot[color=dark_green,  mark=triangle*, line width = 0.5mm, dashdotted,,mark options = {scale= 1.5, solid}] 
coordinates{( 151.0 , 35.0 )( 183.0 , 36.0 )( 223.0 , 49.0 )( 279.0 , 40.0 )( 364.0 , 37.0 )( 500.0 , 37.0 )( 638.0 , 36.0 )( 799.0 , 51.0 )( 942.0 , 53.0 )( 1211.0 , 65.0 )( 1416.0 , 77.0 )( 1658.0 , 92.0 )( 2274.0 , 97.0 )( 2969.0 , 158.0 )( 3614.0 , 83.0 )( 4903.0 , 102.0 )( 6191.0 , 122.0 )( 9087.0 , 111.0 )( 11118.0 , 132.0 )( 13401.0 , 133.0 )( 17823.0 , 142.0 )( 23202.0 , 125.0 )( 29496.0 , 121.0 )( 36529.0 , 136.0 )( 49974.0 , 144.0 )( 64742.0 , 118.0 )( 79739.0 , 113.0 )( 94504.0 , 111.0 )( 110910.0 , 102.0 )( 128613.0 , 97.0 )( 158021.0 , 92.0 )( 184392.0 , 85.0 )( 227287.0 , 75.0 )( 297832.0 , 66.0 )( 363196.0 , 60.0 )( 431915.0 , 56.0 )( 519291.0 , 52.0 )};
\addlegendentry{MUAS}  
\end{semilogxaxis}
\end{tikzpicture}}
\caption{Example~\ref{ex:hemker}: Number of iterations and rejections.}
\label{fig:iterations_hemker_different}
\end{figure}

Concerning the efficiency, the situation is similar as in Example~\ref{ex:hmm86}. 
The simulations with the AFC scheme with Kuzmin limiter and the MUAS method needed
generally a similar number of iterations, see Figure~\ref{fig:iterations_hemker_different}. They were often 
considerably more efficient than the simulations with the AFC scheme with BJK limiter.

\subsection{Summary of the Numerical Studies}

Here, the most important findings of the numerical studies are summarized.
\begin{list}{}{\itemsep0.0ex\parsep0.1ex\topsep0.2ex\leftmargin1.6em
\labelwidth1.3em}
\item[$\bullet$] The global DMP was satisfied for all methods on all grids with hanging nodes. 
On grids with conforming closure, it was always satisfied 
for the AFC scheme with BJK limiter and the 
MUAS method. 
\item[$\bullet$] The AFC method with Kuzmin limiter did not always satisfy the DMP on conforming 
grids with locally very small mesh cells, where the discrete problem is locally 
diffusion-dominated. 
\item[$\bullet$] The AFC scheme with BJK limiter and the MUAS method converge if the 
discrete solution becomes (locally) diffusion-dominated, both on adaptive grids with 
conforming closure and with hanging nodes. 
\item[$\bullet$] If the discrete solution becomes (locally) diffusion-dominated, then 
the AFC method with Kuzmin limiter does not convergence on adaptively refined grids
with conforming closure. 
\item[$\bullet$] The nonlinear problems could be solved often most efficiently for the 
AFC scheme with Kuzmin limiter and the MUAS method. 
\end{list}

\section{Summary}\label{sec:summary}
This paper studied the behavior of algebraic stabilizations for discretizing steady-state
convection-diffusion-reaction equations in simulations on adaptively refined
grids, both with conforming closure and with hanging nodes. 
The AFC 
scheme with BJK limiter and the MUAS method satisfied always the global DMP. It could be demonstrated that the failure 
of the AFC method with Kuzmin limiter to satisfy the DMP on some grids with conforming closure could be removed 
by using grids with hanging nodes. The crucial algorithmic component for a successful application of 
algebraically stabilized schemes on grids with hanging nodes is that the linear system of equations is 
transformed to conforming test and conforming ansatz functions for computing the limiters. 
In summary, taking all the aspects of accuracy, satisfaction of the global DMP, and efficiency into 
account, the MUAS method seems to be the most promising of 
%all 
the three
approaches studied in this paper.

\section*{Acknowledgement}
The work of Petr Knobloch has been
supported through the grant No.~20-01074S of the Czech Science Foundation.

\bibliographystyle{alpha}
\bibliography{hanging_nodes_theory_applications}

\begin{thebibliography}{10}

\bibitem{AR10}
{\sc M.~Ainsworth and R.~Rankin}, {\em Fully computable error bounds for
  discontinuous galerkin finite element approximations on meshes with an
  arbitrary number of levels of hanging nodes}, {SIAM} Journal on Numerical
  Analysis, 47 (2010), pp.~4112--4141, \url{https://doi.org/10.1137/080725945}.

\bibitem{ACJ11}
{\sc M.~Augustin, A.~Caiazzo, A.~Fiebach, J.~Fuhrmann, V.~John, A.~Linke, and
  R.~Umla}, {\em An assessment of discretizations for convection-dominated
  convection{\textendash}diffusion equations}, Computer Methods in Applied
  Mechanics and Engineering, 200 (2011), pp.~3395--3409,
  \url{https://doi.org/10.1016/j.cma.2011.08.012}.

\bibitem{BSW83}
{\sc R.~E. Bank, A.~H. Sherman, and A.~Weiser}, {\em Refinement algorithms and
  data structures for regular local mesh refinement}, in Scientific computing
  ({M}ontreal, {Q}ue., 1982), IMACS Trans. Sci. Comput., I, IMACS, New
  Brunswick, NJ, 1983, pp.~3--17.

\bibitem{BJK16}
{\sc G.~R. Barrenechea, V.~John, and P.~Knobloch}, {\em Analysis of algebraic
  flux correction schemes}, {SIAM} Journal on Numerical Analysis, 54 (2016),
  pp.~2427--2451, \url{https://doi.org/10.1137/15m1018216}.

\bibitem{BJK17}
{\sc G.~R. Barrenechea, V.~John, and P.~Knobloch}, {\em An algebraic flux
  correction scheme satisfying the discrete maximum principle and linearity
  preservation on general meshes}, Mathematical Models and Methods in Applied
  Sciences, 27 (2017), pp.~525--548,
  \url{https://doi.org/10.1142/s0218202517500087}.

\bibitem{BK13}
{\sc M.~Bittl and D.~Kuzmin}, {\em An $hp$-adaptive flux-corrected transport
  algorithm for continuous finite elements}, Computing, 95 (2013), pp.~27--48,
  \url{https://doi.org/10.1007/s00607-012-0223-y}.

\bibitem{BS08}
{\sc S.~C. Brenner and L.~R. Scott}, {\em The Mathematical Theory of Finite
  Element Methods}, Springer New York, 2008,
  \url{https://doi.org/10.1007/978-0-387-75934-0}.

\bibitem{Cia78}
{\sc P.~G. Ciarlet}, {\em The finite element method for elliptic problems},
  North-Holland Publishing Co., Amsterdam-New York-Oxford, 1978.
\newblock Studies in Mathematics and its Applications, Vol. 4.

\bibitem{Dav04}
{\sc T.~A. Davis}, {\em Algorithm 832}, {ACM} Transactions on Mathematical
  Software, 30 (2004), pp.~196--199,
  \url{https://doi.org/10.1145/992200.992206}.

\bibitem{GJM+16}
{\sc S.~Ganesan, V.~John, G.~Matthies, R.~Meesala, A.~Shamim, and
  U.~Wilbrandt}, {\em An object oriented parallel finite element scheme for
  computations of {PDEs}: Design and implementation}, in 2016 {IEEE} 23rd
  International Conference on High Performance Computing Workshops ({HiPCW}),
  {IEEE}, Dec. 2016, \url{https://doi.org/10.1109/hipcw.2016.023}.

\bibitem{GT83}
{\sc D.~Gilbarg and N.~S. Trudinger}, {\em Elliptic Partial Differential
  Equations of Second Order}, Springer Berlin Heidelberg, 2001,
  \url{https://doi.org/10.1007/978-3-642-61798-0}.

\bibitem{Car11}
{\sc C.~Gr{\"a}ser}, {\em Convex minimization and phase field models}, PhD
  thesis, Freie Universit{\"a}t, Berlin, 2011.

\bibitem{Hem96}
{\sc P.~Hemker}, {\em A singularly perturbed model problem for numerical
  computation}, Journal of Computational and Applied Mathematics, 76 (1996),
  pp.~277--285, \url{https://doi.org/10.1016/s0377-0427(96)00113-6}.

\bibitem{HMM86}
{\sc T.~J. Hughes, M.~Mallet, and M.~Akira}, {\em A new finite element
  formulation for computational fluid dynamics: {II}. beyond {SUPG}}, Computer
  Methods in Applied Mechanics and Engineering, 54 (1986), pp.~341--355,
  \url{https://doi.org/10.1016/0045-7825(86)90110-6}.

\bibitem{Jha21}
{\sc A.~Jha}, {\em Hanging nodes for higher-order lagrange finite elements},
  Examples and Counterexamples, 1 (2021), p.~100025,
  \url{https://doi.org/10.1016/j.exco.2021.100025}.

\bibitem{Jha20}
{\sc A.~Jha}, {\em A residual based a posteriori error estimators for {AFC}
  schemes for convection-diffusion equations}, Computers {\&} Mathematics with
  Applications, 97 (2021), pp.~86--99,
  \url{https://doi.org/10.1016/j.camwa.2021.05.031}.

\bibitem{JJ19}
{\sc A.~Jha and V.~John}, {\em A study of solvers for nonlinear {AFC}
  discretizations of convection{\textendash}diffusion equations}, Computers
  {\&} Mathematics with Applications, 78 (2019), pp.~3117--3138,
  \url{https://doi.org/10.1016/j.camwa.2019.04.020}.

\bibitem{JJ18}
{\sc A.~Jha and V.~John}, {\em On basic iteration schemes for nonlinear {AFC}
  discretizations}, in Lecture Notes in Computational Science and Engineering,
  Springer International Publishing, 2020, pp.~113--128,
  \url{https://doi.org/10.1007/978-3-030-41800-7_7}.

\bibitem{JK07_1}
{\sc V.~John and P.~Knobloch}, {\em On spurious oscillations at layers
  diminishing ({SOLD}) methods for convection{\textendash}diffusion equations:
  Part i {\textendash} a review}, Computer Methods in Applied Mechanics and
  Engineering, 196 (2007), pp.~2197--2215,
  \url{https://doi.org/10.1016/j.cma.2006.11.013}.

\bibitem{JK21}
{\sc V.~John and P.~Knobloch}, {\em On algebraically stabilized schemes for
  convection-diffusion-reaction problems}, 2021,
  \url{https://arxiv.org/abs/2111.08697}.
\newblock submitted.

\bibitem{JMT96}
{\sc V.~John, J.~Maubach, and L.~Tobiska}, {\em Nonconforming
  streamline-diffusion-finite-element-methods for convection-diffusion
  problems}, Numerische Mathematik, 78 (1997), pp.~165--188,
  \url{https://doi.org/10.1007/s002110050309}.

\bibitem{JN13}
{\sc V.~John and J.~Novo}, {\em A robust {SUPG} norm a posteriori error
  estimator for stationary convection{\textendash}diffusion equations},
  Computer Methods in Applied Mechanics and Engineering, 255 (2013),
  pp.~289--305, \url{https://doi.org/10.1016/j.cma.2012.11.019}.

\bibitem{Kno17}
{\sc P.~Knobloch}, {\em On the discrete maximum principle for algebraic flux
  correction schemes with limiters of upwind type}, in Lecture Notes in
  Computational Science and Engineering, Springer International Publishing,
  2017, pp.~129--139, \url{https://doi.org/10.1007/978-3-319-67202-1_10}.

\bibitem{KR89}
{\sc R.~Kornhuber and R.~Roitzsch}, {\em On adaptive grid refinement in the
  presence of internal or boundary layers.}, Tech. Report SC-89-05, ZIB,
  Takustr. 7, 14195 Berlin, 1989.

\bibitem{Ku07}
{\sc D.~Kuzmin}, {\em Algebraic flux correction for finite element
  discretizations of coupled systems}, in Proceedings of the Int.~Conf.~on
  Computational Methods for Coupled Problems in Science and Engineering,
  M.~Papadrakakis, E.~O{\~n}ate, and B.~Schrefler, eds., CIMNE, Barcelona,
  2007, pp.~1--5.

\bibitem{Riv84}
{\sc M.-C. Rivara}, {\em Mesh refinement processes based on the generalized
  bisection of simplices}, {SIAM} Journal on Numerical Analysis, 21 (1984),
  pp.~604--613, \url{https://doi.org/10.1137/0721042}.

\bibitem{RST08}
{\sc H.~G. Roos, M.~Stynes, and L.~Tobiska}, {\em Robust numerical methods for
  singularly perturbed differential equations}, vol.~24 of Springer Series in
  Computational Mathematics, Springer-Verlag, Berlin, second~ed., 2008.
\newblock Convection-diffusion-reaction and flow problems.

\bibitem{WB16}
{\sc U.~Wilbrandt, C.~Bartsch, N.~Ahmed, N.~Alia, F.~Anker, L.~Blank,
  A.~Caiazzo, S.~Ganesan, S.~Giere, G.~Matthies, R.~Meesala, A.~Shamim,
  J.~Venkatesan, and V.~John}, {\em {ParMooN}{\textemdash}a modernized program
  package based on mapped finite elements}, Computers {\&} Mathematics with
  Applications, 74 (2017), pp.~74--88,
  \url{https://doi.org/10.1016/j.camwa.2016.12.020}.

\end{thebibliography}
\end{document}